\documentclass[11pt, a4paper]{article}%
\usepackage{amsmath}
\usepackage{amsfonts,xcolor}
\usepackage{amssymb}
\usepackage{graphicx,sw20jart}
\usepackage{caption}
\usepackage{subcaption}
\usepackage{graphicx}%
\providecommand{\U}[1]{\protect\rule{.1in}{.1in}}
\newtheorem{theorem}{Theorem}[section]

\newtheorem{definition}{Definition}[section]

\newtheorem{lemma}[theorem]{Lemma}
\newtheorem{proposition}[theorem]{Proposition}
\newtheorem{remark}{Remark}[section]

\begin{document}

\title{Optimal dividends under a drawdown constraint\\and a curious square-root rule}
\author{Hansj\"{o}rg Albrecher\thanks{Department of Actuarial Science, Faculty of
Business and Economics, University of Lausanne, CH-1015 Lausanne and Swiss
Finance Institute. Supported by the Swiss National Science Foundation Project
200021\_191984.}, Pablo Azcue\thanks{Departamento de Matematicas, Universidad
Torcuato Di Tella. Av. Figueroa Alcorta 7350 (C1428BIJ) Ciudad de Buenos
Aires, Argentina.} and Nora Muler$^{\dag}$}
\date{}
\maketitle

\bigskip\abstract{\begin{quote}
\noindent In this paper we address the problem of optimal dividend payout strategies from a surplus process governed by Brownian
motion with drift under a drawdown constraint, i.e.\ the dividend rate can
never decrease below a given fraction $a$ of its historical maximum. We solve the resulting two-dimensional optimal control problem and identify the
value function as the unique viscosity solution of the corresponding Hamilton-Jacobi-Bellman
equation. We then derive sufficient conditions under which a two-curve strategy is optimal, and show how to determine its concrete form using calculus of variations.
{We establish a smooth-pasting principle and show how it can be used to prove the optimality of two-curve strategies for sufficiently large initial and maximum dividend rate}. We also give {a number of} numerical illustrations in which the optimality of the two-curve strategy can be established for instances with smaller values of the maximum dividend rate, and the concrete form of the curves can be determined. One observes that the resulting drawdown strategies nicely interpolate between the solution for the classical unconstrained dividend problem and the one for a ratcheting constraint as recently studied in \cite{AAM21}. When the maximum allowed dividend rate tends to infinity, we show a surprisingly simple and somewhat intriguing limit result in terms of the parameter $a$ for the surplus level on from which, for sufficiently large current dividend rate, a take-the-money-and-run strategy is optimal in the presence of the drawdown constraint.
\end{quote}}

\section{Introduction and model}

Assume that the surplus process of a company is given by a Brownian motion
with drift%
\begin{equation}
X_{t}=x+\mu t+\sigma W_{t}, \label{model}%
\end{equation}
where $W_{t}$ is a standard Brownian motion, and $\mu,\sigma>0$ are given
constants. Let $(\Omega,\mathcal{F},\left(  \mathcal{F}_{t}\right)  _{t\geq
0},\mathcal{P})$ be the complete probability space generated by the process
$X_{t}$. Assume further that the company uses part of the surplus to pay
dividends to the shareholders with rates in a set $[0,\overline{c}]$, where
$\overline{c}>0$ is the maximum dividend rate possible. Let $D_{t}$ denote the
rate at which the company pays dividends at time $t$, then the controlled
surplus process can be written as%
\begin{equation}
X_{t}^{D}=X_{t}-\int_{0}^{t}D_{s}ds. \label{XtC}%
\end{equation}
It is a classical problem in risk theory to find the dividend strategy
$D=\left(  D_{t}\right)  _{t\geq0}$ that maximizes the expected sum of
discounted dividend payments
\begin{equation}
J(x;D)=\mathbb{E}\left[  \int_{0}^{\tau}e^{-qs}D_{s}ds\right]  \label{JxD}%
\end{equation}
over a set of admissible candidate strategies. Here $q>0$ is a discount factor
and $\tau=\inf\left\{  t\geq0:X_{t}^{D}<0\right\}  $ is the ruin time of the
controlled process. De Finetti \cite{defin} was the first to consider a
problem of this kind for a simple random walk, and Gerber \cite{Ger69,Ger72}
considered extensions, including the diffusion setup \eqref{model} given
above. For a finite maximum dividend rate $\overline{c}$, this problem was
then further investigated by Shreve et al.\ \cite{shreve}, Jeanblanc and
Shiryaev \cite{Jeanblanc}, Radner and Shepp \cite{Radner}, Asmussen and Taksar
\cite{asmtak} and Gerber and Shiu \cite{Gerber2004}. Since then, a lot of
variants of this problem for the process \eqref{model} and more general
underlying risk processes have been considered, see e.g.\ the surveys
\cite{AT} and \cite{avanzi}.\newline For the diffusion model \eqref{model}, in
\cite{AAM21} we recently studied this optimal dividend problem under a
ratcheting constraint, i.e.\ under the assumption that the dividend rate can
never be decreased over the lifetime of the process, which renders the
respective control problem two-dimensional, where the first dimension is the
current surplus and the second dimension is the currently employed dividend
rate. One motivation to consider that constraint was that it may be
psychologically preferable to shareholders to not experience a decrease of
dividend payments, and it is interesting to see to what extent such a
constraint leads to an overall performance loss. \newline

In this paper we would like to go one step further and allow reductions of the
dividend rate over time, but only up to a certain percentage $a$ of the
largest already exercised dividend rate (''drawdown''). More formally, a
\textit{dividend drawdown strategy} $D=\left(  D_{t}\right)  _{t\geq0}$ with
drawdown constraint $a\in\lbrack0,1]$ is one that satisfies $D_{t}\in\lbrack
aR_{t},\overline{c}]$, where $R_{t}$ is the running maximum of the dividend
rates, that is%
\[
R_{t}:=\max\{D_{s}:0\leq s\leq t\}\vee c;
\]
here we denote the initial dividend rate by $R_{0^{-}}=c$. A dividend drawdown
strategy is called \textit{admissible }if it is right-continuous and adapted
with respect to the filtration $\left(  \mathcal{F}_{t}\right)  _{t\geq0}%
.$\newline

Define $\Pi_{x,c,a}^{[0,\overline{c}]}$ as the set of all admissible dividend
drawdown strategies with initial surplus $x\geq0$, initial running maximum
dividend rate $c\in\lbrack0,\overline{c}]$ and drawdown constraint
$a\in\lbrack0,1]$. Given $D\in\Pi_{x,c,a}^{[0,\overline{c}]}$, the value
function of this strategy is given by \eqref{JxD}. Hence, for any initial
surplus $x\geq0$, initial running maximum dividend rate $c\in\lbrack
0,\overline{c}]$ and drawdown constraint $a\in\lbrack0,1]$, our aim in this
paper is to maximize%
\begin{equation}
V_{a}^{\overline{c}}(x,c)=\sup_{D\in\Pi_{x,c,a}^{[0,\overline{c}]}}J(x;D).
\label{Optimal Value Function}%
\end{equation}
Note that the limit case $a=1$ corresponds to the ratcheting case (considered
previously in {\cite{AAM21}}) and the limit case $a=0$ corresponds to the
optimization of bounded dividend rates without any drawdown
constraint.\newline

Drawdown phenomena have been studied in various contexts in the literature. On
the one hand, drawdown times and properties of uncontrolled stochastic
processes were investigated in quite some generality (see for instance
Landriault et al.\ \cite{Land} for the case of L\'evy processes). In the
context of control problems, drawdown constraints on the wealth have been
considered in portfolio problems in the mathematical finance literature, see
for instance Elie and Touzi \cite{elie}, Chen et al.\ \cite{chen} and Kardaras
et al.\ \cite{Kar}. For a minimization of drawdown times of a risk process
through dynamic reinsurance, see Brinker \cite{Brinker} and Brinker and
Schmidli \cite{BrinkerS}. Our context, however, is different, as we are
interested in implementing a drawdown constraint on the payment structure of
the dividend rates, i.e., as a constraint on the admissible dividend policies.
In that sense, our approach is closer related to problems of lifetime
consumption in the mathematical finance literature, see Angoshtari et
al.\ \cite{bayr} who extend the Dusenberry's ratcheting problem of consumption
studied by Dybvig \cite{dybvig} to drawdown constraints. However, the concrete
model setup and embedding, and also the involved techniques there are very
different from dividend problems of the De Finetti-type as studied in this
paper. \newline

After deriving some basic analytic properties of the value function
$V_{a}^{\overline{c}}(x,c)$ of our drawdown problem in Section \ref{sec3}, we
will derive a Hamilton-Jacobi-Bellman equation for $V_{a}^{\overline{c}}(x,c)$
in Section \ref{Hamilton-Jacobi-Bellman equations} and show that
$V_{a}^{\overline{c}}(x,c)$ is its unique viscosity solution with suitable
boundary condition. We then, in Section \ref{sec5}, briefly study in more
detail the value function when one already starts at the maximal dividend rate
$\overline{c}$, which serves as a crucial ingredient in the derivation of
$V_{a}^{\overline{c}}(x,c)$ in Section
\ref{Curve strategies and optimal curve}. Sufficient conditions are given
under which the optimal strategy for bounded dividend rates is a two-curve
strategy in the space $(x,c)\in(0,\infty)\times\lbrack0,\overline{c}]$, which
is partitioned by two curves $\gamma^{\overline{c}}(c)$ and $\zeta
^{\overline{c}}(c)$, with $\gamma^{\overline{c}}(c)<\zeta^{\overline{c}}(c)$
for all $c\in\lbrack0,\overline{c}]$: if for a given $c$, $x<\gamma
^{\overline{c}}(c)$, then dividends are paid at rate $ac$; if $\gamma
^{\overline{c}}(c)\leq x\leq\zeta^{\overline{c}}(c)$, then dividends are paid
at rate $c$; finally, if $x>\zeta^{\overline{c}}(c)$, then the dividend rate
$c$ is increased immediately until $x=\zeta^{\overline{c}}(c_1)$ for some $c_1\in(c,\overline{c})$ {(or $c=\overline{c}$, whichever happens first)} is reached. {We furthermore establish a smooth-pasting principle for these optimal curves. In
Section \ref{Asymptotic Values} it is shown that the limits of $\gamma
^{\overline{c}}(\overline{c})$ and $\zeta^{\overline{c}}(\overline{c})$ as
$\overline{c}\rightarrow\infty$ are finite, and given by surprisingly explicit
formulas:
\begin{equation}
\lim_{\overline{c}\rightarrow\infty}\gamma^{\overline{c}}(\overline{c}%
)=\frac{\mu}{q}\quad\text{and}\quad\lim_{\overline{c}\rightarrow\infty}%
\zeta^{\overline{c}}(\overline{c})=\frac{\mu}{q}\Big(1+\frac{1}{\sqrt{a}%
}\Big).
\label{limitta}%
\end{equation}
This nicely extends the respective limit $2\mu
/q$ of the ratcheting curve that was identified for pure ratcheting ($a=1$) in \cite[Lem.5.21]{AAM21}.\\ In Section \ref{Optimal strategies for c large enough} we then look further into the limiting case, and show that for sufficiently large $\overline{c}$, one has  $\gamma^{\overline{c}}(c)\nearrow\gamma^{\overline{c}}(\overline{c})$ and  $\zeta^{\overline{c}}(c)\searrow\zeta^{\overline{c}}(\overline{c})$ as $c\to \overline{c}$. This enables to establish the general optimality of two-curve strategies whenever the current dividend rate $c$ and the maximal dividend rate $\overline{c}$ are sufficiently large. At the same time, the negative derivative of $\zeta^{\overline{c}}(c)$ close to (sufficiently large) $\overline{c}$ is {notably different from the pure ratcheting case ($a=1$) for which it was shown in \cite{AAM21} that the corresponding derivative is positive for all $c$ close to $\overline{c}$ (and indeed the leading term in the asymptotics of $0<a<1$ breaks down for $a=1$ so that some sort of phase transition happens)}. The simplicity of the right-hand limit in \eqref{limitta} and in particular the appearance of
the square-root of the drawdown coefficient $a$ in the right-hand limit are somewhat intriguing. In the absence of an upper limit for the dividend rate, it identifies the minimum surplus level $x$ on from which, for sufficiently large current dividend rate, it is preferable to pay out all the
surplus $x$ immediately and generate ruin by doing so (a so-called
\textquotedblright take-the-money-and-run\textquotedblright-strategy, see
e.g.\ \cite{LoeffenRenaud}), and that value does not depend on the size of the volatility $\sigma$. Consequently, one can get some intuition on its nature in the much simpler
deterministic model with $\sigma=0$, which we will therefore consider in Section \ref{sec2} before approaching
the general case $\sigma>0$ in the rest of the paper.\\
We give numerical  illustrations in Section
\ref{Numerical examples}, where we establish the optimality of two-curve strategies also for smaller magnitudes of $c$ and $\overline{c}$ by numerically showing that the sufficient conditions from Section
\ref{Curve strategies and optimal curve} are satisfied. We obtain the optimal curves by calculus of variation techniques and discuss the properties of the value functions of the drawdown dividend problem and their comparison to classical and ratcheting solutions for various parameter combinations. Finally, Section \ref{secconcl} concludes, and Section
\ref{Seccion Formulas} collects some longer formulas appearing in the paper in
a compact form.}

\section{Some intuition from the deterministic case}

\label{sec2} Assume in this section for simplicity a completely deterministic
model
\[
X_{t}=x+\mu t
\]
with a positive drift $\mu>0$ (for the study of such a model in another
context in the dividend literature, see e.g.\ \cite{thon}). Then a constant
dividend rate $\overline{c}=\mu$ throughout time will keep the surplus at
level $x$ for all $t\geq0$ and correspondingly
\[
\mathbb{E}\left[  \int_{0}^{\tau}e^{-qs}D_{s}ds\right]  =\mathbb{E}\left[
\int_{0}^{\infty}e^{-qs}\mu\,ds\right]  =\frac{\mu}{q}%
\]
for any $x>0$. Consequently, whenever the initial surplus $x$ is larger than
$\mu/\delta$, paying out all the surplus at the beginning (causing immediate
ruin) will be preferable to any other dividend strategy subject to the
constraint $\overline{c}\leq\mu$. \newline At the same time, if a constant
dividend rate $\overline{c}>\mu$ is applied, the controlled process will lead
to ruin at time $t=x/(\overline{c}-\mu)$ and we obtain instead
\begin{align*}
\mathbb{E}\left[  \int_{0}^{\tau}e^{-qs}D_{s}ds\right]   &  =\mathbb{E}\left[
\int_{0}^{x/(\overline{c}-\mu)}e^{-qs}\overline{c}\,ds\right]  =\frac
{\overline{c}}{q}\left(  1-e^{-qx/(\overline{c}-\mu)}\right) \\
&  =x+x\,\frac{2\mu-qx}{2\overline{c}}+O\left(  \frac{1}{\overline{c}^{2}%
}\right)  .
\end{align*}
The latter shows that whenever $x>2\mu/q$, if allowed to do so, paying out all
the surplus $x$ immediately (and causing immediate ruin) will be preferable to
any other constant dividend strategy with large $\overline{c}$. In other
words, the potential gain from later ruin and therefore more dividend income
(by exploiting the positive drift, without any risk) is outweighed by the
discounting of such later dividend payments. This can also be seen as an
intuitive explanation of the limit $2\mu/q$ in \cite[Lem.5.21]{AAM21}.
\newline

Let us now proceed to the case with drawdown: Assume that we start with
initial capital $x>b$ for some $b$ to be determined and that we pay dividends
at rate $\overline{c}>\mu$ until we reach that lower level $b$ at time
$t=(x-b)/(\overline{c}-\mu)$, from which time on we reduce the dividend
payments to level $a\cdot\overline{c}$ according to our drawdown constraint.
In the deterministic model of this section, this then leads to
\begin{align}
\mathbb{E}\left[  \int_{0}^{\tau}e^{-qs}D_{s}ds\right]   &  =\frac
{\overline{c}}{q}\left(  1-e^{-q \frac{x-b}{\overline{c}-\mu}}\right)  +e^{-q
\frac{x-b}{\overline{c}-\mu}}\cdot\frac{a\,\overline{c}}{q}\,\left(  1-e^{-q
\frac{b}{a\overline{c}-\mu}}\right)  . \label{dadrin}%
\end{align}
Taking the derivative with respect to $b$ and setting it zero gives, after
simple calculations, for large $\overline{c}$, the optimal level
\begin{align}
b^{*}(\overline{c})  &  =\frac{a \overline{c}-\mu}{q}\log\left(  \frac{a
\overline{c}}{a \overline{c}-\mu}\right) \label{still}\\
&  =\frac{\mu}{q}-\frac{\mu^{2}}{2aq\overline{c}}+O\left(  \frac{1}%
{\overline{c}^{2}}\right)  .\nonumber
\end{align}
But if one substitutes that value of $b$ into \eqref{dadrin}, then an
expansion at $\overline{c}=\infty$ gives
\begin{align}
\mathbb{E}\left[  \int_{0}^{\tau}e^{-qs}D_{s}ds\right]   &  =x+\frac
{2axq\mu-ax^{2}q^{2}+\mu^{2}(1-a)}{2 a q \,\overline{c}}+O\left(  \frac
{1}{\overline{c}^{2}}\right)  . \label{lowerorder}%
\end{align}
The numerator in the second term is negative exactly when
\[
x>\frac{\mu}{q}\Big(1+\frac{1}{\sqrt{a}}\Big),
\]
so that in those cases it is preferable to immediately pay $x$ as a lump sum
dividend and go to ruin immediately (if that is allowed) rather than following
the above refracting strategy, as the value $x$ can not be realized at any
later point in time in view of the discounting, despite the continuing
deterministic income with drift $\mu$.\newline One may expect that the size of
the volatility does not matter when $\overline{c}\to\infty$, and indeed, as a
by-product of the results of this paper, it will be shown in Section
\ref{Asymptotic Values} that the same limit can be established for the general
case $\sigma>0$, cf.\ Proposition \ref{Derivada de V cbarra}. Another way to
state this is the following: if one defines $x^{*}(\overline{c})$ as the
surplus value for which, when already currently paying the maximum dividend
rate $\overline{c}$, one is indifferent whether to further increase the
dividend rate or not, then the above result establishes that $\lim
_{\overline{c}\to\infty}x^{*}(\overline{c})=\frac{\mu}{q}(1+{1}/\sqrt{a})$,
and it will be in terms of that notation that the more general result is
proved in Section \ref{Asymptotic Values}.

\section{Basic results \label{Model and basic results}}

\label{sec3} Recall the definition of our optimal value function
$V_{a}^{\overline{c}}(x,c)$ \eqref{Optimal Value Function} and denote by
$V_{a}^{\infty}(x,c)$ the corresponding function when there is no ceiling on
dividend rates, i.e.\ $\overline{c}=\infty$. It is immediate to see that
$V_{a}^{\overline{c}}(0,c)=0$ for all $c\in\lbrack0,\overline{c}]$\textbf{
}and $a\in\lbrack0,1]$.

\begin{remark}
\normalfont\label{Optima sin drawdown acotada} As mentioned in the
introduction, the dividend optimization problem without drawdown constraint
has a long history and, for a finite $\overline{c}$ and the diffusion setup,
was first addressed in Shreve et al.\ \cite{shreve}. Unlike the drawdown
optimization problem, the problem without the drawdown constraint is
one-dimensional. If we denote its optimal value function by $\overline
{V}^{\overline{c}}(x),$ then clearly $V_{0}^{\overline{c}}(x,c)=\overline
{V}^{\overline{c}}(x)$ and $V_{a}^{\overline{c}}(x,c)\leq\overline
{V}^{\overline{c}}(x)$ for all $x\geq0$, $a\in\lbrack0,1]$ and $c\in
\lbrack0,\overline{c}]$. The function $\overline{V}^{\overline{c}}$ is
increasing, concave, twice continuously differentiable with $\overline
{V}^{\overline{c}}(0)=0$ and $\lim_{x\rightarrow\infty}\overline{V}%
^{\overline{c}}(x)=\overline{c}/q$; so it is Lipschitz with Lipschitz constant
$\left(  \overline{V}^{\overline{c}}\right)  ^{\prime}(0).$
\end{remark}

\begin{remark}
\normalfont\label{Optima sin drawdown irrestricta} The dividend optimization
problem without any constraint was addressed by {Gerber and Shiu
\cite{Gerber2004} and Schmidli \cite{Schmidli book 2008}}. If $\overline
{V}(x)$ denotes its optimal value function, we have $\overline{V}%
(x)=V_{0}^{\infty}(x,c)$ for any $c>0$. Clearly $V_{a}^{\infty}(x,c)\leq
\overline{V}(x)$ for all $a\in\lbrack0,1]$. The function $\overline{V}$ is
increasing, concave, twice continuously differentiable with $\overline
{V}(0)=0$ and $x\leq\overline{V}(x)\leq x+\mu/q$; so it is Lipschitz with
Lipschitz constant $\overline{V}^{\prime}(0)$.
\end{remark}

\begin{proposition}
\label{Limite de restricted 1} It holds that $V_{a}^{\overline{c}%
}(x,c)\nearrow V_{a}^{\infty}(x,c)$ as $\overline{c}\rightarrow\infty.$
\end{proposition}

\textit{Proof.} It is straightforward that for any $\overline{c}_{1}\leq$
$\overline{c}_{2}$, $V_{a}^{\overline{c}_{1}}(x,c)\leq V_{a}^{\overline{c}%
_{2}}(x,c)$ $\leq V_{a}^{\infty}(x,c)$ for $0\leq c\leq\overline{c}_{1}$. Take
for any $\varepsilon>0,$ a strategy $D=(D_{t})_{t\geq0}\in\Pi_{x,c,a}%
^{[0,\infty)}$ with ruin time $\tau$, such that $V_{a}^{\infty}(x,c)\leq
J(x;D)+\varepsilon$. Let us consider for an increasing sequence with
$c_{n}\rightarrow\infty$ and $c_{1}>c$, $D^{n}=(D_{t}\wedge c_{n})_{t\geq0}%
\in\Pi_{x,c,a}^{[0,c_{n})}$ and let $\tau^{n}\geq\tau$ be the ruin time of
$D^{n}$. Then, by the theorem of monotone convergence,%

\[
\lim_{n\rightarrow\infty}J(x;D_{t}^{n})=\lim_{n\rightarrow\infty}%
\mathbb{E}\left[  \int_{0}^{\tau^{n}}e^{-qs}D_{s}^{n}ds\right]  \geq
\lim_{n\rightarrow\infty}\mathbb{E}\left[  \int_{0}^{\tau}e^{-qs}D_{s}%
^{n}ds\right]  =J(x;D)
\]
and so we have the result.\hfill$\blacksquare$\newline

We now state a straightforward result regarding the boundedness and
monotonicity of the optimal value functions.

\begin{proposition}
\label{Monotone Optimal Value Function}In the case $\overline{c}<\infty$, the
optimal value function $V_{a}^{\overline{c}}(x,c)$ is bounded by $\overline
{c}/q$ with $\lim_{x\rightarrow\infty}V_{a}^{\overline{c}}(x)=\overline{c}/q$,
non-decreasing in $x$ and non-increasing in $c$.
\end{proposition}

\textit{Proof.} By Remark \ref{Optima sin drawdown acotada} and Theorem 3.3 of
{\cite{AAM21}, we have that}
\[
V_{1}^{\overline{c}}(x,c)\leq V_{a}^{\overline{c}}(x,c)\leq\overline
{V}^{\overline{c}}(x)
\]
with $\lim_{x\rightarrow\infty}V_{1}^{\overline{c}}(x,c)=\lim_{x\rightarrow
\infty}\overline{V}^{\overline{c}}(x)=\overline{c}/q,$ so $V_{a}^{\overline
{c}}$ is bounded by $\overline{c}/q$ with $\lim_{x\rightarrow\infty}%
V_{a}^{\overline{c}}(x,c)=\overline{c}/q.$

On the one hand $V_{a}^{\overline{c}}(x,c)$ is non-increasing in $c$ because
given $c_{1}<c_{2}\leq\overline{c},$ we have $\Pi_{x,c_{2},a}^{[0,\overline
{c}]}\subset\Pi_{x,c_{1,a}}^{[0,\overline{c}]}$ for any $x\geq0$. On the other
hand, given $0\leq x_{1}<x_{2}$ and an admissible ratcheting strategy
$D_{1}\in\Pi_{x_{1},c,a}^{[0,\overline{c}]}$ for any $c\in\lbrack
0,\overline{c}]$, let us define $D_{2}\in\Pi_{x_{2},c,a}^{[0,\overline{c}]}$
as $D_{2,t}=D_{1.t}$ until the ruin time of the controlled process
$X_{t}^{D^{1}}$ with $X_{0}^{D^{1}}=x_{1}$, and pay the maximum rate
$\overline{c}$ afterwards. Thus, $J(x;D_{1})\leq J(x;D_{2})$ and we have the
result. \hfill$\blacksquare$\newline

\begin{proposition}
\label{Limite de restricted 2} $V_{a}^{\infty}(x,c)$ is non-decreasing in $x$
and non-increasing in $c.$ For the case $a>0,\ $we have $\lim_{c\rightarrow
\infty}V_{a}^{\infty}(x,c)=x$. Moreover $x\leq V_{a}^{\infty}(x,c)\leq
x+\mu/q.$
\end{proposition}

\textit{Proof.} By Propositions \ref{Limite de restricted 1} and
\ref{Monotone Optimal Value Function}, we have that $V_{a}^{\infty}(x,c)$ is
non-decreasing in $x$ and non-increasing in $c$. Let us show now that
$V_{a}^{\infty}(x,c)\geq x$. The function $V_{a}^{\infty}(x,c)$ is bounded
from below by the expected discounted dividends resulting from the strategy of
paying a constant rate $n$ up to ruin. Defining ${\tau}_{n}=\inf
\{t:x+(\mu-n)t+\sigma W_{t}=0\}$, one gets
\[
V_{a}^{\infty}(x,c)=\lim_{n\rightarrow\infty}V_{a}^{n}(x,c)\geq\lim
_{n\rightarrow\infty}\mathbb{E}\left[  \int_{0}^{{\tau}_{n}}e^{-qs}%
\,n\,ds\right]  =\lim_{n\rightarrow\infty}\frac{n}{q}(1-\mathbb{E[}e^{-q{\tau
}_{n}}])=x,
\]
where the last equality follows from Formula 2.0.1 on page 295 of Borodin \&
Salminen \cite{BorodinSalminem}.

Finally, let us see that $\lim_{c\rightarrow\infty}V_{a}^{\infty}(x,c)\leq x$.
Take for any $\varepsilon>0$ and for each $c$, $D^{c}=(D_{t}^{c})_{t\geq0}%
\in\Pi_{x,c,a}^{[0,\infty)}$, such that
\[
V_{a}^{\infty}(x,c)\leq J(x;D^{c})+\varepsilon.
\]
The corresponding ruin time is then given by%
\[
\tau^{c}=\inf\left\{  t:x+\mu t+\sigma W_{t}-%
{\textstyle\int\nolimits_{0}^{t}}
D_{s}^{c}ds=0\right\}
\]
and $D_{s}^{c}\geq ac.$ Hence,%

\[%
{\textstyle\int\nolimits_{0}^{\tau^{c}}}
D_{s}^{c}ds=x+\mu\tau^{c}+\sigma W_{\tau^{c}}%
\]
and so
\[
\tau^{c}\leq\inf\left\{  s:x+(\mu-ac)s+\sigma W_{s}\leq0\right\}
=\inf\left\{  s:W_{s}\leq\frac{-x+(ac-\mu)s}{\sigma}\right\}  \text{.}%
\]
Hence, for $c>\mu/a$, $\tau^{c}<\infty$ a.s. and $\mathbb{E}\left[  \tau
^{c}\right]  \rightarrow0$ as $c\rightarrow\infty$. Therefore,%

\[%
\begin{array}
[c]{ll}%
\lim_{c\rightarrow\infty}\mathbb{E}[%
{\textstyle\int\nolimits_{0}^{\tau^{c}}}
e^{-qs}D_{s}^{c}ds] & \leq\lim_{c\rightarrow\infty}\mathbb{E}[%
{\textstyle\int\nolimits_{0}^{\tau^{c}}}
D_{s}^{c}ds]\\
& =\lim_{c\rightarrow\infty}\mathbb{E}[x+\mu\tau^{c}+\sigma W_{\tau^{c}}]\\
& =x+\mu\lim_{c\rightarrow\infty}\mathbb{E}[\tau^{c}]=x
\end{array}
\]
and so we have the result. \hfill$\blacksquare$\newline

The Lipschitz property of the function $\overline{V}$ can now be used to prove
a global Lipschitz result on the regularity of the optimal value function.

\begin{proposition}
\label{Proposition Global Lipschitz zone}In both the restricted case
$\overline{c}<\infty$ and the unrestricted case $\overline{c}=\infty,$ we have
that
\[
0\leq V_{a}^{\overline{c}}(x_{2},c_{1})-V_{a}^{\overline{c}}(x_{1},c_{2})\leq
K\left[  \left(  x_{2}-x_{1}\right)  +\left(  c_{2}-c_{1}\right)  \right]
\]

for all $0\leq x_{1}\leq x_{2}$ and $c_{1},c_{2}\in\lbrack0,\overline{c}]$
with $c_{1}\leq c_{2}$, with $K=\max\{\frac{e^{-1}}{q}a,1\}\overline
{V}^{\prime}(0).$
\end{proposition}

\textit{Proof}. In the case $\overline{c}<\infty$, by Proposition
\ref{Monotone Optimal Value Function}, we have
\begin{equation}
0\leq V_{a}^{\overline{c}}(x_{2},c_{1})-V_{a}^{\overline{c}}(x_{1},c_{2})
\label{Lips1}%
\end{equation}
for all $0\leq x_{1}\leq x_{2}$ and $c_{1},c_{2}\in\lbrack0,\overline{c}]$
with $c_{1}\leq c_{2}$.

Let us show now, that there exists $K_{1}>0$ such that%
\begin{equation}
V_{a}^{\overline{c}}(x_{2},c)-V_{a}^{\overline{c}}(x_{1},c)\leq K_{1}\left(
x_{2}-x_{1}\right)  \label{Lips2}%
\end{equation}
for all $0\leq x_{1}\leq x_{2}$. Take $\varepsilon>0$ and $D\in\Pi_{x_{2}%
,c,a}^{[0,\overline{c}]}$ such that%

\begin{equation}
J(x_{2};D)\geq V_{a}^{\overline{c}}(x_{2},c)-\varepsilon,
\label{casi optima x2}%
\end{equation}
the associated control process is given by%

\[
X_{t}^{D}=x_{2}+\int_{0}^{t}(\mu-D_{s})ds+\sigma W_{t}.
\]
Let $\tau$ be the ruin time of the process $X_{t}^{D}$. Define $\widetilde
{D}\in\Pi_{x_{1},c,a}^{[0,\overline{c}]}$ as $\widetilde{D}_{t}=D_{t}$ and the
associated control process
\[
\text{ }X_{t}^{\widetilde{D}}=x_{1}+\int_{0}^{t}(\mu-D_{s})ds+\sigma W_{t}.
\]
Let $\widetilde{\tau}\leq\tau$ be the ruin time of the process $X_{t}%
^{\widetilde{D}}$; it holds that $X_{t}^{D}-X_{t}^{\widetilde{D}}=x_{2}-x_{1}$
for $t\leq\widetilde{\tau}$.

We can write

{
\begin{equation}%
\begin{array}
[c]{lll}%
J(x_{2};D)-J(x_{1};\widetilde{D}) & = & \mathbb{E}\left[  \int_{\widetilde
{\tau}}^{\tau}e^{-qs}D_{s}ds\right] \\
& = & \mathbb{E}\left[  \mathbb{E}\left[  \left.  \int_{\widetilde{\tau}%
}^{\tau}e^{-qs}D_{s}ds\right\vert \mathcal{F}_{\widetilde{\tau}}\right]
\right] \\
& = & \mathbb{E}\left[  \mathbb{E}\left[  \left.  e^{-q\widetilde{\tau}}%
\int_{0}^{\tau-\widetilde{\tau}}e^{-qu}D_{\widetilde{\tau}+u}du\right\vert
\mathcal{F}_{\widetilde{\tau}}\right]  \right] \\
& \leq & \mathbb{E}\left[  \mathbb{E}\left[  \left.  \int_{0}^{\tau
-\widetilde{\tau}}e^{-qu}D_{\widetilde{\tau}+u}du\right\vert \mathcal{F}%
_{\widetilde{\tau}}\right]  \right]  \ \\
& \leq & V_{a}^{\overline{c}}(x_{2}-x_{1},0).
\end{array}
\label{Nueva desigualdad}%
\end{equation}
}

{ The last inequality of (\ref{Nueva desigualdad}) involves a shift of
stopping times and follows from Theorem 2 of Claisse, Talay and Tan
\cite{Claisse}. Indeed, the assumptions of this theorem are satisfied, because
we can write our controlled process as}%

\[
dX_{s}=b(s,X,D_{s})ds+\sigma(s,X,D_{s})dW_{s},
\]
where $b(s,x,d)=\mu-d$, $\sigma(s,x,d)\equiv\sigma$ and $W_{s}$ is a standard
Brownian motion. Hence we have
\begin{equation}%
\begin{array}
[c]{lll}%
V_{a}^{\overline{c}}(x_{2},c)-V_{a}^{\overline{c}}(x_{1},c) & \leq &
J(x_{2};D)-J(x_{1};\widetilde{D})+\varepsilon\\
& \leq & V^{\overline{c}}(x_{2}-x_{1},0)+\varepsilon\\
& \leq & \overline{V}(x_{2}-x_{1})+\varepsilon\\
& \leq & K_{1}(x_{2}-x_{1})+\varepsilon.
\end{array}
\label{diferencia de vs}%
\end{equation}
So, by Remark \ref{Optima sin drawdown irrestricta}, we have (\ref{Lips2})
with $K_{1}=\overline{V}^{\prime}(0).$

Let us show now that, given $c_{1},c_{2}\in\lbrack0,\overline{c}]$ with
$c_{1}\leq c_{2},$ there exists $K_{2}>0$ such that%
\begin{equation}
V_{a}^{\overline{c}}(x,c_{1})-V_{a}^{\overline{c}}(x,c_{2})\leq K_{2}\left(
c_{2}-c_{1}\right)  . \label{Lips3}%
\end{equation}
Take $\varepsilon>0$ and $D\in\Pi_{x,c_{1,a}}^{[0,\overline{c}]}$ such that%

\begin{equation}
J(x;D)\geq V_{a}^{\overline{c}}(x,c_{1})-\varepsilon\label{casiOptimac1}%
\end{equation}
and denote by $\tau$ the ruin time of the process $X_{t}^{D}$.

Let us consider $\widetilde{D}\in\Pi_{x,c_{2}}^{[0,\overline{c}]}\ $as
$\widetilde{D}_{t}=\max\{D_{t},ac_{2}\}$; denote by $X_{t}^{\widetilde{D}}$
the associated controlled surplus process and by $\overline{\tau}\leq\tau$ the
corresponding ruin time. We have that $\widetilde{D}_{s}-D_{s}\leq
ac_{2}-ac_{1}$ and so $X_{\overline{\tau}}^{D}=X_{\overline{\tau}}%
^{D}-X_{\overline{\tau}}^{\widetilde{D}}\leq a(c_{2}-c_{1})\overline{\tau}$.
By Remark \ref{Optima sin drawdown irrestricta}, we have{
\[%
\begin{array}
[c]{lll}%
\mathbb{E}\left[  \int_{\overline{\tau}}^{\tau}D_{s}e^{-qs}ds\right]  & = &
\mathbb{E}\left[  \mathbb{E}\left[  e^{-q\overline{\tau}}\left.
\int_{\overline{\tau}}^{\tau}D_{s}e^{-q\left(  s-\overline{\tau}\right)
}ds\right\vert \mathcal{F}_{\overline{\tau}}\right]  \right] \\
& \leq & \mathbb{E}\left[  \mathbb{E}\left[  \left.  \int_{0}^{\tau
-\overline{\tau}}D_{u+\overline{\tau}}e^{-qu}du\right\vert \mathcal{F}%
_{\overline{\tau}}\right]  \right] \\
& \leq & \mathbb{E}\left[  V_{a}^{\overline{c}}(X_{\overline{\tau}}%
^{D},0)\right]  .
\end{array}
\]
}As before, the last inequality involves a shift of stopping times and it
follows from Theorem 2 of Claisse, Talay and Tan \cite{Claisse}. Then%

\[%
\begin{array}
[c]{lll}%
\mathbb{E}\left[  \int_{\overline{\tau}}^{\tau}D_{s}e^{-qs}ds\right]  & \leq &
\mathbb{E}\left[  \overline{V}(X_{\overline{\tau}}^{D})\right] \\
& \leq & \mathbb{E}\left[  \overline{V}((c_{2}-c_{1})\overline{\tau}\right] \\
& \leq & K_{1}\mathbb{E}[e^{-q\overline{\tau}}\overline{\tau}(c_{2}-c_{1})].
\end{array}
\]
Hence, we can write,%
\begin{equation}%
\begin{array}
[c]{lll}%
V_{a}^{\overline{c}}(x,c_{1})-V_{a}^{\overline{c}}(x,c_{2}) & \leq &
J(x;D)+\varepsilon-J(x;\widetilde{D})\\
& = & \mathbb{E}\left[  \int_{0}^{\overline{\tau}}\left(  D_{s}-\widetilde
{D}_{s}\right)  e^{-qs}ds\right]  +\mathbb{E}\left[  \int_{\overline{\tau}%
}^{\tau}D_{s}e^{-qs}ds\right]  +\varepsilon\\
& \leq & 0+\mathbb{E}\left[  \int_{\overline{\tau}}^{\tau}D_{s}e^{-qs}%
ds\right]  +\varepsilon\\
& \leq & K_{1}E[ae^{-q\overline{\tau}}\overline{\tau}(c_{2}-c_{1}%
)]+\varepsilon\\
& \leq & K_{2}(c_{2}-c_{1})+\varepsilon.
\end{array}
\label{Lipaux1C}%
\end{equation}
So, we deduce (\ref{Lips3}), taking $K_{2}=K_{1}\max_{t\geq0}\{e^{-qt}%
ta\}=K_{1}\frac{e^{-1}}{q}a\ $and $K=K_{1}\max\{\frac{e^{-1}}{q}a,1\}$. We
conclude the result from (\ref{Lips1}), (\ref{Lips2}) and (\ref{Lips3}).

In the case $\overline{c}=\infty$, the result follows from Proposition
\ref{Limite de restricted 1}.\hfill$\blacksquare$\newline

The following lemma states the dynamic programming principle, its proof is
similar to the one of Lemma 1.2 in Azcue and Muler \cite{AM Libro}{\Large .}

\begin{lemma}
\label{Lemma DPP} Given any stopping time $\widetilde{\tau}$, we can write in
both the restricted case $\overline{c}<\infty$ and the unrestricted case
$\overline{c}=\infty,$
\[
V_{a}^{\overline{c}}(x,c)=\sup\limits_{D\in\Pi_{x,c,a}^{[0,\overline{c}]}%
}\mathbb{E}\left[  \int_{0}^{\tau\wedge\widetilde{\tau}}e^{-qs}D_{s}%
ds+e^{-q(\tau\wedge\widetilde{\tau})}V_{a}^{\overline{c}}(X_{\tau
\wedge\widetilde{\tau}}^{D},R_{\tau\wedge\widetilde{\tau}})\right]  \text{.}%
\]

\end{lemma}

We now show a Lipschitz condition of $h(a)=V_{a}^{\overline{c}}(x,c)$ on the
drawdown constant $a\in\lbrack0,1]$, for fixed $x,$ $c$ and finite
$\overline{c}$.

\begin{proposition}
Given $\overline{c}<\infty$ and $a_{1},a_{2}\in\lbrack0,1]$ with $a_{1}<a_{2}%
$, there exists $K_{3}>0$ such that%
\[
0\leq V_{a_{1}}^{\overline{c}}(x,c)-V_{a_{2}}^{\overline{c}}(x,c)\leq
K_{3}\left(  a_{2}-a_{1}\right)  ,
\]
with $K_{3}=\overline{V}^{\prime}(0)\frac{e^{-1}}{q}\overline{c}$ only
depending on $\overline{c}$. In the case $\overline{c}=\infty$, $V_{a}%
^{\infty}(x,c)$ is continuous in $a\in\lbrack0,1]$.
\end{proposition}

\textit{Proof.} Consider first the case $\overline{c}<\infty$. Take
$\varepsilon>0$ and $D\in\Pi_{x,c,a_{1}}^{[0,\overline{c}]}$ such that%

\[
J(x;D)\geq V_{a_{1}}^{\overline{c}}(x,c)-\varepsilon.
\]
Let us consider $\widetilde{D}\in\Pi_{x,c,a_{2}}^{[0,\overline{c}]}\ $defined
as $\widetilde{D}_{t}=\max\{D_{t},a_{2}R_{t}\}$. Denote by $X_{t}%
^{\widetilde{D}}$ the associated controlled surplus process and by
$\overline{\tau}\leq\tau$ the corresponding ruin time. We have that
$0\leq\widetilde{D}_{s}-D_{s}\leq(a_{2}-a_{1})R_{s}$ and so
\[
X_{\overline{\tau}}^{D}=X_{\overline{\tau}}^{D}-X_{\overline{\tau}%
}^{\widetilde{D}}\leq\int_{0}^{\overline{\tau}}(a_{2}-a_{1})R_{s}%
ds=(a_{2}-a_{1})\overline{\tau}\,\overline{c}.
\]
We can write{
\[%
\begin{array}
[c]{lll}%
V_{a_{1}}^{\overline{c}}(x,c)-V_{a_{2}}^{\overline{c}}(x,c) & = &
J(x;D)-J(x;\widetilde{D})+\varepsilon\\
& = & \mathbb{E}\left[  \int_{0}^{\overline{\tau}}e^{-qs}\left(
D_{s}-\widetilde{D}_{s}\right)  ds\right]  +\mathbb{E}\left[  \int
_{\overline{\tau}}^{\tau}e^{-qs}D_{s}ds\right]  +\varepsilon\\
& \leq & 0+\mathbb{E}\left[  \mathbb{E}\left[  \left.  \int_{\overline{\tau}%
}^{\tau}e^{-qs}D_{s}ds\right\vert \mathcal{F}_{\overline{\tau}}\right]
\right]  +\varepsilon\\
& = & \mathbb{E}\left[  \mathbb{E}\left[  \left.  e^{-q\overline{\tau}}%
\int_{0}^{\tau-\overline{\tau}}e^{-qu}D_{\overline{\tau}+u}du\right\vert
\mathcal{F}_{\overline{\tau}}\right]  \right]  +\varepsilon\\
& \leq & \mathbb{E[}e^{-q\overline{\tau}}\overline{V}((a_{2}-a_{1}%
)\overline{\tau}\overline{c})]+\varepsilon\\
& \leq & \mathbb{E[}e^{-q\overline{\tau}}\overline{V}^{\prime}(0)((a_{2}%
-a_{1})\overline{\tau}\overline{c}]+\varepsilon\\
& \leq & \overline{V}^{\prime}(0)\frac{e^{-1}}{q}\overline{c}(a_{2}%
-a_{1})+\varepsilon,
\end{array}
\]
}and one can conclude the result defining $K_{3}=\overline{V}^{\prime}%
(0)\frac{e^{-1}}{q}\overline{c}$.

In the case $\overline{c}=\infty$, we want to show that given $\varepsilon>0$
and $a_{1}\geq0,$ there exists $\delta>0$ such that, if $0<a_{2}-a_{1}<\delta$
then $V_{a_{1}}^{\infty}(x,c)-V_{a_{2}}^{\infty}(x,c)<\varepsilon$. Take
$\overline{c}_{0}$ large enough such that $V_{a_{1}}^{\infty}(x,c)-V_{a_{1}%
}^{\overline{c}_{0}}(x,c)<\varepsilon/2$ and $\delta$ $=\varepsilon
/(2\overline{V}^{\prime}(0)\frac{e^{-1}}{q}\overline{c}_{0})$. Given any
$a_{2}\in(a_{1},a_{1}+\delta),$ we have%
\[%
\begin{array}
[c]{lll}%
V_{a_{1}}^{\infty}(x,c)-V_{a_{2}}^{\infty}(x,c) & = & V_{a_{1}}^{\infty
}(x,c)-V_{a_{1}}^{\overline{c}_{0}}(x,c)+V_{a_{1}}^{\overline{c}_{0}%
}(x,c)-V_{a_{2}}^{\overline{c}_{0}}(x,c)+V_{a_{2}}^{\overline{c}_{0}%
}(x,c)-V_{a_{2}}^{\infty}(x,c)\\
& \leq & \varepsilon/2+\overline{V}^{\prime}(0)\frac{e^{-1}}{q}\overline
{c}_{0}(a_{2}-a_{1})+0\\
& \leq & \varepsilon.
\end{array}
\]
\hfill$\blacksquare$

\begin{remark}
\normalfont
\label{Problemas en a=0 y cbarra=infinito} Note that in the case $a=0,$
Proposition \ref{Limite de restricted 2} does not hold. Indeed, $V_{0}%
^{\infty}(x,c)=\overline{V}(x)$, so that $\lim_{c\rightarrow\infty}%
V_{0}^{\infty}(x,c)=\overline{V}(x)>x.$ Although $\lim_{c\rightarrow\infty
}V_{a}^{\infty}(x,c)=x$ for $a\in(0,1]$ and $\lim_{a\rightarrow0^{+}}%
V_{a}^{\infty}(x,c)=V_{0}^{\infty}(x,c)$ by the previous proposition, the lack
of the Lipschitz property of $V_{a}^{\infty}(x,c)$ at $a=0$ enables the
iterated limits
\[
\lim_{c\rightarrow\infty}\left(  \lim_{a\rightarrow0^{+}}V_{a}^{\infty
}(x,c)\right)  =\overline{V}(x)\;\;\text{and }\lim_{a\rightarrow0^{+}}\left(
\lim_{c\rightarrow\infty}V_{a}^{\infty}(x,c)\right)  =x
\]
to not coincide.
\end{remark}

In the next proposition, we study the continuity of $V_{a}^{\overline{c}%
}(x,c)$ with respect to $\overline{c}.$

\begin{proposition}
\label{Cambios con cbarra}Given $\overline{c}_{1},\overline{c}_{2}\in
\lbrack0,\infty)$ with $\overline{c}_{1}<\overline{c}_{2}<\infty,$ there
exists a $K_{2}>0$ such that%
\[
0\leq V_{a}^{\overline{c}_{2}}(x,c)-V_{a}^{\overline{c}_{1}}(x,c)\leq\frac
{1}{q}\left(  \overline{c}_{2}-\overline{c}_{1}\right)
\]
for $c\leq\overline{c}_{1}.$
\end{proposition}

\textit{Proof.} Take $\varepsilon>0$ and $D\in\Pi_{x,c,a}^{[0,\overline{c}%
_{2}]}$ such that%
\[
J(x;D)\geq V_{a}^{\overline{c}_{2}}(x,c)-\varepsilon,
\]
and denote the ruin time of the process $X_{t}^{D}$ by $\tau$. Let us consider
$\widetilde{D}\in\Pi_{x,c,a}^{[0,\overline{c}_{1}]}\ $as $\widetilde{D}%
_{t}=\min\{D_{t},\overline{c}_{1}\}=\overline{c}_{1}I_{D_{t}>\overline{c}_{1}%
}+D_{t}I_{D_{t}\leq\overline{c}_{1}}$ for $t\leq\tau$ and $\widetilde{D}%
_{t}=\overline{c}_{1}$ for $t>\tau$, denote by $X_{t}^{\widetilde{D}}$ the
associated controlled surplus process and by $\overline{\tau}\geq\tau$ the
corresponding ruin time. We then have $D_{s}-\widetilde{D}_{s}\leq\overline
{c}_{2}-\overline{c}_{1}$ and one can deduce%
\[%
\begin{array}
[c]{lll}%
V_{a}^{\overline{c}_{2}}(x,c)-V_{a}^{\overline{c}_{1}}(x,c) & \leq &
J(x;D)+\varepsilon-J(x;\widetilde{D})\\
& = & \mathbb{E}\left[  \int_{0}^{\tau}\left(  D_{s}-\widetilde{D}_{s}\right)
e^{-qs}ds\right]  -\mathbb{E}\left[  \int_{\tau}^{\overline{\tau}}D_{s}%
e^{-qs}ds\right]  +\varepsilon\\
& \leq & \mathbb{E}\left[  \int_{0}^{\tau}\left(  D_{s}-\widetilde{D}%
_{s}\right)  e^{-qs}ds\right]  +\varepsilon\\
& \leq & \mathbb{E}\left[  \int_{0}^{\tau}\left(  \overline{c}_{2}%
-\overline{c}_{1}\right)  e^{-qs}ds\right]  +\varepsilon\\
& = & \frac{\left(  \overline{c}_{2}-\overline{c}_{1}\right)  }{q}%
\mathbb{E}\left[  1-e^{-q\tau}\right]  +\varepsilon\\
& \leq & \frac{\left(  \overline{c}_{2}-\overline{c}_{1}\right)  }%
{q}+\varepsilon.
\end{array}
\]
\hfill$\blacksquare$

\section{The Hamilton-Jacobi-Bellman equation
\label{Hamilton-Jacobi-Bellman equations}}

In this section we introduce the Hamilton-Jacobi-Bellman (HJB) equation of the
drawdown problem. We show that the optimal value function $V$ defined in
(\ref{Optimal Value Function}) is the unique viscosity solution of the
corresponding HJB equation with suitable boundary conditions.

As we stated in the previous section, the limit case $a=0$ (no drawdown
restriction) has been studied for both $\overline{c}<\infty$ and $\overline
{c}=\infty$, and the case $a=1$ (ratcheting) for $\overline{c}<\infty$.

Define
\begin{equation}
\mathcal{L}^{d}(W)(x,c):=\frac{\sigma^{2}}{2}\partial_{xx}W(x,c)+(\mu
-d)\partial_{x}W(x,c)-qW(x,c)+d. \label{Lc}%
\end{equation}
The HJB equation associated to (\ref{Optimal Value Function}) for both
$\overline{c}<\infty$ and $\overline{c}=\infty$ is given by%

\begin{equation}
\max\{\max_{d\in\lbrack ac,c]}\mathcal{L}^{d}(u)(x,c),\partial_{c}%
u(x,c)\}=0\text{ for }x\geq0\ \text{and }0\leq c<\overline{c}\text{. }
\label{HJB equation 2}%
\end{equation}

\noindent Note that an alternative equivalent formulation is
\begin{equation}
\max\{\mathcal{L}^{c}(u)(x,c),\mathcal{L}^{ac}(u)(x,c),\partial_{c}%
u(x,c)\}=0\text{ for }x\geq0\ \text{and }0\leq c<\overline{c}\text{.}
\label{HJB equation}%
\end{equation}

\noindent For the ratcheting case $a=1$, the HJB equation correspondingly
simplifies to
\[
\max\{\mathcal{L}^{c}(u)(x,c),\partial_{c}u(x,c)\}=0\text{ for }%
x\geq0\ \text{and }0\leq c<\overline{c}.
\]
Let us introduce the usual notion of viscosity solution for the HJB equation
in both cases $0<\overline{c}<\infty$ or $\overline{c}=\infty$.

\begin{definition}
\label{Viscosity}

(a) A locally Lipschitz function $\overline{u}:[0,\infty)\times\lbrack
0,\overline{c})\rightarrow{\mathbb{R}}$\ is a viscosity supersolution of
(\ref{HJB equation})\ at $(x,c)\in(0,\infty)\times\lbrack0,\overline{c})$,\ if
any (2,1)-differentiable function $\varphi:[0,\infty)\times\lbrack
0,\overline{c})\rightarrow{\mathbb{R}}\ $with $\varphi(x,c)=\overline{u}(x,c)$
such that $\overline{u}-\varphi$\ reaches the minimum at $\left(  x,c\right)
$\ satisfies
\[
\max\left\{  \mathcal{L}^{c}(\varphi)(x,c),\mathcal{L}^{ac}(\varphi
)(x,c),\partial_{c}\varphi(x,y)\right\}  \leq0.\
\]
The function $\varphi$ is called a \textbf{test function for supersolution} at
$(x,c)$.

(b) A function $\underline{u}:$ $[0,\infty)\times\lbrack0,\overline
{c})\rightarrow{\mathbb{R}}\ $\ is a viscosity subsolution\ of
(\ref{HJB equation})\ at $(x,c)\in(0,\infty)\times\lbrack0,\overline{c})$,\ if
any (2,1)-differentiable function $\psi:[0,\infty)\times\lbrack0,\overline
{c})\rightarrow{\mathbb{R}}\ $with $\psi(x,c)=\underline{u}(x,c)$ such that
$\underline{u}-\psi$\ reaches the maximum at $\left(  x,c\right)  $ satisfies
\[
\max\left\{  \mathcal{L}^{c}(\psi)(x,c),\mathcal{L}^{ac}(\psi)(x,c),\partial
_{c}\psi(x,c)\right\}  \geq0\text{.}%
\]
The function $\psi$ is called a \textbf{test function for subsolution} at
$(x,c)$.

(c) A function $u:[0,\infty)\times\lbrack0,\overline{c})\rightarrow
{\mathbb{R}}$ which is both a supersolution and subsolution at $(x,c)\in
\lbrack0,\infty)\times\lbrack0,\overline{c})$ is called a viscosity solution
of (\ref{HJB equation})\ at $(x,c)$.
\end{definition}

\subsection{HJB equation with bounded dividend rates
\label{HJB equation with bounded dividend rates}}

Given $a\in(0,1]$ and $\overline{c}<\infty$, we denote for in the sequel, for
simplicity of exposition,
\begin{equation}
\Pi_{x,c}:=\Pi_{x,c,a}^{[0,\overline{c}]}\text{ and }V:=V_{a}^{\overline{c}}.
\label{Optimal Value Function Bounded}%
\end{equation}
Here the state variables are the current surplus and the running maximum
dividend rate. The results of this subsection for the case $a=1$ (ratcheting
dividend constraint) were already proved in \cite{AAM21}.

In the next proposition we state that $V$ is a viscosity solution of the
corresponding HJB equation.

\begin{proposition}
\label{Proposicion Viscosidad} $V$ is a viscosity solution of
(\ref{HJB equation}) in $(0,\infty)\times\lbrack0,\overline{c})$.
\end{proposition}

\textit{Proof.} Let us show first that $V$ is a viscosity supersolution in
$(0,\infty)\times\lbrack0,\overline{c})$ . By Proposition
\ref{Monotone Optimal Value Function}, $\partial_{c}V\leq0$ in $(0,\infty
)\times\lbrack0,\overline{c})$ in the viscosity sense.

Consider now $(x,c)\in(0,\infty)\times\lbrack0,\overline{c})$ and the
admissible strategy $D\in\Pi_{x,c}$, which pays dividends at constant rate
$d\in\lbrack ac,c]$ up to the ruin time $\tau$. Let $X_{t}^{D}$ be the
corresponding controlled surplus process and suppose that there exists a test
function $\varphi$ for supersolution (\ref{HJB equation}) at $(x,c),$ then
$\varphi\leq V$ and $\varphi(x,c)=V(x,c)$. {We want to prove that
$\mathcal{L}^{d}\mathcal{(}\varphi)(x,c)\leq0$. For that purpose, we consider
an auxiliary test function for the supersolution $\tilde{\varphi}$ in such a
way that $\tilde{\varphi}\leq\varphi\leq V$ in $[0,\infty)\times
\lbrack0,\overline{c}]$, $\tilde{\varphi}=\varphi$ in $[0,2x]$ (so
$\mathcal{L}^{d}\mathcal{(}\varphi)(x,c)=\mathcal{L}^{d}\mathcal{(}%
\tilde{\varphi})(x,c)$) and }$\mathcal{L}^{d}\mathcal{(}{\tilde{\varphi}%
})(\cdot,c)${ is bounded in $[0,\infty)$. We introduce $\tilde{\varphi}$
because }$\mathcal{L}^{d}\mathcal{(}\varphi)(\cdot,c)${ may be unbounded in
$[0,\infty)$. We construct $\tilde{\varphi}$ as follows: take $g:[0,\infty
)\rightarrow\lbrack0,1]$ twice continuously differentiable with $g=0$ in
$[2x+1,\infty)$ and $g=1$ in $[0,2x]$, and define $\tilde{\varphi}(y,d)=$
$\varphi(y,d)g(y)$.}

Using Lemma \ref{Lemma DPP}, we obtain for $h>0$%

\[%
\begin{array}
[c]{lll}%
\tilde{\varphi}(x,c) & = & V(x,c)\\
& \geq & \mathbb{E}\left[  \int\nolimits_{0}^{\tau\wedge h}de^{-q\,s}%
\,ds\right]  +\mathbb{E}\left[  e^{-q(\tau\wedge h)}\tilde{\varphi}%
(X_{\tau\wedge h}^{D},c)\right]  \text{.}%
\end{array}
\]
Hence, we get using It\^{o}'s formula%

\[%
\begin{array}
[c]{lll}%
0 & \geq & \mathbb{E}\left[  \int\nolimits_{0}^{\tau\wedge h}e^{-q\,s}%
\,ds\right]  +\mathbb{E}\left[  e^{-q(\tau\wedge h)}\tilde{\varphi}%
(X_{\tau\wedge h}^{D},c)-\tilde{\varphi}(x,c)\right] \\
& = & \mathbb{E}\left[  \int\nolimits_{0}^{\tau\wedge h}de^{-q\,s}%
\,\,ds\right]  +\mathbb{E}\left[  \int\nolimits_{0}^{\tau\wedge h}%
e^{-q\,s}(\frac{\sigma^{2}}{2}\partial_{xx}\tilde{\varphi}(X_{s}%
^{D},c)+\partial_{x}\tilde{\varphi}(X_{s}^{D},c)(\mu-d)-q\tilde{\varphi}%
(X_{s}^{D},c))ds\right] \\
&  & +\mathbb{E}\left[  \int_{0}^{\tau\wedge h}\partial_{x}\tilde{\varphi
}(X_{s}^{D},c)\sigma dWs~\right] \\
& = & \mathbb{E}\left[  \int\nolimits_{0}^{\tau\wedge h}e^{-q\,s}%
\mathcal{L}^{d}\mathcal{(}\tilde{\varphi})(X_{s}^{D},c)ds\right]  \text{.}%
\end{array}
\]
Since $\tau>0$ a.s.,%
\[
\left\vert \frac{1}{h}\int\nolimits_{0}^{\tau\wedge h}e^{-q\,s}\mathcal{L}%
^{d}\mathcal{(}\tilde{\varphi})(X_{s}^{D},c)ds\right\vert \leq\sup
_{y\in\lbrack0,\infty)}\left\vert \mathcal{L}^{d}\mathcal{(}\tilde{\varphi
})(y,c)\right\vert ,
\]
and
\[
\lim_{h\rightarrow0^{+}}\frac{1}{h}\int\nolimits_{0}^{\tau\wedge h}%
e^{-q\,s}\mathcal{L}^{d}\mathcal{(}\tilde{\varphi})(X_{s}^{D},c)ds=\mathcal{L}%
^{d}\mathcal{(}\tilde{\varphi})(x,c)~\text{a.s..}%
\]
We conclude, using the bounded convergence theorem, that $\mathcal{L}%
^{d}\mathcal{(}\varphi)(x,c)=\mathcal{L}^{d}\mathcal{(}\tilde{\varphi
})(x,c)\leq0$ for any $d\in\lbrack ac,c]$; so $V$ is a viscosity supersolution
at $(x,c)$.

We skip the proof that $V$ is a viscosity subsolution in $(0,\infty
)\times\lbrack0,\overline{c})$, because it is similar to the one of
Proposition 3.1 in \cite{AAM21}. \hfill$\blacksquare$\newline

Let us consider the function
\begin{equation}
v^{\overline{c}}(x):=V(x,\overline{c}):[0,\infty)\rightarrow\lbrack0,\infty).
\label{Optimal tapa}%
\end{equation}
In the next proposition, we state a comparison result for the viscosity
solutions of (\ref{HJB equation}) for $\overline{c}>0$. The proof is similar
to the one of Lemma 3.2 of \cite{AAM21}.

\begin{lemma}
\label{Lema para Unicidad} Assume that (i) $\underline{u}$ is a viscosity
subsolution and $\overline{u}$ is a viscosity supersolution of the HJB
equation (\ref{HJB equation}) for all $x>0$ and for all $c\in\lbrack
0,\overline{c})$, (ii) $\underline{u}$ and $\overline{u}$ are non-decreasing
in the variable $x$ and Lipschitz in $[0,\infty)\times\lbrack0,\overline{c}]$,
(iii) $\underline{u}(0,c)=\overline{u}(0,c)=0$, $\lim_{x\rightarrow\infty
}\underline{u}(x,c)\leq\overline{c}/q\leq\lim_{x\rightarrow\infty}\overline
{u}(x,c)$ and (iv) $\underline{u}(x,\overline{c})\leq v^{\overline{c}}%
(x)\leq\overline{u}(x,\overline{c})$ for $x\geq0$. Then $\underline{u}%
\leq\overline{u}$ in $[0,\infty)\times\lbrack0,\overline{c}).$
\end{lemma}

The following characterization theorem is a direct consequence of the previous
lemma and Propositions \ref{Monotone Optimal Value Function} and
\ref{Proposicion Viscosidad}.

\begin{theorem}
\label{Caracterizacion Continua}The optimal value function $V$ is the unique
function non-decreasing in $x$ that is a viscosity solution of
(\ref{HJB equation}) in $(0,\infty)\times\lbrack0,\overline{c})$ with
$V(0,c)=0,$ $V(x,\overline{c})=v^{\overline{c}}(x)$ and $\lim_{x\rightarrow
\infty}$ $V(x,c)=\overline{c}/q$ for $c\in\lbrack0,\overline{c}).$
\end{theorem}

From Definition \ref{Optimal Value Function}, Lemma \ref{Lema para Unicidad},
and Proposition \ref{Monotone Optimal Value Function} together with
Proposition \ref{Proposicion Viscosidad}, we also get the following
verification theorem.

\begin{theorem}
\label{verification result} Consider a family of strategies
\[
\left\{  C_{x,c}\in\Pi_{x,c}:(x,c)\in\lbrack0,\infty)\times\lbrack
0,\overline{c}]\right\}  .
\]
If the function $W(x,c):=J(x;C_{x,c})$ is a viscosity supersolution of the HJB
equation (\ref{HJB equation}) in $(0,\infty)\times\lbrack0,\overline{c})$ with
$\lim_{x\rightarrow\infty}W(x,c)=$ $\overline{c}/q,$ then $W$ is the optimal
value function $V$. Also, if for each $k\geq1$ there exists a family of
strategies $\left\{  C_{x,c}^{k}\in\Pi_{x,c}:(x,c)\in\lbrack0,\infty
)\times\lbrack0,\overline{c}]\right\}  $ such that $W(x,c):=\lim
_{k\rightarrow\infty}J(x;C_{x,c}^{k})$ is a viscosity supersolution of the HJB
equation (\ref{HJB equation}) in $(0,\infty)\times\lbrack0,\overline{c})$ with
$\lim_{x\rightarrow\infty}W(x,c)=$ $\overline{c}/q$, then $W$ is the optimal
value function $V$.
\end{theorem}

\subsection{HJB equation with unbounded dividend rates}

Let us now consider the case $\overline{c}=\infty$ with $a\in(0,1]$. Since $a$
is fixed, we denote $V^{\infty}=V_{a}^{\infty}$. The proof of the following
proposition is similar to the one of the case with bounded dividend rate.

\begin{proposition}
\label{Proposicion Viscosidad unbounded} $V^{\infty}$ is a viscosity solution
of (\ref{HJB equation}) for any $(x,c)\in$ $(0,\infty)\times\lbrack0,\infty)$.
\end{proposition}

We now state a comparison result for the unbounded case.

\begin{lemma}
\label{Lemma Unicidad Unbounded} Assume that (i) $\underline{u}$ is a
viscosity subsolution and $\overline{u}$ is a viscosity supersolution of the
HJB equation (\ref{HJB equation}) for all $x>0$ and for all $c\in
\lbrack0,\infty)$, (ii) $\underline{u}$ and $\overline{u}$ are non-decreasing
in the variable $x$ and Lipschitz in $[0,\infty)\times\lbrack0,\infty)$, (iii)
$\underline{u}(0,c)=\overline{u}(0,c)=0$, (iv) $~\underline{u}(x,c)\leq
x+{\mu}/{q}$, $x\leq\overline{u}(x,c)$ and (v) $\lim_{c\rightarrow\infty
}\underline{u}(x,c)\leq x\leq\lim_{c\rightarrow\infty}\overline{u}(x,c)$ for
$x\geq0$. Then $\underline{u}\leq\overline{u}$ in $[0,\infty)\times
\lbrack0,\infty).$
\end{lemma}

\textit{Proof.} Suppose that there is a point $(x_{0},c_{0})\in(0,\infty
)\times(0,\infty)$ such that $\underline{u}(x_{0},c_{0})-\overline{u}%
(x_{0},c_{0})>0$. We prove here that the $\max\nolimits_{x\geq0,c\geq0}\left(
\underline{u}(x,c)-\overline{u}^{s_{0}}(x,c)\right)  $ is achieved in a
bounded set. From this we get a contradiction following the arguments of the
proof of Lemma 3.2 of \cite{AAM21}.

Let us define
\[
h(c)=1+(\frac{\underline{u}(x_{0},c_{0})-\overline{u}(x_{0},c_{0})}%
{2\overline{u}(x_{0},c_{0})})e^{-{c}}>1\ \text{and }\overline{u}%
^{s}(x,c)=s\,h(c)\,\overline{u}(x,c)
\]
for any $s>1$. We have
\[%
\begin{array}
[c]{lll}%
\underline{u}(x_{0},c_{0})-\overline{u}^{s}(x_{0},c_{0}) & = & \underline
{u}(x_{0},c_{0})-\left(  1+\frac{\underline{u}(x_{0},c_{0})-\overline{u}%
(x_{0},c_{0})}{2\overline{u}(x_{0},c_{0})}e^{-{c}}\right)  s\,\,\overline
{u}(x_{0},c_{0})\\
& = & (1-\frac{e^{-{c}}\,s}{2})\left(  \underline{u}(x_{0},c_{0}%
)-s\overline{u}(x_{0},c_{0})\right) \\
& > & 0
\end{array}
\]
for $s\in(1,2)$.

Let us show now that $\overline{u}^{s}$ is a strict supersolution. We have
that $\varphi$ is a test function for the supersolution of $\overline{u}$ at
$(x,c)$ if and only if $\varphi^{s}:=s\,h(c)\,\varphi$ is a test function for
the supersolution of $\overline{u}^{s}$ at $(x,c)$. Moreover,
\begin{equation}
\mathcal{L}^{d}(\varphi^{s})(x,c)=sh(c)\mathcal{L}^{d}(\varphi
)(x,c)+d(1-sh(c))<0, \label{Desigualdad L1 no acotado}%
\end{equation}
for $d\in\lbrack ac,c]$ and%
\begin{equation}
\partial_{c}\varphi^{s}(x,c)\leq-s(h(c)-1)\varphi(x,c)<0
\label{Desigualdad L2 no acotado}%
\end{equation}
since $\varphi(x,c)=\overline{u}(x,c)\geq x>0$.

Take $s_{0}>1$ such that $\underline{u}(x_{0},c_{0})-\overline{u}^{s_{0}%
}(x_{0},c_{0})>0$. We define%
\begin{equation}
M:=\sup\limits_{x\geq0,c\geq0}\left(  \underline{u}(x,c)-\overline{u}^{s_{0}%
}(x,c)\right)  . \label{Definicion de M no acotado}%
\end{equation}

Let us show that
\begin{equation}
\arg\max_{x\geq0,c\geq0}\left(  \underline{u}(x,c)-\overline{u}^{s_{0}%
}(x,c)\right)  \in(0,b)\times(0,c_{1}) \label{(x,c) en compacto}%
\end{equation}
for some positive $b$ and $c_{1}$. Since$~\underline{u}(x,c)\leq x+\frac{\mu
}{q}$ and $x\leq\overline{u}(x,c),$%

\[%
\begin{array}
[c]{lll}%
\underline{u}(x,c)-\overline{u}^{s_{0}}(x,c) & \leq & \left(  x+\frac{\mu}%
{q}\right)  -s_{0}h(c)x\\
& < & x(1-s_{0})+\frac{\mu}{q}\\
& < & 0
\end{array}
\]
for $x$ large enough, so there exists a $b>x_{0}$ such that $\arg\max
_{x\geq0,c\geq0}\left(  \underline{u}(x,c)-\overline{u}^{s_{0}}(x,c)\right)
\in(0,b)\times(0,\infty)$. Besides, we have that the function
\[
g(c):=\max_{x\geq0}\{\underline{u}(x,c)-\overline{u}^{s_{0}}(x,c)\}=\max
_{x\in(0,b)}\{\underline{u}(x,c)-\overline{u}^{s_{0}}(x,c)\}
\]
satisfies that $\limsup$ $_{c\rightarrow\infty}$ $g(c)\leq0$ because
$\lim_{c\rightarrow\infty}\underline{u}(x,c)\leq x\leq\lim_{c\rightarrow
\infty}\overline{u}(x,c)$ for $x\geq0$, so there exists a $c_{1}>0$ such that
$g(c)\leq\frac{M}{2}$ for $c\geq c_{1}$ and then we conclude
(\ref{(x,c) en compacto}).

Hence, we obtain that the maximum is achieved in a bounded set, that is%

\[
0<\underline{u}(x_{0},c_{0})-\overline{u}^{s_{0}}(x_{0},c_{0})\leq
M=\max\limits_{x\in(0,b)\times(0,c_{1})}\left(  \underline{u}(x,c)-\overline
{u}^{s_{0}}(x,c)\right)  .
\]
\hfill$\blacksquare$\newline

As for bounded dividend rates, the following characterization theorem is a
direct consequence of the previous lemma, Remark \ref{Limite de restricted 2}
and Proposition \ref{Proposicion Viscosidad unbounded}.

\begin{theorem}
\label{Caracterizacion Continua Unbounded}The optimal value function
$V^{\infty}$ is the unique function non-decreasing in $x$ that is a viscosity
solution of (\ref{HJB equation}) in $(0,\infty)\times\lbrack0,\infty)$ with
$V^{\infty}(0,c)=0,$ $V^{\infty}(x,\overline{c})-x\ $bounded and
$\lim_{c\rightarrow\infty}$ $V^{\infty}(x,c)=x.$
\end{theorem}

From Definition \ref{Optimal Value Function}, Lemma
\ref{Lemma Unicidad Unbounded}, and Remark \ref{Limite de restricted 2}
together with Proposition \ref{Proposicion Viscosidad unbounded}, we also get
the following verification theorem.

\begin{theorem}
\label{verification result Unbounded} Consider a family of strategies
\[
\left\{  C_{x,c}\in\Pi_{x,c}:(x,c)\in\lbrack0,\infty)\times\lbrack
0,\infty)\right\}  .
\]
If the function $W(x,c):=J(x;C_{x,c})$ is a viscosity supersolution of the HJB
equation (\ref{HJB equation}) in $(0,\infty)\times\lbrack0,\infty)$ with
$W(x,c)\geq x$, then $W$ is the optimal value function $V^{\infty}$. Also, if
for each $k\geq1$ there exists a family of strategies $\left\{  C_{x,c}^{k}%
\in\Pi_{x,c}:(x,c)\in\lbrack0,\infty)\times\lbrack0,\infty)\right\}  $ such
that $W(x,c):=\lim_{k\rightarrow\infty}J(x;C_{x,c}^{k})$ is a viscosity
supersolution of the HJB equation (\ref{HJB equation}) in $(0,\infty
)\times\lbrack0,\infty)$ with $W(x,c)\geq x$, then $W$ is the optimal value
function $V^{\infty}$.
\end{theorem}

\section{Refracting dividend strategies and $v^{\overline{c}}$}

\label{sec5}

In the case $0<\overline{c}<\infty$ and $a\in(0,1),$ we now want to
investigate further the function $v^{\overline{c}}$ (defined in
(\ref{Optimal tapa})) of paying dividends with rates in $d\in\lbrack
a\overline{c},\overline{c}]$ in an optimal way. The following characterization
is the one-dimensional version of the results of Section
\ref{HJB equation with bounded dividend rates}.

\begin{proposition}
\label{HJB tapa} The function $v^{\overline{c}}:[0,\infty)\rightarrow
{\mathbb{R}}$ is the unique viscosity solution of
\[
\max\left\{  \mathcal{L}^{\overline{c}}(W)(x),\mathcal{L}^{a\overline{c}%
}(W)(x)\right\}  =0
\]
with boundary conditions $W(0)=0$ and $\lim_{x\rightarrow\infty}%
W(x)=\overline{c}/q.$
\end{proposition}

We present in this section a formula for $v^{\overline{c}}$, which turns out
to be the value function of the optimal refracting strategy as derived in
\cite{ABB}.

The functions $W$ that satisfy $\mathcal{L}^{d}(W)=0$ are given by
\begin{equation}
\frac{d}{q}+a_{1}e^{\theta_{1}(d)x}+a_{2}e^{\theta_{2}(d)x}\text{ with }%
a_{1},a_{2}\in{\mathbb{R}}, \label{Solucion General L=0}%
\end{equation}
where $\theta_{1}(d)>0$ and $\theta_{2}(d)<0$ are the roots of the
characteristic equation
\[
\frac{\sigma^{2}}{2}z^{2}+(\mu-d)z-q=0
\]
associated to the operator $\mathcal{L}^{d}$, that is%
\begin{equation}
\theta_{1}(d):=\frac{d-\mu+\sqrt{(d-\mu)^{2}+2q\sigma^{2}}}{\sigma^{2}}%
,\quad\theta_{2}(d):=\text{$\frac{d-\mu-\sqrt{(d-\mu)^{2}+2q\sigma^{2}}%
}{\sigma^{2}}$.} \label{Definicion tita1 tita2}%
\end{equation}
The following are basic properties of $\theta_{1}(d)$ and $\theta_{2}(d)$:

\begin{enumerate}
\item $\theta_{1}(d)=-\theta_{2}(d)$ if $d=\mu$ and $\theta_{1}^{2}%
(d)\geq\theta_{2}^{2}$$(d)$ if, and only if, $d-\mu\geq0.$

\item $\theta_{1}^{\prime}(d)=\frac{1}{\sigma^{2}}(1+\frac{d-\mu}{\sqrt
{(d-\mu)^{2}+2q\sigma^{2}}})$ and $\theta_{2}^{\prime}(d)=\frac{1}{\sigma^{2}%
}(1-\frac{d-\mu}{\sqrt{(d-\mu)^{2}+2q\sigma^{2}}}),$ so $\theta_{1}^{\prime
}(d),\theta_{2}^{\prime}(d)\in(0,\frac{2}{\sigma^{2}}).$
\end{enumerate}

\noindent The solutions of $\mathcal{L}^{d}(W)=0$ with boundary condition
$W(0)=0$ are then of the more specific form%
\begin{equation}
\frac{d}{q}\left(  1-e^{\theta_{2}(d)x}\right)  +A(e^{\theta_{1}%
(d)x}-e^{\theta_{2}(d)x})\ \text{with }A\in{\mathbb{R}}.
\label{Solucion de L=0 con condicion en 0}%
\end{equation}
Finally, the unique solution of $\mathcal{L}^{d}(W)=0$ with boundary
conditions $W(0)=0$ and $\lim_{x\rightarrow\infty}$ $W(x)=d/q\ $corresponds to
$A=0$, so that%
\begin{equation}
W(x)=\frac{d}{q}\left(  1-e^{\theta_{2}(d)x}\right)  .
\label{Formula de V con c constante}%
\end{equation}
We have that $W$ is increasing and concave in $[0,\infty)$.

In \cite[Th.3.1]{ABB}, the value function of a 'refracting strategy' that pays
$a\overline{c}$ when the current surplus is below a refracting threshold $b$
and pays $\overline{c}$ when the current surplus is above $b$ was shown to be%
\begin{equation}
v(x,\overline{c},b)=\Big(B(\overline{c},b)W_{0}(x,\overline{c})+\frac
{a\overline{c}}{q}(1-e^{\theta_{2}(a\overline{c})x})\Big)I_{x<b}%
+\Big(\frac{\overline{c}}{q}+D(\overline{c},b)e^{\theta_{2}(\overline{c}%
)x}\Big)I_{x\geq b}, \label{Formula value function refraction en b}%
\end{equation}
where%
\[
W_{0}(x,\overline{c})=\frac{e^{\theta_{1}(a\overline{c})x}-e^{\theta
_{2}(a\overline{c})x}}{\sqrt{(\mu-a\overline{c})^{2}+2q\sigma^{2}}},
\]%
\begin{equation}
B(\overline{c},b)=\frac{1}{q}\frac{a\overline{c}e^{\theta_{2}(a\overline{c}%
)b}\left(  \theta_{2}(a\overline{c})-\theta_{2}(\overline{c})\right)
-(1-a)\overline{c}\theta_{2}(\overline{c})}{\partial_{x}W_{0}(b,\overline
{c})-\theta_{2}(\overline{c})W_{0}(b,\overline{c})}, \label{Definicion B}%
\end{equation}
and%
\[
D(\overline{c},b)=B(\overline{c},b)e^{-\theta_{2}(\overline{c})b}%
W_{0}(b,\overline{c})-\frac{a\overline{c}}{q}e^{(\theta_{2}(a\overline
{c})-\theta_{2}(\overline{c}))b}-\frac{(1-a)\overline{c}}{q}e^{-\theta
_{2}(\overline{c})b}.
\]
The optimal threshold $b^{\ast}(\overline{c})$ corresponds to
\begin{equation}
b^{\ast}(\overline{c})=\arg\max_{b\geq0}v(x,\overline{c},b).
\label{Ecuacion bestrella}%
\end{equation}
In case it is positive, by (\ref{Formula value function refraction en b}) this
is the value of $b$ satisfying
\begin{equation}
\partial_{b}B(\overline{c},b)=0. \label{Ecuacion implicita de bestrella}%
\end{equation}

\noindent From \cite{ABB}, we know that the threshold can be characterized as
the \textbf{unique} $b$ such that $v(x,\overline{c},b)$ is twice continuously
differentiable in $x=b$. \ Hence, since $v(x,\overline{c},b^{\ast}%
(\overline{c}))$ is twice continuously differentiable with $v(0,\overline
{c},b^{\ast}(\overline{c}))=0$, $\lim_{x\rightarrow\infty}v(x,\overline
{c},b^{\ast}(\overline{c}))=\overline{c}/q$ and it is also a solution of
\[
\max\left\{  \mathcal{L}^{\overline{c}}(W)(x),\mathcal{L}^{a\overline{c}%
}(W)(x)\right\}  =0,
\]
by Proposition \ref{HJB tapa} we have that
\begin{equation}
v^{\overline{c}}(x)=v(x,\overline{c},b^{\ast}(\overline{c})).
\label{V cbarra formula}%
\end{equation}
That is, the strategy achieving $v^{\overline{c}}$
has a 'refracting' threshold structure with optimal threshold $b^{\ast}(c).$
%

Note also, that since $v^{\overline{c}}$ is twice continuously differentiable
at $b^{\ast}(\overline{c})$ and $\mathcal{L}(v^{\overline{c}})(b^{\ast
}(\overline{c}))=\mathcal{L}^{a\overline{c}}(v^{\overline{c}})(b^{\ast
}(\overline{c}))=0,$ then $\partial_{x}v^{\overline{c}}(b^{\ast}(\overline
{c}))=1$. Also, since
\[
\mathcal{L}^{a\overline{c}}(v^{\overline{c}})(x)-\mathcal{L}^{\overline{c}%
}(v^{\overline{c}})(x)=c(1-a)\left(  \partial_{x}v^{\overline{c}}(x)-1\right)
\]
we obtain
\begin{equation}
\partial_{x}v^{\overline{c}}(x)\geq1\text{ for }x\leq b^{\ast}(\overline
{c})\text{ and }\partial_{x}v^{\overline{c}}(x)\leq1\text{ for }x\geq b^{\ast
}(\overline{c})_{.} \label{Derivada en tapa}%
\end{equation}

\section{Curve strategies and the optimal two-curve strategy for bounded
dividend rates \label{Curve strategies and optimal curve}}

Using the formulas of the previous section, we can find the optimal value
function defined in (\ref{Optimal Value Function Bounded}).

\begin{remark}
\label{Rermark caso interesante}\normalfont Before proceeding, note that this
problem is only interesting for $\overline{c}>q\sigma^{2}/(2\mu)$, as for
smaller values of $\overline{c}$ we know from \cite[Eqn.1.8]{asmtak}
(translated to our notation) that even without a drawdown constraint it is
optimal to pay dividends at maximal rate $\overline{c}$ until the time of
ruin. This then also is the optimal strategy in our situation, as the drawdown
constraint does not affect its applicability. Indeed, and as a self-contained
derivation of this result in the present context, the value function of that
strategy fulfills
\begin{equation}%
\begin{array}
[c]{lll}%
\mathcal{L}^{d}(\frac{\overline{c}}{q}\left(  1-e^{\theta_{2}(\overline{c}%
)x}\right)  )(x) & = & (\overline{c}-d)(-\frac{\overline{c}}{q}\theta
_{2}(\overline{c})e^{\theta_{2}(\overline{c})x}-1)\\
& \leq & (\overline{c}-d)(-\frac{\overline{c}}{q}\theta_{2}(\overline{c})-1)\\
& \leq & 0
\end{array}
\label{Techo1}%
\end{equation}
for both $d=ac$ and $d=c$. So, by Proposition \ref{HJB tapa}, $v^{\overline
{c}}(x)=\frac{\overline{c}}{q}\left(  1-e^{\theta_{2}(\overline{c})x}\right)
$ and $b^{\ast}(\overline{c})=0$. With the notation $U(x,c):=\frac
{\overline{c}}{q}\left(  1-e^{\theta_{2}(\overline{c})x}\right)  $, by Theorem
\ref{Caracterizacion Continua} it is then sufficient to prove that
\[
\max\{\mathcal{L}^{ac}(U)(x,c),\mathcal{L}^{c}(U)(x,c),\partial_{c}%
U(x,c)\}\leq0
\]
for any $c\in\lbrack0,\overline{c})$, but this indeed follows from
(\ref{Techo1}).
\end{remark}


\noindent In the rest of this paper, we will therefore assume that
$\overline{c}>\frac{q\sigma^{2}}{2\mu}.$\newline

\noindent The way in which the optimal value function $V(x,c)$ solves the HJB
equation (\ref{HJB equation}) suggests that the state space $[0,\infty
)\times\lbrack0,\overline{c}]$ is partitioned into two regions: a
\textit{non-change running maximum dividend region }$\mathcal{NC}^{\ast}$ in
which the running maximum dividend rate $c$ does not change and a
\textit{change dividend region }$\mathcal{CH}^{\ast}$ in which the dividend
rate exceeds $c$ (so the running maximum dividend rate increases). Moreover,
the region $\mathcal{NC}^{\ast}$ splits into two connected subregions:
$\mathcal{NC}_{ac}^{\ast}$ in which the dividends are paid at constant rate
$\ ac$ and $\mathcal{NC}_{c}^{\ast}$ in which the dividends are paid at
constant rate $c$.

Roughly speaking, the interior of the region $\mathcal{NC}_{ac}^{\ast}$
consists of the points $(x,c)$ in the state space where $\mathcal{L}%
^{ac}(V)(x,c)=0$, $\mathcal{L}^{c}(V)(x,c)<0$ and $\partial_{c}V<0$; ~the
interior of the region $\mathcal{NC}_{c}^{\ast}$ consists of the points
$(x,c)$ in the state space where $\mathcal{L}^{c}(V)(x,c)=0,$ $\mathcal{L}%
^{ac}(V)(x,c)<0$ and $\partial_{c}V<0;$ and the interior of $\mathcal{CH}%
^{\ast}$ consists of the points where $\partial_{c}V=0,$ $\mathcal{L}%
^{c}(V)(x,c)<0$ and $\mathcal{L}^{ac}(V)(x,c)<0$. We introduce a family of
stationary strategies (or limit of stationary strategies) where the different
dividend payment regions are connected and split by two free boundary curves.

\begin{figure}[tbh]
\par
\begin{center}
\includegraphics[width=7cm]{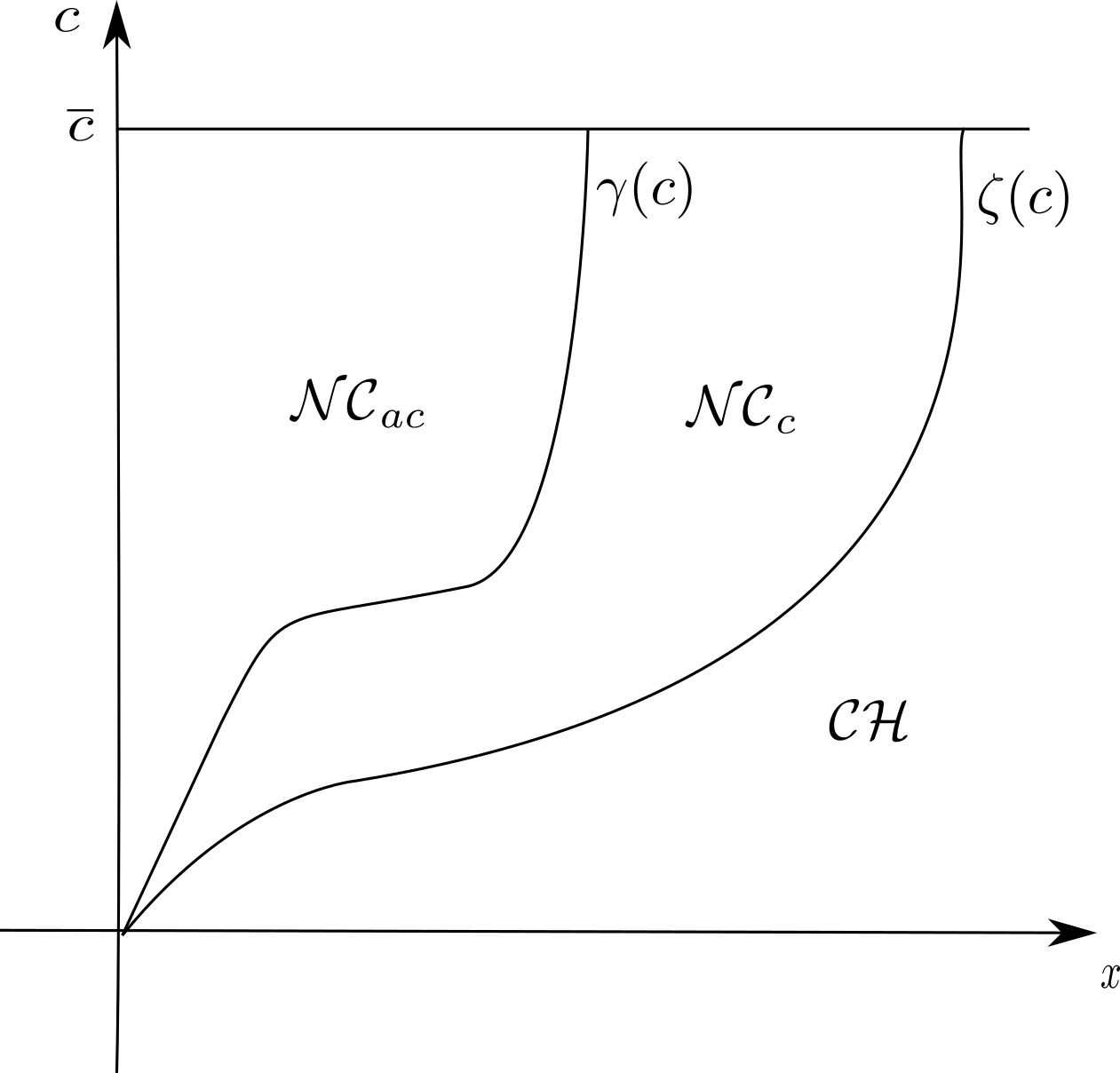}
\end{center}
\caption{{\protect\small A two-curve strategy with its regions.}}%
\label{FigRegions}%
\end{figure}

Let us consider the two functions $\gamma:[0,\overline{c}]\rightarrow
(0,\infty)~$continuously~differentiable, $~\zeta:[0,\overline{c}%
]\rightarrow(0,\infty)~$ bounded, Riemann integrable and c\`{a}dl\`{a}g, and
let us define the set%
\begin{equation}
\mathcal{B}=\{(\gamma,\zeta)~\text{s.t. }\gamma\leq\zeta\text{ and}%
\lim_{c\rightarrow\overline{c}^{-}}\zeta(c)=\zeta(\overline{c})\}.
\label{Conjunto B}%
\end{equation}
In the first part of this section, we define the \textit{function}
$W^{\gamma,\zeta}:[0,\infty)\times\lbrack0,\overline{c}]\rightarrow
\lbrack0,\infty)$ for each $(\gamma,\zeta)\in$ $\mathcal{B}$. We will see
that, in some sense, $W^{\gamma,\zeta}(x,c)$ is a value function of the
\textit{two-curve strategy} which pays dividends at constant rate $ac$ for the
points to the left of the curve $\mathcal{R}(\gamma)$, pays dividends at
constant rate $c$ in between the curves $\mathcal{R}(\gamma)$ and
$\mathcal{R}(\zeta)$ and pays more than $c$ as dividend rate otherwise, where
\[
\mathcal{R}(g)=\left\{  (g(c),c):c\in\lbrack0,\overline{c}]\right\}  .
\]
Hence, the curves $\mathcal{R}(\gamma)$ and $\mathcal{R}(\zeta)$ split the
state space $[0,\infty)\times\lbrack0,\overline{c})$ into three connected
regions:%
\[
\mathcal{NC}_{ac}\mathcal{(\gamma},\mathcal{\zeta)=}\{(x,c)\in\lbrack
0,\infty)\times\lbrack0,\overline{c}):0\leq x<\gamma(c)\}
\]
where dividends are paid with constant rate $ac$,%
\[
\mathcal{NC}_{c}\mathcal{(\mathcal{\gamma}},\mathcal{\mathcal{\zeta}%
)=}\{(x,c)\in\lbrack0,\infty)\times\lbrack0,\overline{c}):\gamma(c)\leq
x<\zeta(c)\}
\]
where dividends are paid with constant rate $c$, and
\[
\mathcal{CH(\mathcal{\mathcal{\gamma}}},\mathcal{\mathcal{\mathcal{\zeta}}%
)=}\{(x,c)\in\lbrack0,\infty)\times\lbrack0,\overline{c}):x\geq\zeta(c)\},
\]
cf. Figure \ref{FigRegions}. Let us call $\mathcal{NC(\gamma},\mathcal{\zeta
)=NC}_{ac}\mathcal{(\gamma},\mathcal{\zeta)\cup NC}_{c}\mathcal{(\gamma
},\mathcal{\zeta)}.$ 

In the second part of the section, we use calculus of variations to look for
$\left(  \mathcal{\gamma}^{0},\zeta^{0}\right)  \in$ $\mathcal{B}$ which
maximizes the value function $W^{\mathcal{\mathcal{\gamma}}%
,\mathcal{\mathcal{\zeta}}}$ among all $\left(  \mathcal{\gamma},\zeta\right)
\in$ $\mathcal{B}.$\newline

Let us consider the set
\[
T:=\left\{  (y,z):0<y\leq z\right\}  ,
\]
and the following auxiliary functions $b_{0}$, $b_{1}:T\times\lbrack
0,\infty)\times\lbrack0,\overline{c}]\rightarrow{\mathbb{R}}$
\begin{equation}%
\begin{array}
[c]{ccc}%
b_{0}(y,z,w,c) & := & \dfrac{1}{q\left(  \theta_{1}(c)-\theta_{2}(c)\right)
d(y,z,c)}(b_{00}(y,z,c)+w(e^{(z-y)\theta_{1}(c)}-e^{(z-y)\theta_{2}%
(c)})~b_{01}(y,c)),\\
b_{1}(y,z,w,c) & := & \dfrac{1}{\left(  \theta_{1}(c)-\theta_{2}(c)\right)
d(y,z,c)}(b_{10}(y,z,c)+w(e^{(z-y)\theta_{1}(c)}-e^{(z-y)\theta_{2}%
(c)})~b_{11}(y,c)),
\end{array}
\label{b0 y b1}%
\end{equation}
where the functions $b_{00},b_{01},b_{10},b_{11}$ and $d$ are defined in
Section \ref{Seccion Formulas}.

\begin{lemma}
\label{Lema Denominador Positivo}The function $d(y,z,c)$ defined in Section
\ref{Seccion Formulas} is positive in $T\times\lbrack0,\overline{c}]$, and so
$b_{0}$ and $b_{1}$ are well defined.
\end{lemma}

\textit{Proof.} Using that $\theta_{1}>0>\theta_{2}$ and $\theta_{1}^{\prime
},\theta_{2}^{\prime}>0$, let us define the function
$g(y,h,c)=d(y,y+h,c)/e^{h\theta_{2}(c)}$. We have that
\[
g(y,0,c)=\left(  e^{y\theta_{1}(ac)}-e^{y\theta_{2}(ac)}\right)  \left(
\theta_{1}(c)-\theta_{2}(c)\right)  >0
\]
for $y>0$ and
\[%
\begin{array}
[c]{lll}%
\partial_{h}g(y,h,c) & = & (\theta_{2}(c)-\theta_{1}(c))e^{h(\theta
_{1}(c)-\theta_{2}(c))}\left(  e^{y\theta_{1}(ac)}(\theta_{2}(c)-\theta
_{1}(ac))+e^{y\theta_{2}(ac)}(\theta_{2}(ac)-\theta_{2}(c))\right) \\
& > & 0
\end{array}
~
\]
for $y\geq0.$ So the result holds.\hfill$\blacksquare$\newline

In order to define $W^{\mathcal{\mathcal{\gamma}},\mathcal{\mathcal{\zeta}}}$
in the non-change regions $\mathcal{NC}_{ac}\mathcal{(\gamma},\mathcal{\zeta
)}$ and $\mathcal{NC}_{c}\mathcal{(\mathcal{\mathcal{\gamma}}}%
,\mathcal{\mathcal{\mathcal{\zeta}})}$, we will define and study in the next
technical lemma the functions $H^{\mathcal{\mathcal{\gamma}}%
,\mathcal{\mathcal{\zeta}}}$ and $A^{\mathcal{\mathcal{\gamma}}%
,\mathcal{\mathcal{\zeta}}}$ .

\begin{lemma}
\label{Funcion de valor Non Change Region} Given $(\gamma
,\mathcal{\mathcal{\zeta}})\in\mathcal{B}$, the unique continuous function
$H^{\gamma,\zeta}:[0,\infty)\times\lbrack0,\overline{c}]\rightarrow
\lbrack0,\infty)$ , with $H^{\gamma,\zeta}(\cdot,c)$ continuously
differentiable which satisfies for any $c\in\lbrack0,\overline{c})$ that
\[
\mathcal{L}^{ac}(H^{\gamma,\zeta})(x,c)=0\text{ for }0\leq x<\gamma
(c)\text{,}~\mathcal{L}^{c}(H^{\gamma,\zeta})(x,c)=0\text{ for }\gamma(c)\leq
x
\]
with boundary conditions $H^{\gamma,\zeta}(0,c)=0$, $H^{\gamma,\zeta
}(x,\overline{c})=v(x,\overline{c},\gamma(\overline{c}))$ and $\partial
_{c}H^{\gamma,\zeta}(\zeta(c),c)=0\ $ at the points of continuity of $\zeta$
is given by%
\begin{equation}
H^{\gamma,\zeta}(x,c)=(f_{10}(x,c)+f_{11}(x,c)A^{\gamma
,\mathcal{\mathcal{\zeta}}}(c))I_{x<\gamma(c)}+(f_{20}(\gamma(c),x,c)+f_{21}%
(\gamma(c),x,c)A^{\gamma,\mathcal{\mathcal{\zeta}}}(c))I_{x\geq\gamma(c)},
\label{Definicion de Hz}%
\end{equation}
where $f_{10},$ $f_{11},f_{20},f_{21}$ are defined in Section
\ref{Seccion Formulas},
\begin{equation}
A^{\mathcal{\mathcal{\gamma}},\mathcal{\mathcal{\zeta}}}%
(c)=A^{\mathcal{\mathcal{\gamma}},\mathcal{\mathcal{\zeta}}}(\overline
{c})e^{-\int_{c}^{\overline{c}}b_{1}(\gamma(s),\zeta(s),\gamma^{\prime
}(s),s)ds}-\int_{c}^{\overline{c}}e^{-\int_{c}^{t}b_{1}(\gamma(s),\zeta
(s),\gamma^{\prime}(s),s)ds}b_{0}(\gamma(t),\zeta(t),\gamma^{\prime}(t),t)dt,
\label{Definicion A(z)}%
\end{equation}
and%
\begin{equation}
A^{\mathcal{\mathcal{\gamma}},\mathcal{\mathcal{\zeta}}}(\overline{c}%
)=\frac{B(\overline{c},\gamma(\overline{c}))}{\sqrt{(\mu-a\overline{c}%
)^{2}+2q\sigma^{2}}}, \label{Definicion A en cbarra}%
\end{equation}
where the function $B$ is defined in (\ref{Definicion B}) and the function
$b_{0}$ and $b_{1}$ are defined in (\ref{b0 y b1}). Moreover,
$A^{\mathcal{\mathcal{\gamma}},\mathcal{\mathcal{\zeta}}}$ is differentiable
and satisfies
\begin{equation}
\left(  A^{\mathcal{\mathcal{\gamma}},\mathcal{\mathcal{\zeta}}}\right)
^{\prime}(c)=b_{0}(\gamma(c),\zeta(c),\gamma^{\prime}(c),c)+b_{1}%
(\gamma(c),\zeta(c),\gamma^{\prime}(c),c)A^{\mathcal{\mathcal{\gamma}%
},\mathcal{\mathcal{\zeta}}}(c), \label{Ecuacion Diferencial de Az}%
\end{equation}
at the points where $\zeta$ is continuous and satisfies the boundary condition
(\ref{Definicion A en cbarra}).
\end{lemma}

\textit{Proof.} Since $H^{\gamma,\zeta}(\cdot,c)$ is continuously
differentiable at $x=\gamma(c)$ and satisfies $H^{\gamma,\zeta}(0,c)=0$,
$\mathcal{L}^{ac}(H^{\gamma,\zeta})(x,c)=0$ for $0\leq x<\gamma(c)$ and
$~\mathcal{L}^{c}(H^{\gamma,\zeta})(x,c)=0$ for $\gamma(c)\leq x,$ we can
write, using (\ref{Solucion General L=0}) and
(\ref{Solucion de L=0 con condicion en 0}),
\[
H^{\gamma,\zeta}(x,c)=(f_{10}(x,c)+f_{11}(x,c)A(c))I_{x<\gamma(c)}%
+(f_{20}(\gamma(c),x,c)+f_{21}(\gamma(c),x,c)A(c))I_{x\geq\gamma(c)}%
\]
for some function $A(c)$. Since
\[
H^{\gamma,\zeta}(x,\overline{c})=v(x,\overline{c},\gamma(\overline{c})),
\]
we obtain, by (\ref{Formula value function refraction en b}), that
\[
A(\overline{c})=\frac{B(\overline{c},\gamma(\overline{c}))}{\sqrt
{(\mu-a\overline{c})^{2}+2q\sigma^{2}}}.
\]

Let us find $A^{\mathcal{\mathcal{\gamma}},\mathcal{\mathcal{\zeta}}%
}:[0,\overline{c}]\rightarrow\mathbb{R}$, the function $A(c)$ such that
$\left.  \partial_{c}H^{\gamma,\zeta}(x,c)\right\vert _{x=\zeta(c)}=0$ for all
$c\in\lbrack0,\overline{c}]$. Since $\zeta(c)>\gamma(c),$
\[%
\begin{array}
[c]{lll}%
0=\left.  \partial_{c}H^{\gamma,\zeta}(x,c)\right\vert _{x=\zeta(c)} & = &
\left.  \frac{d}{dc}(f_{20}(\gamma(c),x,c)+f_{21}(\gamma
(c),x,c)A^{\mathcal{\mathcal{\gamma}},\mathcal{\mathcal{\zeta}}}%
(c))\right\vert _{x=\zeta(c)}\\
& = & \left.  \frac{d}{dc}(f_{20}(\gamma(c),x,c))+\frac{d}{dc}\left(
f_{21}(\gamma(c),x,c)\right)  A^{\mathcal{\mathcal{\gamma}}%
,\mathcal{\mathcal{\zeta}}}(c))+f_{21}(\gamma
(c),x,c)A^{\mathcal{\mathcal{\gamma}},\mathcal{\mathcal{\zeta}}\prime
}(c)\right\vert _{x=\zeta(c)},
\end{array}
\]
and so, since by Lemma \ref{Lema Denominador Positivo}, $f_{21}(y,x,c)=\frac
{d(y,x,c)}{\text{$\theta_{1}(c)$}-\text{$\theta_{2}(c)$}}>0$ for $x>y,$ we obtain%

\[%
\begin{array}
[c]{ccc}%
(A^{\mathcal{\mathcal{\gamma}},\mathcal{\mathcal{\zeta}}})^{\prime}(c) & = &
\left.  \frac{-\frac{d}{dc}(f_{20}(\gamma(c),x,c))}{f_{21}(\gamma
(c),x,c)}\right\vert _{x=\zeta(c)}+\left.  \frac{-\frac{d}{dc}(f_{21}%
(\gamma(c),x,c))}{f_{21}(\gamma(c),x,c)}\right\vert _{x=\zeta(c)}%
A^{\mathcal{\mathcal{\gamma}},\mathcal{\mathcal{\zeta}}}(c)\\
& = & b_{0}(\gamma(c),\zeta(c),\gamma^{\prime}(c),c)+b_{1}(\gamma
(c),\zeta(c),\gamma^{\prime}(c),c)A^{\mathcal{\mathcal{\gamma}}%
,\mathcal{\mathcal{\zeta}}}(c),
\end{array}
~
\]
at the points where $\zeta$ is continuous, where $b_{0}$ and $b_{1}$ are
defined in (\ref{b0 y b1}). Since $\zeta$ is Riemann integrable, it is
differentiable almost everywhere. Note that the function
$A^{\mathcal{\mathcal{\gamma}},\mathcal{\mathcal{\zeta}}}$ defined in
(\ref{Definicion A(z)}) is the unique solution of this ODE. Hence, we have the
result. \hfill$\blacksquare$\newline

Given $(\gamma,\zeta)\in\mathcal{B}$, we define
\begin{equation}
W^{\gamma,\zeta}(x,c):=\left\{
\begin{array}
[c]{lll}%
H^{\gamma,\zeta}(x,c)\text{ } & \text{if} & (x,c)\in\mathcal{NC(\zeta)},\\
H^{\gamma,\zeta}(x,\ell(x,c)) & \text{if} & (x,c)\in\mathcal{CH(\zeta)},
\end{array}
\right.  \label{Definicion de Wz}%
\end{equation}
where $H^{\gamma,\zeta}$ is defined in Lemma
\ref{Funcion de valor Non Change Region} and%
\begin{equation}
\ell(x,c):=\max\{h\in\lbrack c,\overline{c}]:\zeta(d)\leq x{\large \ }%
\text{{\large for}}{\large \ }d\in\lbrack c,h)\} \label{Definicion de c(x)}%
\end{equation}
for $x\geq\zeta(c)$ and $c\in\lbrack0,\overline{c})$.\newline

In the next propositions we will show first that in the case that $\zeta$ is a
step function, $W^{\gamma,\zeta}$ is the value function of a two-curve
strategy and in the general case $W^{\gamma,\zeta}$ is the limit of value
functions of two-curve strategies.

When $\zeta$ is a step function, that is%
\[
\zeta(c):=%
{\displaystyle\sum\limits_{i=1}^{n-1}}
{z}_{i}{I}_{[c_{i},c_{i+1})}{(c),}%
\]
with $0=c_{1}<c_{2}<\cdots<c_{n}=\overline{c}$ and $z_{i}>0$, then the
two-curve strategy, starting with an initial surplus $x$ and initial running
maximum dividend rate $c,$ is given by

(1) if $0\leq x<\zeta(c)$, that is $(x,c)\in
\mathcal{NC(\mathcal{\mathcal{\gamma}}},\mathcal{\mathcal{\mathcal{\zeta}})},$
follow the refracting strategy which pays $ac$ when the current surplus is
below a refracting threshold $\gamma(c)$ and pays $c$ when the current surplus
is above $\gamma(c)$ until either reaching the curve $\mathcal{R(\zeta)}$ or
ruin (whatever comes first),

(2) If $x\geq\zeta(c)$, that is $(x,c)\in\mathcal{CH(\mathcal{\mathcal{\gamma
}}},\mathcal{\mathcal{\mathcal{\zeta}})},$ increase immediately the dividend
rate to $\ell(x,c)\in\{c_{2},\cdots,c_{n}\};$ note that
\[
\ell(x,c)=\max\{c_{i}\geq c:z_{k}\leq x~\text{for }l(c)\leq k\leq
i-1\},\ \text{where }l(c):=\min\{l:c_{l}\geq c\}.
\]

If $\zeta$ is a step function, we denote this stationary strategy as
$\mathbf{\pi}^{(\gamma,\zeta)}.$

\begin{proposition}
\label{Extended threshold strategy} Consider $(\gamma,\zeta)\in\mathcal{B}$
with $\zeta$ being a step function. Let $D^{x,c}\in\Pi_{x,c}^{[0,\overline
{c}]}$ be the admissible strategy corresponding to the stationary strategy
$\mathbf{\pi}^{(\gamma,\zeta)}$ starting in $(x,c)$. Calling
$j(x,c):=J(x;D^{x,c})$, we obtain that $j$ is continuous in $[0,\infty
)\times\lbrack0,\overline{c}]$ and $j(x,c)=W^{\gamma,\zeta}(x,c).$
\end{proposition}

\textit{Proof.} Let us prove inductively that $j(x,c_{i})\ $is continuous in
$x$ for $i=1,\ldots,n$. $j(\cdot,c_{i})$ is differentiable in $[0,z_{i})$
because it corresponds to a value function of a refracting dividend strategy
at $x=\gamma(c_{i})$ with a given boundary condition at $x=z_{i}$ (see for
instance \cite[Th.3.1]{ABB}).\ In the case $i=n$, $j(x,c_{n})=v(x,c,\gamma
(\overline{c}))$ which is continuous in $x$; in the case $i<n,$ $j(x,c_{i}%
)\ $is continuous in $x$ for $x\leq z_{i}$ because
\[
j(x,c_{i})=(f_{10}(x,c_{i})+f_{11}(x,c_{i})A_{i})I_{x<\gamma(c_{i})}%
+(f_{20}(\gamma(c_{i}),x,c_{i})+f_{21}(\gamma(c_{i}),x,c)A_{i})I_{x\geq
\gamma(c_{i})}%
\]
for some constant $A_{i}$, and $j(x,c_{i})=j(x,c_{i+1})$ for $x\geq z_{i}.$
Since $j(x,c)=j(x,c_{i+1})$ for {$c\in(c_{i},c_{i+1})$, we conclude that }$j$
is continuous in $[0,\infty)\times\lbrack0,\overline{c}]$ $.$

Let us show now that $j(x,c)$ satisfies the assumptions of Lemma
\ref{Funcion de valor Non Change Region} and so $j(x,c)=H^{\gamma,\zeta
}(x,c)=W^{\gamma,\zeta}(x,c)$ for $0\leq x\leq\zeta(c)$. Indeed, it {is
straightforward that }$j(\cdot,\overline{c})${$=H^{\gamma,\zeta}%
(\cdot,\overline{c})=$}$v(\cdot,c,\gamma(\overline{c}))${, }$j(\cdot,c)$ is
continuously differentiable for any $c\in\lbrack0,\overline{c})${,
}$\mathcal{L}^{ac}(j)(x,c)=0$ for $0\leq x<\gamma(c)$,$~\mathcal{L}%
^{c}(j)(x,c)=0$ for $\gamma(c)\leq x\leq\zeta(c)$ and $j(0,c)=0$. Also
{$\partial_{c}j(\zeta(c),c)=0$ }at the points of continuity of $\zeta
\ $because $j(x,c)${$=j(x,c_{i+1}{)}$ for $x\geq\zeta(c)=z$}$_{{i}}$ i{n the
case $c\in(c_{i},c_{i+1})$.}

From definition of $\mathbf{\pi}^{(\gamma,\zeta)},$ it is straightforward that
$j(x,c)=H^{\gamma,\zeta}(x,\ell(x,c))$ if $x\geq\zeta(c)$, so we get the
result. \hfill$\blacksquare$\newline

In the next proposition we show that for any $(\gamma,\zeta)\in\mathcal{B}$,
the function $W^{\gamma,\zeta}$ is the limit of value functions of curve
strategies where $\zeta_{k}$ are step functions with $\zeta_{k}\rightarrow
\zeta$ uniformly.

\begin{proposition}
\label{Limite de value functions} Given $(\gamma,\zeta)\in\mathcal{B}$, there
exists a sequence of right-continuous step functions $\zeta_{k}:[0,\overline
{c}]\rightarrow\lbrack0,\infty)$ such that $W^{\gamma,\zeta_{k}}(x,c)$
converges uniformly to $W^{\gamma,\zeta}(x,c)$.
\end{proposition}

\textit{Proof.} Since $\zeta$ is a Riemann integrable c\`{a}dl\`{a}g function,
we can approximate it uniformly by right-continuous step functions. Namely,
take a sequence of finite sets $\mathcal{S}^{k}=\{c_{1}^{k},c_{2}^{k}%
,\cdots,c_{n_{k}}^{k}\}$ with $0=c_{1}^{k}<c_{2}^{k}<\cdots<c_{n_{k}}%
^{k}=\overline{c}$, and consider the right-continuous step functions%
\[
\zeta_{k}(c)=%
{\displaystyle\sum\limits_{i=1}^{n_{k}-1}}
\zeta(c_{i}^{k})I_{[c_{i}^{k},c_{i+1}^{k})},
\]
such that $\delta(\mathcal{S}^{k})=\max_{i=1,\cdots,n_{k}-1}(c_{i+1}^{k}%
-c_{i}^{k})\rightarrow0$. We have that $\zeta_{k}\rightarrow\zeta$ uniformly,
and so both $A^{\gamma,\zeta_{k}}(c)\rightarrow A^{\gamma,\zeta}(c)$ and
$W^{\gamma,\zeta_{k}}(x,c)$ $\rightarrow$ $W^{\gamma,\zeta}(x,c)$ uniformly.
\hfill$\blacksquare$\newline

\begin{remark}
\normalfont
\label{Remark general two-cuves strategy} Given a $(\gamma,\zeta
)\in\mathcal{B}$ where $\zeta$ is not a step function, we say that
$W^{\gamma,\zeta}$ is the value function of the \textit{two-curve} stationary
strategy $\mathbf{\pi}^{(\gamma,\zeta)}$ which, starting with an initial
surplus $x$ and initial running maximum dividend rate $c$,

(1) in the case $0\leq x<\zeta(c),$ it follows the refracting strategy which
pays $ac$ when the current surplus is below a refracting threshold $\gamma(c)$
and pays $c$ when the current surplus is above $\gamma(c)$ until either
reaching the curve $\mathcal{R(\zeta)}$ or ruin (whatever comes first),

(2) in the case $x>\zeta(c),$ increase immediately the divided rate from $c$
to $\ell(x,c)$,

(3) in the case $x=\zeta(c)$, it can be seen as the limit of admissible
strategies in $\mathbf{\pi}_{x,c}^{(\gamma,\zeta_{k})}\in\Pi_{x,c}$ arising
from Proposition \ref{Limite de value functions}.
\end{remark}

We now look for the maximum of $W^{\gamma,\zeta}$ among $(\gamma,\zeta
)\in\mathcal{B}$. We will show later that, if there exists a pair $(\gamma
_{0},\zeta_{0})\in\mathcal{B}$ such that%
\begin{equation}
A^{\gamma_{0},\zeta_{0}}(0)=\max\{A^{\gamma,\zeta}(0):(\gamma,\zeta
)\in\mathcal{B}\}, \label{Definicion z0}%
\end{equation}
then $W^{\gamma_{0},\zeta_{0}}(x,c)\geq W^{\gamma,\zeta}(x,c)$ for all
$(x,c)\in\lbrack0,\infty)\times\lbrack0,\overline{c}]\ $and $(\gamma,\zeta
)\in\mathcal{B}.$

From Lemma \ref{Lema Denominador Positivo} and $\theta_{2}>0>\theta_{1}$, we
obtain that $f_{11}$ and $f_{21}$ defined in (\ref{Definicion de Hz}) are
positive and so
\[
\arg\max_{(\gamma,\zeta)\in\mathcal{B}}W^{\gamma,\zeta}(x,c)=\arg\max
_{(\gamma,\zeta)\in\mathcal{B}}A^{\gamma,\zeta}(c).
\]
This follows from (\ref{Definicion de Hz}) and the next lemma, in which we
prove that the function $\zeta_{0}$ which maximizes (\ref{Definicion z0}) also
maximizes $A^{\gamma,\zeta}(c)$ for any $c\in\lbrack0,\overline{c})$.

\begin{lemma}
\label{Lema cortar en c} For a given $c\in\lbrack0,\overline{c})$, consider
the functions $\gamma:[c,\overline{c}]\rightarrow(0,\infty)~$%
continuously~differentiable,$~\zeta:[c,\overline{c}]\rightarrow(0,\infty)~$is
bounded, Riemann integrable and c\`{a}dl\`{a}g, and let us define the set%
\[
\mathcal{B}_{c}=\{(\gamma,\zeta)~\text{s.t. }\gamma\leq\zeta\text{ in
}[c,\overline{c}]~\text{and }\lim_{c\rightarrow\overline{c}^{-}}\zeta
(c)=\zeta(\overline{c}^{-})\}.
\]
If $(\gamma_{0},\zeta_{0})\in\mathcal{B}$ satisfies (\ref{Definicion z0}),
then for any $c\in\lbrack0,\overline{c})$
\[
A^{\gamma_{0},\zeta_{0}}(c)=\max\{A^{\gamma,\zeta}(c):\zeta\in\mathcal{B}%
_{c}\}.\
\]

\end{lemma}

\textit{Proof.} Given $(\gamma,\zeta)~\in\mathcal{B},$ we can write%
\[
A^{\mathcal{\mathcal{\gamma}},\mathcal{\mathcal{\zeta}}}%
(c)=A^{\mathcal{\mathcal{\gamma}},\mathcal{\mathcal{\zeta}}}(\overline
{c})e^{-\int_{c}^{\overline{c}}b_{1}(\gamma(s),\zeta(s),\gamma^{\prime
}(s),s)ds}-\int_{c}^{\overline{c}}e^{-\int_{c}^{t}b_{1}(\gamma(s),\zeta
(s),\gamma^{\prime}(s),s)ds}b_{0}(\gamma(t),\zeta(t),\gamma^{\prime}(t),t)dt,
\]%
\[
A^{\gamma,\zeta}(0)=-\int_{0}^{c}e^{-\int_{0}^{t}b_{1}(\gamma(s),\zeta
(s),\gamma^{\prime}(s),s)ds}b_{0}(\gamma(t),\zeta(t),\gamma^{\prime
}(t),t)dt+\left(  e^{-\int_{0}^{c}b_{1}(\gamma(s),\zeta(s),\gamma^{\prime
}(s),s)ds}\right)  A^{\gamma,\zeta}(c).
\]
{ Note that }$\left(  {\left.  \gamma_{0}\right\vert _{[c,\overline{c}%
)},{\left.  \zeta_{0}\right\vert _{[c,\overline{c})}}}\right)  $%
{$\in\mathcal{B}_{c}$; take any }$(\gamma,\zeta)~${$\in\mathcal{B}_{c}$ and
any function }$\chi:[0,\overline{c}]\rightarrow\lbrack0,1]$ continuously
differentiable with{ }$\chi=0$ in $[0,c]$ and define {{ }}${{\gamma}}$%
{{$_{1}(s)={{\gamma}}_{0}(s)I_{\{0\leq s<c\}}+{{\gamma}}(s)\chi(s)I_{\{c\leq
s\leq\overline{c}\}}$} and {$\zeta_{1}(s)=\zeta_{0}(s)I_{\{0\leq s<c\}}%
+\zeta(s)\chi(s)I_{\{c\leq s\leq\overline{c}\}}$ }then }${({\gamma}_{1}%
,{\zeta_{1}})}${{$\in\mathcal{B}$}},{%
\[
A^{\gamma_{0},\zeta_{0}}(0)\geq A^{\mathcal{\mathcal{\gamma}}_{1}\zeta_{1}%
}(0)=-\int_{0}^{c}e^{-\int_{0}^{t}b_{1}(\gamma_{0}(s),\zeta_{0}(s),s)ds}%
b_{0}(\gamma_{0}(t),\zeta_{0}(t),t)dt+\left(  e^{-\int_{0}^{c}b_{1}(\gamma
_{0}(s),\zeta_{0}(s),s)ds}\right)  A^{{{\gamma}}\chi,\zeta\chi}(c).
\]
Hence,%
\[%
\begin{array}
[c]{lll}%
A^{\gamma_{0},\zeta_{0}}(0) & \geq & -\int_{0}^{c}e^{-\int_{0}^{t}b_{1}%
(\gamma_{0}(s),\zeta_{0}(s),s)ds}b_{0}(\gamma_{0}(t),\zeta_{0}(t),t)dt+\left(
e^{-\int_{0}^{c}b_{1}(\gamma_{0}(s),\zeta_{0}(s),s)ds}\right)  \sup
_{(\gamma,\zeta)\in\mathcal{B}_{c},\chi}A^{{{\gamma}}\chi,\zeta\chi}(c)\text{
}\\
& = & -\int_{0}^{c}e^{-\int_{0}^{t}b_{1}(\gamma_{0}(s),\zeta_{0}(s),s)ds}%
b_{0}(\gamma_{0}(t),\zeta_{0}(t),t)dt+\left(  e^{-\int_{0}^{c}b_{1}(\gamma
_{0}(s),\zeta_{0}(s),s)ds}\right)  \sup_{(\gamma,\zeta)\in\mathcal{B}_{c}%
}A^{{{\gamma}},\zeta}(c)\\
& \geq & -\int_{0}^{c}e^{-\int_{0}^{t}b_{1}(\gamma_{0}(s),\zeta_{0}%
(s),s)ds}b_{0}(\gamma_{0}(t),\zeta_{0}(t),t)dt+\left(  e^{-\int_{0}^{c}%
b_{1}(\gamma_{0}(s),\zeta_{0}(s),s)ds}\right)  A^{\mathcal{\mathcal{\gamma}%
}_{0}\zeta_{0}}(c)\\
& = & A^{\mathcal{\mathcal{\gamma}}_{0}\zeta_{0}}(0),
\end{array}
\]
and so we have that }$\sup_{(\gamma,\zeta)\in\mathcal{B}_{c}}A^{{{\gamma}%
},\zeta}(c)=A^{\mathcal{\mathcal{\gamma}}_{0},\zeta_{0}}(c)$.\hfill
$\blacksquare$\newline

\noindent Let us now find the implicit equation for the function
$A^{\gamma_{0},\zeta_{0}}$ for $(\gamma_{0},\zeta_{0})$ satisfying
(\ref{Definicion z0}).

\begin{proposition}
\label{Condicion de Optimo} If the pair $(\gamma_{0},\zeta_{0})$ defined in
(\ref{Definicion z0}) exists, then $A^{\gamma_{0},\zeta_{0}}(c)$ satisfies%
\begin{equation}
b_{1z}(c)A^{\gamma_{0},\zeta_{0}}(c)+b_{0z}(c)=0~\text{for any }c\in
\lbrack0,\overline{c}), \label{First equation}%
\end{equation}%
\begin{equation}
b_{1w}(c)A^{\gamma_{0},\zeta_{0}}(c)+b_{0w}(c)=0~~\text{for all }c\in
\lbrack0,\overline{c}] \label{Second equation}%
\end{equation}
and%
\begin{equation}
b_{1y}(c)A^{\gamma_{0},\zeta_{0}}(c)+b_{0y}(c)=0~~\text{for all }c\in
\lbrack0,\overline{c}], \label{Third equation}%
\end{equation}
where%
\[%
\begin{array}
[c]{ll}%
b_{i}(s):=b_{i}(\gamma_{0}(s),\zeta_{0}(s),\gamma_{0}^{\prime}(s),s)\text{,} &
\text{ }b_{iy}(s):=\partial_{y}b_{i}(\gamma_{0}(s),\zeta_{0}(s),\gamma
_{0}^{\prime}(s),s)\text{,}\\
b_{iz}(s):=\partial_{z}b_{i}(\gamma_{0}(s),\zeta_{0}(s),\gamma_{0}^{\prime
}(s),s) & \text{and~}b_{iw}(s):=\partial_{w}b_{i}(\gamma_{0}(s),\zeta
_{0}(s),\gamma_{0}^{\prime}(s),s)
\end{array}
\]
$\ $for $i=0,1$.

Moreover, $\gamma_{0}(\overline{c})=b^{\ast}(\overline{c})$ is the optimal
threshold defined in \eqref{Ecuacion bestrella} and
\begin{equation}
A^{\gamma_{0},\zeta_{0}}(\overline{c})=\frac{B(\overline{c},b^{\ast}%
(\overline{c}))}{\sqrt{(\mu-ac)^{2}+2q\sigma^{2}}}\text{.}
\label{A en la tapa}%
\end{equation}

\end{proposition}

\textit{Proof.} Consider any function $(\gamma_{1},\zeta_{1})\in\mathcal{B}$
with $\gamma_{1}(\overline{c})=\zeta_{1}(\overline{c})=0$. Then%
\[%
\begin{array}
[c]{l}%
A^{\gamma_{0}+\eta\gamma_{1},\zeta_{0}+\varepsilon\zeta_{1}}(0)\\%
\begin{array}
[c]{ll}%
= & A^{\gamma_{0},\zeta_{0}}(\overline{c})e^{-\int_{0}^{\overline{c}}%
b_{1}(\gamma_{0}(t)+\eta\gamma_{1}(t),\zeta_{0}(t)+\varepsilon\zeta
_{1}(t),\gamma_{0}^{\prime}(t)+\eta\gamma_{1}^{\prime}(t),t)dt}\\
& -\int_{0}^{\overline{c}}e^{-\int_{0}^{s}b_{1}(\gamma_{0}(u)+\eta\gamma
_{1}(u),\zeta_{0}(u)+\varepsilon\zeta_{1}(u),\gamma_{0}^{\prime}(u)+\eta
\gamma_{1}^{\prime}(u),u)du}b_{0}(\gamma_{0}(s)+\eta\gamma_{1}(s),\zeta
_{0}(s)+\varepsilon\zeta_{1}(s),\gamma_{0}^{\prime}(s)+\eta\gamma_{1}^{\prime
}(s),s)ds\text{.}%
\end{array}
\end{array}
\]
Taking the derivative with respect to $\varepsilon$ at $\eta=\varepsilon=0$,
we get%

\[%
\begin{array}
[c]{lll}%
0=\left.  \partial_{\varepsilon}A^{\gamma_{0}+\eta\gamma_{1},\zeta
_{0}+\varepsilon\zeta_{1}}(0)\right\vert _{\eta=0,\varepsilon=0} & = &
\int_{0}^{\overline{c}}\left(  e^{-\int_{0}^{s}b_{1}(u)du}b_{0}(s)\left(
\int_{0}^{s}\zeta_{1}(u)b_{1z}(u)du\right)  \right)  ds\\
&  & -\int_{0}^{\overline{c}}(e^{-\int_{0}^{s}b_{1}(u)du}\zeta_{1}%
(s)b_{0z}(s))ds\\
&  & -A^{\gamma_{0},\zeta_{0}}(\overline{c})e^{-\int_{0}^{\overline{c}}%
b_{1}(t)dt}\int_{0}^{\overline{c}}\zeta_{1}(s)b_{0z}(s))ds.
\end{array}
\]
Using integration by parts we obtain%
\[
\int_{0}^{\overline{c}}\left(  e^{-\int_{0}^{s}b_{1}(u)du}b_{0}(s)\left(
\int_{0}^{s}\zeta_{1}(u)b_{1z}(u)du\right)  \right)  ds=\int_{0}^{\overline
{c}}\left(  \zeta_{1}(s)b_{1z}(s)\left(  \int_{s}^{\overline{c}}e^{-\int
_{0}^{t}b_{1}(u)du}b_{0}(t)dt\right)  \right)  ds,
\]
and so%
\[
0=\int_{0}^{\overline{c}}\zeta_{1}(s)\left(  b_{1z}(s)\int_{s}^{\overline{c}%
}e^{-\int_{0}^{t}b_{1}(u)du}b_{0}(t)dt-e^{-\int_{0}^{s}b_{1}(u)du}%
b_{0z}(s)-A^{\gamma_{0},\zeta_{0}}(\overline{c})e^{-\int_{0}^{\overline{c}%
}b_{1}(t)dt}b_{0z}(s)\right)  ds.
\]
Since this holds for any $\zeta_{1}$ with $\zeta_{1}(\overline{c})=0$, we get
using (\ref{Definicion A(z)}) that for any $c\in\lbrack0,\overline{c})$%

\[%
\begin{array}
[c]{lll}%
0 & = & b_{1z}(c)\int_{c}^{\overline{c}}e^{-\int_{0}^{t}b_{1}(u)du}%
b_{0}(t)dt-e^{-\int_{0}^{c}b_{1}(u)du}b_{0z}(c)-A^{\gamma_{0},\zeta_{0}%
}(\overline{c})e^{-\int_{0}^{\overline{c}}b_{1}(t)dt}b_{0z}(c)\\
& = & e^{-\int_{0}^{c}b_{1}(u)du}(-b_{1z}(c)A^{\gamma_{0},\zeta_{0}}%
(c)-b_{0z}(c)),
\end{array}
\]
and so we conclude (\ref{First equation}).

Taking the derivative with respect to $\eta$ at $\eta=\varepsilon=0$, we get%

\[%
\begin{array}
[c]{lll}%
0 & = & \left.  \partial_{\eta}A^{\gamma_{0}+\eta\gamma_{1},\zeta
_{0}+\varepsilon\zeta_{1}}(0)\right\vert _{\eta=0,\varepsilon=0}\\
& = & \int_{0}^{\overline{c}}\left(  e^{-\int_{0}^{s}b_{1}(u)du}%
b_{0}(s)\left(  \int_{0}^{s}\gamma_{1}^{\prime}(u)b_{1w}(u)du\right)  \right)
ds\\
&  & +\int_{0}^{\overline{c}}\left(  e^{-\int_{0}^{s}b_{1}(u)du}%
b_{0}(s)\left(  \int_{0}^{s}\gamma_{1}(u)b_{1y}(u)du\right)  \right)  ds\\
&  & -\int_{0}^{\overline{c}}\left(  e^{-\int_{0}^{s}b_{1}(u)du}\gamma
_{1}^{\prime}(s)b_{0w}(s)\right)  ds-\int_{0}^{\overline{c}}\left(
e^{-\int_{0}^{s}b_{1}(u)du}\gamma_{1}(s)b_{0y}(s)\right)  ds\\
&  & -A^{\gamma_{0},\zeta_{0}}(\overline{c})e^{-\int_{0}^{\overline{c}}%
b_{1}(t)dt}\int_{0}^{\overline{c}}\gamma_{1}^{\prime}(t)b_{1z}(t))dt-A^{\gamma
_{0},\zeta_{0}}(\overline{c})e^{-\int_{0}^{\overline{c}}b_{1}(t)dt}\int
_{0}^{\overline{c}}\gamma_{1}(t)b_{1y}(t))ds.
\end{array}
\]
Using integration by parts, we obtain%
\[%
\begin{array}
[c]{lll}%
0 & = & \gamma_{1}(0)\left(  -b_{1w}(0)\int_{0}^{\overline{c}}e^{-\int_{0}%
^{t}b_{1}(u)du}b_{0}(t)dt+b_{0w}(0)+e^{-\int_{0}^{\overline{c}}b_{1}%
(t)dt}A^{\gamma_{0},\zeta_{0}}(\overline{c})b_{1w}(0)\right) \\
&  & +\gamma_{1}(\overline{c})\left(  -e^{-\int_{0}^{\overline{c}}b_{1}%
(t)dt}b_{0w}(\overline{c})-e^{-\int_{0}^{\overline{c}}b_{1}(t)dt}A^{\gamma
_{0},\zeta_{0}}(\overline{c})b_{1w}(\overline{c})\right) \\
&  & +\int_{0}^{\overline{c}}\gamma_{1}(s)\left(  \frac{d}{ds}b_{1w}(s)\left(
e^{-\int_{s}^{\overline{c}}b_{1}(t)dt}A^{\gamma_{0},\zeta_{0}}(\overline
{c})-\int_{s}^{\overline{c}}e^{-\int_{s}^{t}b_{1}(u)du}b_{0}(t)dt\right)
\right)  ds\\
&  & +\int_{0}^{\overline{c}}\gamma_{1}(s)\left(  b_{1y}(s)\left(  \int
_{s}^{\overline{c}}e^{-\int_{s}^{t}b_{1}(u)du}b_{0}(t)dt-e^{-\int
_{s}^{\overline{c}}b_{1}(t)dt}A^{\gamma_{0},\zeta_{0}}(\overline{c})\right)
\right)  ds\\
&  & +\int_{0}^{\overline{c}}\gamma_{1}(s)\left(  b_{1w}(s)b_{0}(s)+\frac
{d}{ds}b_{1w}(s)-b_{0w}(s)b_{1}(s)-b_{0y}(s)\right)  ds.
\end{array}
\]
Using that $\gamma_{1}(\overline{c})=0$ and (\ref{Definicion A(z)}), we get%
\[%
\begin{array}
[c]{lll}%
0 & = & \gamma_{1}(0)\left(  b_{1w}(0)A^{\gamma_{0},\zeta_{0}}(0)+b_{0w}%
(0)\right) \\
&  & +\int_{0}^{\overline{c}}\gamma_{1}(s)e^{-\int_{0}^{s}b_{1}(u)du}\left(
\left(  \frac{d}{ds}b_{1w}(s)-b_{1y}(s)\right)  A^{\gamma_{0},\zeta_{0}%
}(s)+\frac{d}{ds}b_{0w}(s)-b_{0y}(s)+b_{1w}(s)b_{0}(s)-b_{0w}(s)b_{1}%
(s)\right)  ds.
\end{array}
\]
Since this holds for any $\gamma_{1}$ with $\gamma_{1}(\overline{c})=0$, we
obtain%
\begin{equation}
\left(  \frac{d}{ds}b_{1w}(s)-b_{1y}(s)\right)  A^{\gamma_{0},\zeta_{0}%
}(s)+\frac{d}{ds}b_{0w}(s)-b_{0y}(s)+b_{1w}(s)b_{0}(s)-b_{0w}(s)b_{1}%
(s)=0~\text{for all }c\in\lbrack0,\overline{c}] \label{ecuacion irrelevante}%
\end{equation}
and
\[
b_{1w}(0)A^{\gamma_{0},\zeta_{0}}(0)+b_{0w}(0)=0.
\]
By Lemma \ref{Lema cortar en c}, we also obtain, taking the derivative
$0=\left.  \partial_{\eta}A^{\gamma_{0}+\eta\gamma_{1},\zeta_{0}%
+\varepsilon\zeta_{1}}(c)\right\vert _{\eta=0,\varepsilon=0}$, that
(\ref{Second equation}) holds. Note that with (\ref{Second equation}), we have%

\[%
\begin{array}
[c]{lll}%
0 & = & \frac{d}{ds}\left(  b_{1w}(s)A^{\gamma_{0},\zeta_{0}}(s)+b_{0w}%
(s)\right) \\
& = & \frac{d}{ds}b_{1w}(s)A^{\gamma_{0},\zeta_{0}}(s)+\frac{d}{ds}%
b_{0w}(s)+b_{0}(s)b_{1w}(s)-b_{1}(s)b_{0w}(s)
\end{array}
\]

\bigskip

and so, from (\ref{ecuacion irrelevante})
\[%
\begin{array}
[c]{lll}%
0 & = & \frac{d}{ds}b_{1w}(s)A^{\gamma_{0},\zeta_{0}}(s)-b_{1y}(s)A^{\gamma
_{0},\zeta_{0}}(s)+\frac{d}{ds}b_{0w}(s)-b_{0y}(s)+b_{1w}(s)b_{0}%
(s)-b_{0w}(s)b_{1}(s)\\
& = & \left[  \frac{d}{ds}b_{1w}(s)A^{\gamma_{0},\zeta_{0}}(s)+\frac{d}%
{ds}b_{0w}(s)+b_{1w}(s)b_{0}(s)-b_{0w}(s)b_{1}(s)\right]  -b_{1y}%
(s)A^{\gamma_{0},\zeta_{0}}(s)-b_{0y}(s)\\
& = & \frac{d}{ds}\left(  b_{1w}(s)A^{\gamma_{0},\zeta_{0}}(s)+b_{0w}%
(s)\right)  -(b_{1y}(s)A^{\gamma_{0},\zeta_{0}}(s)+b_{0y}(s))\\
& = & -(b_{1y}(s)A^{\gamma_{0},\zeta_{0}}(s)+b_{0y}(s)),
\end{array}
\]
from which we conclude (\ref{Third equation}). \hfill$\blacksquare$\newline

\begin{proposition}
\label{condiciones para zo derivable} Consider the functions $C_{0}$ and
$C_{ij}$ for $i=1,2$ and $j=0,1,2$ defined in Section \ref{Seccion Formulas}.
If $(\gamma_{0},\zeta_{0})\in\mathcal{B}\ $defined in (\ref{Definicion z0})
satisfies that $\zeta_{0}$ {is continuous} and
\begin{equation}
C_{11}(\gamma_{0}(c),\zeta_{0}(c),c)\cdot C_{22}(\gamma_{0}(c),\zeta
_{0}(c),c)\neq0 \label{Condicion 2 para z derivable}%
\end{equation}
for $c\in\lbrack0,\overline{c}]$, then $\gamma_{0}\ $and $\zeta_{0}$ are
infinitely differentiable and $(\gamma_{0},\zeta_{0})$ is a solution of the
system of ODE's%
\begin{equation}
\left\{
\begin{array}
[c]{lll}%
\gamma^{\prime}(c) & = & \dfrac{C_{10}(\gamma(c),\zeta(c),c)}{C_{11}%
(\gamma(c),\zeta(c),c)}\\
\zeta^{\prime}(c) & = & \dfrac{C_{20}(\gamma(c),\zeta(c),c)C_{11}%
(\gamma(c),\zeta(c),c)-C_{21}(\gamma(c),\zeta(c),c)C_{10}(\gamma
(c),\zeta(c),c)}{C_{11}(\gamma(c),\zeta(c),c)C_{22}(\gamma(c),\zeta(c),c)}%
\end{array}
\right.  \label{Ecuacion diferencial de z0}%
\end{equation}
with boundary conditions
\begin{equation}
\gamma_{0}(\overline{c})=b^{\ast}(\overline{c})~\text{and~}C_{0}(b^{\ast
}(\overline{c}),\zeta_{0}(\overline{c}),\overline{c}))=0,
\label{Condicion de Optimo en cbarra}%
\end{equation}
where $b^{\ast}(\overline{c})$ is the optimal threshold defined in \eqref{Ecuacion bestrella}.
\end{proposition}

\textit{Proof.} From (\ref{First equation}) and (\ref{Second equation}), we
have%
\[
(b_{1w}b_{0z}-b_{0w}b_{1z})(\gamma_{0}(c),\zeta_{0}(c),\gamma_{0}^{\prime
}(c),c)=0\text{~for }c\in\lbrack0,\overline{c}].
\]
From (\ref{b0 y b1}) we can write
\[%
\begin{array}
[c]{lll}%
(b_{1w}b_{0z}-b_{0w}b_{1z})(y,z,w,c) & = & \tfrac{(e^{(z-y)\theta_{1}%
(c)}-e^{(z-y)\theta_{2}(c)})~}{q\left(  \theta_{1}(c)-\theta_{2}(c)\right)
^{2}d(y,z,c)}\left(  b_{11}(y,c)\partial_{z}\left(  \frac{b_{00}%
(y,z,c)}{d(y,z,c)}\right)  -b_{01}(y,c)\partial_{z}\left(  \frac
{b_{10}(y,z,c)}{d(y,z,c)}\right)  \right) \\
& = & \tfrac{(e^{(z-y)\theta_{1}(c)}-e^{(z-y)\theta_{2}(c)})~}{q\left(
\theta_{1}(c)-\theta_{2}(c)\right)  ^{2}d(y,z,c)}C_{0}(y,z,c),
\end{array}
\]
which does not depend on $w$. So we conclude, that
\begin{equation}
C_{0}(\gamma_{0}(c),\zeta_{0}(c),c)=0\text{~for }c\in\lbrack0,\overline{c}].
\label{Ecuacion con w y z}%
\end{equation}
Moreover, from Proposition \ref{Condicion de Optimo}, we get
(\ref{Condicion de Optimo en cbarra}).

From (\ref{Second equation}) and (\ref{Third equation}), we have%
\[
(b_{1w}b_{0y}-b_{0w}b_{1y})(\gamma_{0}(c),\zeta_{0}(c),\gamma_{0}^{\prime
}(c),c)=0\text{~for }c\in\lbrack0,\overline{c}],
\]
and we can write from (\ref{b0 y b1}),
\[
(b_{1w}b_{0y}-b_{0w}b_{1y})(y,z,w,c)=\frac{(e^{(z-y)\theta_{1}(c)}%
-e^{(z-y)\theta_{2}(c)})}{q\left(  \theta_{1}(c)-\theta_{2}(c)\right)
^{2}d(y,z,c)}(wC_{11}(y,z,c)-C_{10}(y,z,c)).
\]
So, since $C_{11}(\gamma_{0}(c),\zeta_{0}(c),c)\neq0,$ we get the first
equation of (\ref{Ecuacion diferencial de z0}).

Taking the derivative of (\ref{Ecuacion con w y z}) with respect to $c$ and
using that $\zeta_{0}$ {is continuous, }$\gamma_{0}$ is continuously
differentiable and $C_{22}(\gamma_{0}(c),\zeta_{0}(c),c)\neq0$, we obtain that
$\zeta_{0}$ is continuously differentiable and%

\[%
\begin{array}
[c]{lll}%
0 & = & \partial_{y}C_{0}(\gamma_{0}(c),\zeta_{0}(c),c)\gamma_{0}^{\prime
}(c)+\partial_{z}C_{0}(\gamma_{0}(c),\zeta_{0}(c),c)\zeta_{0}^{\prime
}(c)+\partial_{c}C_{0}(\gamma_{0}(c),\zeta_{0}(c),c)\\
& = & C_{21}(\gamma_{0}(c),\zeta_{0}(c),c)\gamma_{0}^{\prime}(c)+C_{22}%
(\gamma_{0}(c),\zeta(c),c)\zeta_{0}^{\prime}(c)-C_{20}(\gamma_{0}(c),\zeta
_{0}(c),c).
\end{array}
\]
Using the first equation of (\ref{Ecuacion diferencial de z0}), we get the
second equation of (\ref{Ecuacion diferencial de z0}). By a recursive
argument, we finally obtain that $\gamma_{0}$ and $\zeta_{0}$ are infinitely
differentiable. \hfill $\blacksquare$\newline

Let us study the uniqueness of the solution of
(\ref{Ecuacion diferencial de z0}) with boundary condition
(\ref{Condicion de Optimo en cbarra}). We know that if $\left(  \gamma
,\zeta\right)  \in\mathcal{B}$ is a solution, then $\gamma(\overline
{c})=b^{\ast}(\overline{c})$, the optimal threshold defined in
\eqref{Ecuacion bestrella}. In order to obtain $\zeta(\overline{c}),$ we have
to find a zero of $C_{0}(b^{\ast}(\overline{c}),\cdot,\overline{c})$ in
$(b^{\ast}(\overline{c}),\infty).$ Let us assume that there exists a unique
zero $z^{\ast}(\overline{c})$ of $C_{0}(b^{\ast}(\overline{c}),\cdot
,\overline{c})$ in $(b^{\ast}(\overline{c}),\infty)$. In the next proposition
we show that, under this assumption, the existence of a solution $\left(
\gamma,\zeta\right)  $ of (\ref{Ecuacion diferencial de z0}) implies uniqueness.

In Section \ref{Asymptotic Values}, we will show that there is a unique zero
$z^{\ast}(\overline{c})$ of $C_{0}(b^{\ast}(\overline{c}),\cdot,\overline{c})$
in $(b^{\ast}(\overline{c}),\infty)$ for $\overline{c}$ large enough. Also, we
check this assumption in the numerical examples for each set of parameters.

\begin{proposition}
\label{Unicidad de ecuacion diferencial} Let us assume that there exists a
unique zero $z^{\ast}(\overline{c})$ of $C_{0}(b^{\ast}(\overline{c}%
),\cdot,\overline{c})$ in $(b^{\ast}(\overline{c}),\infty)$. If $\left(
\gamma_{1},\zeta_{1}\right)  \in\mathcal{B}\ $and $\left(  \gamma_{2}%
,\zeta_{2}\right)  \in\mathcal{B}$ are two solutions of the system of
differential equations (\ref{Ecuacion diferencial de z0}) satisfying the
boundary condition (\ref{Condicion de Optimo en cbarra}), then $\left(
\gamma_{1},\zeta_{1}\right)  =\left(  \gamma_{2},\zeta_{2}\right)  $.
\end{proposition}

\textit{Proof.} Consider%
\[
c_{m}=\min\left\{  c\in\lbrack0,\overline{c}]:\left(  \gamma_{1}(d),\zeta
_{1}(c)\right)  =\left(  \gamma_{2}(d),~\zeta_{2}(d)\right)  \text{ for }%
d\in\lbrack c,\overline{c}]\right\}  .
\]
Let us call%
\[
F_{1}(y,z,c)=\left(  C_{10}(y,z,c)C_{22}(y,z,c),C_{20}(y,z,c)C_{11}%
(y,z,c)-C_{21}(y,z,c)C_{10}(y,z,c)\right)  ,
\]%
\[
F_{2}(y,z,c)=C_{11}(y,z,c)C_{22}(y,z,c),
\]
and $F(y,z,c)=F_{1}(y.z,c)/F_{2}(y.z,c).$ Note that $F_{1},$ $F_{2}$ and are
infinitely differentiable,
\[
\left(  \gamma_{i}^{\prime}(c),\zeta_{i}^{\prime}(c)\right)  =F\left(
\gamma_{i}(c),\zeta_{i}(c),c\right)  \text{ for }c\in\lbrack0,\overline{c}]
\]
and $\left(  \gamma_{i}(c_{m}),\zeta_{i}(c_{m})\right)  =\left(  \gamma
_{i}(c_{m}),\zeta_{i}(c_{m})\right)  \ $for $i=1,2\ $and $F_{2}(\gamma
_{1}(c_{m}),\zeta_{1}(c_{m}),c_{m})=F_{2}(\gamma_{2}(c_{m}),\zeta_{2}%
(c_{m}),c_{m})\neq0$. If $c_{m}=0$, we have the result. On the other hand, for
$c_{m}>0$, using the Picard-Lindel\"{o}f theorem we have that there exists a
unique solution of (\ref{Ecuacion diferencial de z0}) with boundary condition
$\zeta(c_{m})=\zeta_{1}(c_{m})$ in $[\max\{c_{m}-\delta,0\},c_{m}]$ for some
$\delta>0,$ which is a contradiction. \hfill $\blacksquare$\newline

Let us now introduce a lower bound $\underline{c}$ for the dividend rate (to be specified later), and denote by 
$\left(  \overline{\gamma},\overline{\zeta}\right)  $ a solution of
(\ref{Ecuacion diferencial de z0}) in $[\underline{c},\overline{c}]$ with
boundary conditions (\ref{Condicion de Optimo en cbarra}).
\begin{remark}
Since the functions
$C_{ij}$ are infinitely differentiable, a recursive argument establishes that $\overline{\gamma}$ and
$\overline{\zeta}$ are also infinitely differentiable.
\end{remark}

In the next proposition, we state that the value function $W^{\overline{\gamma},\overline{\zeta}}$ satisfies a
smooth-pasting property on the two free-boundary curves. {Note that this extends \cite[Prop.5.13]{AAM21} from the ratcheting case with one free boundary to our present drawdown case. For a general account on conditions for smooth-pasting when the value function is not necessarily smooth, see e.g.\ Guo and Tomecek \cite{Guo}. 
}
\begin{proposition}
	\label{smooth pasting} If a pair of infinitely differentiable functions
	$(\gamma,\zeta)\in\mathcal{B}$ satisfies
	\[
	\partial_{xx}W^{\gamma,\zeta}(\gamma(c)^{+},c)=\partial_{xx}W^{\gamma,\zeta
	}(\gamma(c)^{-},c)\text{ for }c\in\lbrack\underline{c},\overline{c}]
	\]
	and
	\[
	\partial_{cx}W^{\gamma,\zeta}(\zeta(c),c)=\partial_{cc}W^{\gamma,\zeta}%
	(\zeta(c),c)=0\text{ for }c\in\lbrack\underline{c},\overline{c}],
	\]
	then $(\gamma,\zeta)\ ${is a solution of both (\ref{First equation}) and}
	(\ref{Second equation}) in $[\underline{c},\overline{c}]$ with boundary
	conditions (\ref{Condicion de Optimo en cbarra}). Moreover $\partial
	_{x}(W^{\gamma,\zeta})(\gamma(c),c)=1$ for $c\in\lbrack\underline{c}%
	,\overline{c}]$. Conversely, let $(\bar{\gamma},\bar{\zeta})\ ${be a solution
		of} (\ref{Ecuacion diferencial de z0}) in $[\underline{c},\overline{c}]$ with
	boundary conditions (\ref{Condicion de Optimo en cbarra}), then $W^{\bar
		{\gamma},\bar{\zeta}}$ satisfies the smooth-pasting properties
	\[
	\partial_{xx}W^{\bar{\gamma},\bar{\zeta}}(\overline{\gamma}(c)^{+}%
	,c)=\partial_{xx}W^{\bar{\gamma},\bar{\zeta}}(\overline{\gamma}(c)^{-}%
	,c)\text{ for }c\in\lbrack\underline{c},\overline{c}]
	\]
	and
	\[
	\partial_{cx}W^{\bar{\gamma},\bar{\zeta}}(\bar{\zeta}(c),c)=\partial
	_{cc}W^{\bar{\gamma},\bar{\zeta}}(\bar{\zeta}(c),c)=0\text{ for }c\in
	\lbrack\underline{c},\overline{c}]\text{.}%
	\]
	
\end{proposition}

\textit{Proof.} Take a pair of infinitely differentiable functions $(\gamma,\zeta
)\in\mathcal{B}$, and let us consider the function $H^{\gamma,\zeta}(x,c)$
introduced in Lemma \ref{Funcion de valor Non Change Region}. Firstly, note
that $H^{\gamma,\zeta}$ satisfies {$\mathcal{L}^{ac}(H^{\gamma,\zeta})(x,c)=0$
} for $0\leq x\leq\gamma(c),$ {$\mathcal{L}^{c}(H^{\gamma,\zeta})(x,c)=0$} for
$x\geq\gamma(c)$, $H^{\gamma,\zeta}(0,c)=0,$ $\left.  \partial_{c}%
H^{\gamma,\zeta}(x,c)\right\vert _{x=\zeta(c)}=0$ and $H^{\gamma,\zeta
}(x,\overline{c})=v\left(  x,\overline{c},\gamma(\overline{c})\right)  $. So
we have, for $x>\gamma(c)$,%
\[%
\begin{array}
	[c]{lll}%
	\partial_{c}H^{\gamma,\zeta}(x,c) & = & f_{21}(\gamma(c),x,c)\left(
	-b_{0}(\gamma(c),x,\gamma^{\prime}(c),c)-A^{\gamma,\zeta}(c)b_{1}%
	(\gamma(c),x,\gamma^{\prime}(c),c)\right)
\end{array}
\]
and%

\[
\left.  \partial_{cx}H^{\gamma,\zeta}(x,c)\right\vert _{x=\zeta(c)}%
=f_{21}(\gamma(c),\zeta(c),c)\left(  -\left.  \partial_{x}b_{0}(\gamma
(c),x,\gamma^{\prime}(c),c)\right\vert _{x=\zeta(c)}-A^{\gamma,\zeta
}(c)\left.  \partial_{x}b_{1}(\gamma(c),x,\gamma^{\prime}(c),c)\right\vert
_{x=\zeta(c)}\right)  .
\]
Since, by Lemma \ref{Lema Denominador Positivo}, $f_{21}%
(y,x,c)=d(y,x,c)/(\theta_{1}(c)$$-\theta_{2}(c))>0$ for $x>y$, we obtain that
$\left.  \partial_{cx}H^{\gamma,\zeta}(x,c)\right\vert _{x=\zeta(c)}=0$ if and
only if (\ref{First equation}) holds for $(\gamma,\zeta)$ in $[\underline
{c},\overline{c}]$. As $W^{\gamma,\zeta}(x,c)=H^{\gamma,\zeta}(x,c)$ for
$x<\zeta(c)$ and $W^{\gamma,\zeta}(x,c)=H^{\gamma,\zeta}(x,C(x,c))$ for
$x\geq\zeta(c)$, we get $\partial_{c}W^{\gamma,\zeta}(x,c)=0$ for $x\geq
\zeta(c)$ and consequently $\partial_{cx}W^{\gamma,\zeta}(\zeta(c),c)=0$.
Moreover, $\partial_{c}H^{\gamma,\zeta}(\zeta(c),c)=0$ for $c\in
\lbrack\underline{c},\overline{c}],$ and so%
\[%
\begin{array}
	[c]{lll}%
	0 & = & \frac{d}{dc}(\partial_{c}H^{\gamma,\zeta}(\zeta(c),c))\text{ }\\
	& = & \partial_{cc}H^{\gamma,\zeta}(\zeta(c),c)+\partial_{cx}H^{\gamma,\zeta
	}(\zeta(c),c)\zeta^{\prime}(c)\\
	& = & \partial_{cc}H^{\gamma,\zeta}(\zeta(c),c).
\end{array}
\]
Altogether, since $W^{\gamma,\zeta}(x,c)=H^{\gamma,\zeta}(x,C(x,c))$ if
$x\geq\zeta(c)$, we get $\partial_{cc}W^{\gamma,\zeta}(x,c)=0$ if $x\geq
\zeta(c)$ and so $\partial_{cc}W^{\gamma,\zeta}(\zeta(c),c)=0$.

Secondly, by the definitions in Section \ref{Seccion Formulas}, we have that
$\partial_{y}f_{11}(y,c)\neq0$ and $b_{11}(y,z,c)>0$%
\[
-\dfrac{b_{01}(y,c)}{qb_{11}(y,c)}=\dfrac{1-\partial_{y}f_{10}(y,c)}%
{\partial_{y}f_{11}(y,c)}.
\]

So
\[
A^{\gamma,\zeta}(c)=-\frac{\partial_{w}b_{0}(\gamma(c),\zeta(c),\gamma
	^{\prime}(c),c)}{\partial_{w}b_{1}(\gamma(c),\zeta(c),\gamma^{\prime}%
	(c),c)}=-\frac{b_{01}(\gamma(c),c)}{qb_{11}(\gamma(c),c)}=\dfrac
{1-\partial_{y}f_{10}(\gamma(c),c)}{\partial_{y}f_{11}(\gamma(c),c)}.
\]
And then \bigskip%
\[
\partial_{x}(W^{\gamma,\zeta})(\gamma(c),c)=\partial_{x}(W^{\gamma,\zeta
})(\gamma(c)^{-},c)=\partial_{y}f_{10}(\gamma(c),c))+A^{\gamma,\zeta
}(c)\partial_{y}f_{11}(\gamma(c),c)=1
\]
if, and only if, (\ref{Second equation}) holds for $(\gamma,\zeta)$ in
$[\underline{c},\overline{c}]$. Note that, since $\mathcal{L}^{c}%
W^{\gamma,\zeta}(\gamma(c)^{+},c)=\mathcal{L}^{ac}W^{\gamma,\zeta}%
(\gamma(c)^{-},c)=0$ and $\partial_{x}(W^{\gamma,\zeta})(\gamma(c),c)=1$, we
obtain
\[
0=\mathcal{L}^{c}W^{\gamma,\zeta}(\gamma(c)^{+},c)-\mathcal{L}^{ac}%
W^{\gamma,\zeta}(\gamma(c)^{-},c)=\frac{\sigma^{2}}{2}(\partial_{xx}%
W^{\gamma,\zeta}(\gamma(c)^{+},c)-\partial_{xx}W^{\gamma,\zeta}(\gamma
(c)^{-},c)).\newline%
\]
{Therefore, when} $(\bar{\gamma},\bar{\zeta})\ ${is a solution of}
(\ref{Ecuacion diferencial de z0}), it satisfies
(\ref{First equation}) and (\ref{Second equation}), and so%
\[
\partial_{xx}W^{\bar{\gamma},\bar{\zeta}}(\overline{\gamma}(c)^{+}%
,c)=\partial_{xx}W^{\bar{\gamma},\bar{\zeta}}(\overline{\gamma}(c)^{-}%
,c)\text{ for }c\in\lbrack\underline{c},\overline{c}]
\]
and
\[
\partial_{cx}W^{\bar{\gamma},\bar{\zeta}}(\bar{\zeta}(c),c)=\partial
_{cc}W^{\bar{\gamma},\bar{\zeta}}(\bar{\zeta}(c),c)=0\text{ for }c\in
\lbrack\underline{c},\overline{c}].
\]
\hfill $\blacksquare$\\

\noindent In the next proposition, we show more regularity for $W^{\overline{\gamma}%
,\overline{\zeta}}$ in the case that $\overline{\zeta}$ is strictly monotone.

\begin{proposition}
\label{z creciente} If $(\overline{\gamma},\overline{\zeta})\ ${is a solution of}
(\ref{Ecuacion diferencial de z0}) in $[\underline{c},\overline{c}]$ with
boundary conditions (\ref{Condicion de Optimo en cbarra}) and $\overline{\zeta
}^{\prime}(c)\neq0$ in $[\underline{c},\overline{c}]$, then $W^{\overline{\gamma
},\overline{\zeta}}$ is (2,1)-differentiable.
\end{proposition}

\textit{Proof.} It holds that%
\[%
\begin{array}
[c]{rl}%
\left.  f_{10}(x,c)\right\vert _{x=y}=\left.  f_{20}(y,x,c)\right\vert
_{x=y}, & \left.  f_{11}(x,c)\right\vert _{x=y}=\left.  f_{21}%
(y,x,c)\right\vert _{x=y},\\
\left.  \partial_{x}f_{10}(x,c)\right\vert _{x=y}=\left.  \partial_{x}%
f_{20}(y,x,c)\right\vert _{x=y}, & \left.  \partial_{x}f_{11}(x,c)\right\vert
_{x=y}=\left.  \partial_{x}f_{21}(y,x,c)\right\vert _{x=y},\\
\left.  \partial_{c}f_{10}(x,c)\right\vert _{x=y}=\left.  \partial_{c}%
f_{20}(y,x,c)\right\vert _{x=y}, & \left.  \partial_{c}f_{11}(x,c)\right\vert
_{x=y}=\left.  \partial_{c}f_{21}(y,x,c)\right\vert _{x=y},\\
\left.  \partial_{cx}f_{10}(x,c)\right\vert _{x=y}=\left.  \partial_{cx}%
f_{20}(y,x,c)\right\vert _{x=y}, & \left.  \partial_{cx}f_{11}(x,c)\right\vert
_{x=y}=\left.  \partial_{cx}f_{21}(y,x,c)\right\vert _{x=y},\\
\left.  \partial_{y}f_{20}(y,x,c)\right\vert _{x=y}=\left.  \partial_{y}%
f_{21}(y,x,c)\right\vert _{x=y}=0, & \left.  \partial_{cy}f_{20}%
(y,x,c)\right\vert _{x=y}=\left.  \partial_{cy}f_{21}(y,x,c)\right\vert
_{x=y}=0.
\end{array}
\]
By Proposition \ref{smooth pasting}, $W_{xx}^{\overline{\gamma},\overline{\zeta}%
}(x,c)$ is continuous at $x=$ $\overline{\gamma}(c)$ and so $W^{\overline{\gamma
},\overline{\zeta}}$ is (2,1)-differentiable for $x$ $<\overline{\zeta}(c).$

In the case that $\overline{\zeta}^{\prime}(c)>0$ in $[\underline{c},\overline{c}%
]$, the inverse $\overline{\zeta}^{-1}$ exists and $\ell(x,c)$ can be written as
\[
\ell(x,c)=\left\{
\begin{array}
[c]{lll}%
\overline{c} & \text{if} & \overline{\zeta}(\overline{c})\leq x,\\
\overline{\zeta}^{-1}(x) & \text{if} & \overline{\zeta}(c)\leq x<\overline{\zeta}%
(\overline{c}).
\end{array}
\right.
\]
In order to show that $W^{\overline{\gamma},\overline{\zeta}}$ is (2,1)-differentiable,
it is enough to prove that $\partial_{xx}W^{\overline{\zeta}}(x^{+},c)=$
$\partial_{xx}W^{\overline{\zeta}}(x^{-},c)$ for $\overline{\zeta}(c)\leq x<\overline{\zeta
}(\overline{c})$. We have, by Proposition \ref{smooth pasting},%
\[%
\begin{array}
[c]{lll}%
\partial_{x}W^{\overline{\gamma},\overline{\zeta}}(x^{+},c) & = & \partial_{x}%
H^{\overline{\gamma},\overline{\zeta}}(x,\overline{\zeta}^{-1}(x))+\partial_{c}H^{\bar
{\gamma},\overline{\zeta}}(x,\overline{\zeta}(x))\left(  \overline{\zeta}^{-1}\right)
^{\prime}(x)\\
& = & \partial_{x}H^{\overline{\gamma},\overline{\zeta}}(x,\overline{\zeta}^{-1}(x))\\
& = & \partial_{x}W^{\overline{\gamma},\overline{\zeta}}(x^{-},c).
\end{array}
\]
Consequently,
\[%
\begin{array}
[c]{lll}%
\partial_{xx}W^{\overline{\gamma},\overline{\zeta}}(x^{+},c) & = & \partial_{xx}%
H^{\overline{\gamma},\overline{\zeta}}(x,\overline{\zeta}^{-1}(x))+\partial_{cx}%
H^{\overline{\gamma},\overline{\zeta}}(x,\overline{\zeta}^{-1}(x))\left(  \overline{\zeta}%
^{-1}\right)  ^{\prime}(x)\\
& = & \partial_{xx}H^{\overline{\gamma},\overline{\zeta}}(x,\overline{\zeta}^{-1}(x))\\
& = & \partial_{xx}W^{\overline{\gamma},\overline{\zeta}}(x^{-},c)\text{. }%
\end{array}
\]

In the case that $\overline{\zeta}^{\prime}(c)<0$ in $[\underline{c},\overline{c}%
]$, $\ell(x,c)=\overline{c}$ for $x\geq\overline{\zeta}(c)$. In order to show that
$W^{\overline{\gamma},\overline{\zeta}}$ is (2,1)-differentiable, it is sufficient to prove
that $\partial_{xx}W^{\overline{\zeta}}(x^{+},c)=$ $\partial_{xx}W^{\overline{\zeta}%
}(x^{-},c)$ for $x=\overline{\zeta}(c).$ Since $\ell(x,c)=\overline{c},$%

\[
H^{\overline{\gamma},\overline{\zeta}}(\overline{\zeta}(c),c)=W^{\overline{\gamma},\overline{\zeta}%
}(\overline{\zeta}(c),c)=v^{\overline{c}}(\overline{\zeta}(c)).
\]
That is, we have,%
\[
H_{x}^{\overline{\gamma},\overline{\zeta}}(\overline{\zeta}(c),c)\overline{\zeta}^{\prime
}(c)+H_{c}^{\overline{\gamma},\overline{\zeta}}(\overline{\zeta}(c),c)=(v^{\overline{c}%
})^{\prime}(\overline{\zeta}(c))\overline{\zeta}^{\prime}(c).
\]
Since $H_{c}^{\overline{\gamma},\overline{\zeta}}(\overline{\zeta}(c),c)=W_{c}%
^{\overline{\gamma},\overline{\zeta}}(\overline{\zeta}(c),c)=0$ and $\overline{\zeta}^{\prime
}(c)<0$ we get
\[
W_{x}^{\overline{\gamma},\overline{\zeta}}((\overline{\zeta}(c))^{-},c)=H_{x}^{\overline{\gamma
},\overline{\zeta}}(\overline{\zeta}(c),c)=(v^{\overline{c}})^{\prime}(\overline{\zeta
}(c))=W_{x}^{\overline{\gamma},\overline{\zeta}}((\overline{\zeta}(c))^{+},c),
\]
so that $H_{x}^{\overline{\gamma},\overline{\zeta}}(\overline{\zeta}(c),c)=(v^{\overline{c}%
})^{\prime}(\overline{\zeta}(c))$. Hence, taking the derivative one more time with respect to
$c$ we get
\[
H_{xx}^{\overline{\gamma},\overline{\zeta}}(\overline{\zeta}(c),c)\overline{\zeta}^{\prime
}(c)+H_{xc}^{\overline{\gamma},\overline{\zeta}}(\overline{\zeta}(c),c)=(v^{\overline{c}%
})^{\prime\prime}(\overline{\zeta}(c))\overline{\zeta}^{\prime}(c).
\]
By virtue of Proposition \ref{smooth pasting}, $H_{xc}^{\overline{\gamma}%
,\overline{\zeta}}(\overline{\zeta}(c),c)=H_{cx}^{\overline{\gamma},\overline{\zeta}}(\overline{\zeta
}(c),c)=0$ and $\overline{\zeta}^{\prime}(c)<0$, and we obtain
\[
W_{xx}^{\overline{\gamma},\overline{\zeta}}((\overline{\zeta}(c))^{-},c)=H_{xx}^{\overline{\gamma
},\overline{\zeta}}(\overline{\zeta}(c),c)=(v^{\overline{c}})^{\prime\prime}(\overline{\zeta
}(c))=W_{xx}^{\overline{\gamma},\overline{\zeta}}((\overline{\zeta}(c))^{+},c)\text{.
}
\]
\hfill $\blacksquare$

\begin{proposition}
\label{Verificacion Curva} Let $(\overline{\gamma},\overline{\zeta})\ ${be a solution
of} (\ref{Ecuacion diferencial de z0}) in $[\underline{c},\overline{c}]$ with
boundary conditions (\ref{Condicion de Optimo en cbarra}) such that the
function $W^{\overline{\gamma},\overline{\zeta}}$ is (2,1)-differentiable and satisfies
\[
\partial_{c}W^{\overline{\gamma},\overline{\zeta}}(x,c)\leq0\text{ for }x\in
\lbrack0,\overline{\zeta}(c))\text{ }%
\]
and%
\[
\partial_{x}W^{\overline{\gamma},\overline{\zeta}}(x,c)\geq1\text{ for }x\in
\lbrack0,\overline{\gamma}(c))\text{ and }\partial_{x}W^{\overline{\gamma}%
,\overline{\zeta}}(x,c)\leq1\text{ for }x\in\lbrack\overline{\gamma}(c),\overline
{\zeta}(c)]
\]
for $c\in\lbrack\underline{c},\overline{c})$, then $W^{\overline{\gamma},\bar
{\zeta}}=V$.
\end{proposition}

{ }\textit{Proof.} Since $\overline{\zeta}\ $is continuous in $[\underline{c},\overline{c}],$
there exists $M=\max_{c\in\lbrack\underline{c},\overline{c}]}\overline{\zeta
}(c).$ By definition, if $x\geq M$ then $\ell(x,c)=\overline{c}~$and
$W^{\overline{\gamma},\overline{\zeta}}(x,c)=v^{\overline{c}}(x),$ so $\lim
_{x\rightarrow\infty}W^{\overline{\gamma},\overline{\zeta}}(x,c)=\lim_{x\rightarrow
\infty}v^{\overline{c}}(x)=\overline{c}/q.$

By (\ref{Derivada en tapa}), we have that
\begin{equation}
\partial_{x}v^{\overline{c}}(x)\leq1~\text{for }x\geq\overline{\zeta
}(\overline{c})\text{ and }\partial_{x}v^{\overline{c}}(b^{\ast}(\overline
{c}))=1. \label{Vx -1 negativo en la tapa}%
\end{equation}
Since
\[
\mathcal{L}^{ac}(W^{\overline{\gamma},\overline{\zeta}})(x,c)-\mathcal{L}^{c}%
(W^{\overline{\gamma},\overline{\zeta}})(x,c)=(c-ac)\partial_{x}W^{\overline{\gamma}%
,\overline{\zeta}}(x,c)+(ac-c)=c(1-a)\left(  \partial_{x}W^{\overline{\gamma},\bar
{\zeta}}(x,c)-1\right)  ,
\]
we get $\mathcal{L}^{ac}(W^{\overline{\gamma},\overline{\zeta}})(x,c)\leq\mathcal{L}%
^{c}(W^{\overline{\gamma},\overline{\zeta}})(x,c)=0$ for $x\in\lbrack\overline{\gamma
}(c),\overline{\zeta}(c)]$ and $\mathcal{L}^{c}(W^{\overline{\gamma},\overline{\zeta}%
})(x,c)\leq\mathcal{L}^{ac}(W^{\overline{\gamma},\overline{\zeta}})(x,c)=0$ for
$x\in\lbrack0,\overline{\gamma}(c)]$.

By Theorem \ref{Caracterizacion Continua}, it remains to prove that
$\mathcal{L}^{ac}W^{\overline{\gamma},\overline{\zeta}}(x,c)\leq0$ and $\mathcal{L}%
^{c}W^{\overline{\gamma},\overline{\zeta}}(x,c)\leq0$ for $x\geq\overline{\zeta}(c),$
$c\in\lbrack\underline{c},\overline{c})$. We have that
\[
\ell(x,c)=\max\{h\in\lbrack c,\overline{c}]:\overline{\zeta}(d)\leq
x{\large \ }\text{{\large for}}{\large \ }d\in\lbrack c,h)\}
\]
satisfies $\ell(x,c)\geq c$, and also either $\ell(x,c)=\overline{c}$ or
$\overline{\zeta}(\ell(x,c))=x.$ So, we obtain $\left.  \mathcal{L}^{\alpha
}W^{\overline{\gamma},\overline{\zeta}}(x,\alpha)\right\vert _{\alpha=\ell(x,c)}=0$ and then%

\[%
\begin{array}
[c]{lll}%
\mathcal{L}^{c}W^{\overline{\gamma},\overline{\zeta}}(x,c) & = & \mathcal{L}^{\ell
(x,c)}W^{\overline{\gamma},\overline{\zeta}}(x,c)+(\ell(x,c)-c)(\partial_{x}%
W^{\overline{\gamma},\overline{\zeta}}(x,\ell(x,c))-1)\\
& = & (\ell(x,c)-c)(\partial_{x}W^{\overline{\gamma},\overline{\zeta}}(x,\ell
(x,c))-1)\leq0\text{,}%
\end{array}
\]
because we have, from (\ref{Vx -1 negativo en la tapa}) and $\partial
_{x}W^{\overline{\gamma},\overline{\zeta}}(\overline{\zeta}(c),c)\leq1$ for $c\in
\lbrack\underline{c},\overline{c}]$, that $\partial_{x}W^{\overline{\gamma}%
,\overline{\zeta}}(x,\ell(x,c))\leq1$. Also,
\[
\mathcal{L}^{ac}(W^{\overline{\gamma},\overline{\zeta}})(x,c)-\mathcal{L}^{c}%
(W^{\overline{\gamma},\overline{\zeta}})(x,c)=c(1-a)\left(  \partial_{x}W^{\overline{\gamma
},\overline{\zeta}}(x,c)-1\right)  =c(1-a)\left(  \partial_{x}W^{\overline{\gamma}%
,\overline{\zeta}}(x,\ell(x,c))-1\right)  \leq0
\]
for $x\geq\overline{\zeta}(c),$ $c\in\lbrack\underline{c},\overline{c})$.
\hfill $\blacksquare$

\begin{remark}
\normalfont\label{Conjetura dos curvas} We conjecture that there is always a
unique zero $z^{\ast}(\overline{c})$ of $C_{0}(b^{\ast}(\overline{c}%
),\cdot,\overline{c})$ in $(b^{\ast}(\overline{c}),\infty)$ for $\overline
{c}>$ $q\sigma^{2}/(2\mu)$, that there exists a solution $\left(
\overline{\gamma},\overline{\zeta}\right)  \in\mathcal{B}$ of the system of
differential equations (\ref{Ecuacion diferencial de z0}) satisfying the
boundary condition (\ref{Condicion de Optimo en cbarra}), and that the value
function $W^{\overline{\gamma},\overline{\zeta}}$ is a viscosity supersolution
of the HJB equation (\ref{HJB equation}). In such a case, $(\overline{\gamma
},\overline{\zeta})=(\gamma_{0},\zeta_{0})$ and $W^{\overline{\gamma
},\overline{\zeta}}$ is the optimal value function $V$. Moreover, the optimal
strategy is then a two-curve strategy. In Section \ref{Optimal strategies for c large enough} we will show that that this
conjecture in any case holds in $\left[  \underline{c},\overline{c}\right]  $ for
$\overline{c}$ large enough and some suitable  $\underline{c}<\overline{c}$, and in Section \ref{Numerical examples} we also show numerically that this conjecture is true for further instances. 
\end{remark}

\section{Asymptotic values as $\overline{c}\rightarrow\infty$
\label{Asymptotic Values}}

The symbolic computations of this section are highly involved, so we use the
Wolfram Mathematica software to obtain Taylor expansions. {Note that all results of this section are derived for $0<a<1$, and the resulting expressions may not necessarily be applicable for the limit to $a=1$, as dominant terms in the asymptotics may change.}

Recall the boundary condition $C_{0}(b^{\ast}(\overline{c}),\cdot,\overline
{c})=0$ of the differential equation of the last section,
cf.\ \eqref{Condicion de Optimo en cbarra}. Note that for $\overline{c}%
>\frac{q\sigma^{2}}{2\mu},$ we have from Remark \ref{Rermark caso interesante}
that $b^{\ast}(\overline{c})$ is the unique positive $b$ satisfying
(\ref{Ecuacion implicita de bestrella}). In this section, we show that there
is a unique zero $z^{\ast}(\overline{c})$ of $C_{0}(b^{\ast}(\overline
{c}),\cdot,\overline{c})$ in $(b^{\ast}(\overline{c}),\infty)$ for
$\overline{c}$ large enough and that
\[
\lim_{\overline{c}\rightarrow\infty}\left(  b^{\ast}(\overline{c}),z^{\ast
}(\overline{c})\right)  =\left(  \frac{\mu}{q},\frac{\mu}{q}\big(1+\frac
{1}{\sqrt{a}}\big)\right)  .
\]
We also show that $\lim_{\overline{c}\rightarrow\infty}V_{a}^{\overline{c}%
}(x,\overline{c})\searrow x$ for $0<x<\lim_{\overline{c}\rightarrow\infty
}z^{\ast}(\overline{c})=\frac{\mu}{q}(1+\frac{1}{\sqrt{a}})$ and
$\lim_{\overline{c}\rightarrow\infty}V_{a}^{\overline{c}}(x,\overline
{c})\nearrow x$ for $x>\lim_{\overline{c}\rightarrow\infty}z^{\ast}%
(\overline{c})=\frac{\mu}{q}(1+\frac{1}{\sqrt{a}})$.

In the rest of the section we denote $V_{a}^{\overline{c}}$ by $V^{\overline
{c}}$.

\begin{proposition}
\label{Limite de condiciones en cbarra} It holds that $\lim_{\overline
{c}\rightarrow\infty}b^{\ast}(\overline{c})=\mu/q$. Moreover precisely, the
Taylor expansion of $b^{\ast}(\overline{c})$ at $\overline{c}=\infty$ is given
by
\begin{equation}
b^{\ast}(\overline{c})=\frac{\mu}{q}-\frac{\mu^{2}+aq\sigma^{2}}{2aq}\frac
{1}{\overline{c}}+O\left(  \frac{1}{\overline{c}^{2}}\right)  \text{.}
\label{expansion de bestrella en c}%
\end{equation}

\end{proposition}

\textit{Proof.} We have from (\ref{Ecuacion bestrella}) that%
\begin{equation}
\partial_{b}B(\overline{c},b^{\ast}(\overline{c}))=0. \label{Condicion b*}%
\end{equation}
But%
\[
\partial_{b}B(\overline{c},b)=\frac{c\sqrt{(\mu-a\overline{c})^{2}%
+2q\sigma^{2}}}{q\left(  e^{\theta_{1}(a\overline{c})b}(\theta_{1}%
(a\overline{c})-\theta_{2}(\overline{c}))+e^{\theta_{2}(a\overline{c}%
)b}(\theta_{2}(\overline{c})-\theta_{2}(a\overline{c}))\right)  ^{2}}\cdot
E(\overline{c},b),
\]
where
\begin{align}
E(\overline{c},b)  &  =e^{\theta_{1}(a\overline{c})b}(a-1)\theta_{2}%
(\overline{c})(\theta_{2}(\overline{c})-\theta_{1}(a\overline{c}))\theta
_{1}(a\overline{c})\label{E(c,b)}\\
&  +e^{\theta_{2}(a\overline{c})b}(1-a)\theta_{2}(\overline{c})(\theta
_{2}(\overline{c})-\theta_{2}(a\overline{c}))\theta_{2}(a\overline
{c})\nonumber\\
&  +e^{\left(  \theta_{1}(a\overline{c})+\theta_{2}(a\overline{c})\right)
b}a(\theta_{2}(\overline{c})-\theta_{2}(a\overline{c}))(\theta_{2}%
(\overline{c})-\theta_{1}(a\overline{c}))(\theta_{2}(a\overline{c})-\theta
_{1}(a\overline{c})).\nonumber
\end{align}
Let us define $F_{0}(\overline{c},b):=E(\overline{c},b)/e^{\theta
_{1}(a\overline{c})b}$. The Taylor expansions of $\theta_{1}(c)$ and
$\theta_{2}(c)$ at $c=\infty\ $are given by
\begin{equation}
\theta_{1}(c)=\frac{2}{\sigma^{2}}c-\frac{2\mu}{\sigma^{2}}+q\frac{1}%
{c}+O(\frac{1}{c^{2}})\ \text{and }\theta_{2}(c)=-q\frac{1}{c}-q\mu\frac
{1}{c^{2}}+O(\frac{1}{c^{3}}). \label{Taylor titas}%
\end{equation}
Let us prove first that there is not a sequence $b^{\ast}(c_{n})\rightarrow
\infty$ with $c_{n}\rightarrow\infty$. Using (\ref{Taylor titas}), we obtain%

\[
\lim_{n\rightarrow\infty}F_{0}(c_{n},b^{\ast}(c_{n}))=\lim_{n\rightarrow
\infty}\frac{4(a-1)a^{2}q}{\sigma^{4}}c_{n}(1-e^{-\frac{qb^{\ast}(c_{n}%
)}{ac_{n}}}).
\]
Firstly, let us assume that $b^{\ast}(c_{n})=c_{n}\ \alpha_{n}\ $with
$\alpha_{n}\rightarrow\overline{\alpha}\in(0,\infty)$. Then, since
$(1-e^{-\frac{q\overline{\alpha}}{a}})>0$ and $a<1,$
\[
0=\lim_{n\rightarrow\infty}F_{0}(c_{n},b^{\ast}(c_{n}))=-\infty,
\]
which is a contradiction. Secondly, let us assume that $b^{\ast}(c_{n}%
)=c_{n}\ \alpha_{n}\ $with $\alpha_{n}\rightarrow\infty$. Then, since
$e^{-\frac{qb^{\ast}(c_{n})}{ac_{n}}}\rightarrow0$, we have $0=\lim
_{n\rightarrow\infty}F_{0}(c_{n},b^{\ast}(c_{n}))=-\infty$ which is also a
contradiction. Finally, let us assume that $b^{\ast}(c_{n})=c_{n}\ \alpha
_{n}\ \ $with $\alpha_{n}\rightarrow0^{+}$.
\[
0=\lim_{n\rightarrow\infty}F_{0}(c_{n},b^{\ast}(c_{n}))=\lim_{n\rightarrow
\infty}\frac{4(a-1)a^{2}q}{\sigma^{4}}(\frac{1-e^{-\frac{q\alpha_{n}}{a}}%
}{\alpha_{n}})b^{\ast}(c_{n})=-\infty.
\]
Hence, $\lim\sup_{\overline{c}\rightarrow\infty}b^{\ast}(\overline{c})<\infty
$.\newline

Let us define the function $H_{0}(u,b):[0,\infty)\times(0,\infty
)\rightarrow\mathbb{R}$ as
\[
H_{0}(u,b):=\left\{
\begin{array}
[c]{lll}%
\frac{4(a-1)aq(qb-\mu)}{\sigma^{4}} & \text{for} & u=0\\
F_{0}(\frac{1}{u},b) & \text{for} & u>0.
\end{array}
\right.
\]
$H_{0}(u,b)$ is infinitely continuously differentiable because it is
infinitely continuously differentiable for $u>0$ and $\lim_{u\rightarrow0^{+}%
}F_{0}(\frac{1}{u},b)=4(a-1)aq(qb-\mu)/\sigma^{4}<\infty$. Moreover, its
first-order Taylor expansion is given by%
\[
H_{0}(u,b)=\frac{4(a-1)aq(qb-\mu)}{\sigma^{4}}+\left(  \frac{2(a-1)q-q^{2}%
b^{2}+2(1-a)\mu^{2}+aq(2b\mu+\sigma^{2}))}{\sigma^{4}}\right)  u+O(u^{2}%
)~\text{,}%
\]
From (\ref{Condicion b*}), we obtain $H_{0}(u,b^{\ast}(1/u))=0~$for $u>0$. Let
us show that $\lim_{u\rightarrow0^{+}}b^{\ast}(1/u)=\mu/q$. We have already
seen that $b^{\ast}(1/u)\ $is bounded for $u$ $\in\lbrack0,\varepsilon)$ for
some $\varepsilon>0.$ Take any convergent sequence $u_{n}\rightarrow0^{+}$
with $\lim_{n\rightarrow\infty}b^{\ast}(1/u_{n})=b_{0}<\infty$, then
\[
\lim_{n\rightarrow\infty}H_{0}(u_{n},b^{\ast}(1/u_{n}))=H_{0}(0,b_{0}%
)=\frac{4(a-1)aq(qb_{0}-\mu)}{\sigma^{4}}=0
\]
and so $b_{0}=\mu/q.$ Using that $\partial_{b}H_{0}(0,b)=\frac{4(a-1)aq^{2}%
}{\sigma^{4}}\neq0$, we conclude by the implicit function theorem, that the
function $h(u):[0,\infty)\rightarrow\mathbb{R}$ defined as $h(0)=\frac{\mu}%
{q}$ and $h(u)=b^{\ast}(\frac{1}{u})$ for $u>0$ is infinitely continuously
differentiable and the result follows. \hfill$\blacksquare$

\begin{proposition}
\label{Limites de condiciones en z barra} There exists a unique zero $z^{\ast
}(\overline{c})$ of $C_{0}(b^{\ast}(\overline{c}),\cdot,\overline{c})$ in
$(b^{\ast}(\overline{c}),\infty)$ for $\overline{c}$ large enough with
$\lim_{\overline{c}\rightarrow\infty}z^{\ast}(\overline{c})=\frac{\mu}%
{q}(1+\frac{1}{\sqrt{a}})$. More precisely, $z^{\ast}(\overline{c})$ is
infinitely continuously differentiable for $\overline{c}$ large enough and its
first-order Taylor expansion at $\overline{c}=\infty$ is given by
\begin{equation}
\label{zstar}z^{\ast}(\overline{c})=\frac{\mu}{q}\Big(1+\frac{1}{\sqrt{a}%
}\Big)+\frac{\left(  1-2\sqrt{a}-3a\right)  \mu^{2}-3(1+\sqrt{a^{3}}%
/2)q\sigma^{2}}{3q\sqrt{a^{3}}}\frac{1}{\overline{c}}+O(\frac{1}{\overline
{c}^{2}}).
\end{equation}

\end{proposition}

\textit{Proof.} Considering the function $C_{0}(y,z,c)$ defined in Section
\ref{Seccion Formulas} and the function $E(c,y)$ defined in (\ref{E(c,b)}), we
define%
\[%
\begin{array}
[c]{lll}%
\tilde{C}_{0}(y,z,c) & = & \left(  ce^{2(z-y)\theta_{1}(c)+y\theta_{1}%
(ac)}(\theta_{2}(c)-\theta_{1}(c))(\theta_{1}(ac)-\theta_{2}(c))\theta
_{1}^{\prime}(c)\right)  E(c,y)\\
&  & -\left(  ce^{2(z-y)\theta_{1}(c)+y\theta_{2}(ac)}(\theta_{2}%
(ac)-\theta_{2}(c))(\theta_{2}(c)-\theta_{1}(c))\theta_{1}^{\prime
}(c))\right)  E(c,y)\\
&  & +\frac{C_{0}(y,z,c)}{\left(  \text{$\theta_{1}(c)$}-\text{$\theta_{2}%
(c)$}\right)  }.
\end{array}
\]
Since $E(\overline{c},b^{\ast}(\overline{c}))=0$, $\theta_{2}(c)$$-\theta
_{1}(c)<0$ and $d(y,z,c)>0$, we have that $C_{0}(b^{\ast}(\overline
{c}),z^{\ast}(\overline{c}),\overline{c})=0$ is equivalent to $\tilde{C}%
_{0}(b^{\ast}(\overline{c}),z^{\ast}(\overline{c}),\overline{c})=0.$ We can
write
\begin{equation}
\tilde{C}_{0}(y,z,c)=%
{\displaystyle\sum\limits_{i=1}^{16}}
m_{i}(y,z,c)e^{g_{i}(y,z,c)}, \label{Formula C0sombrero}%
\end{equation}
where $m_{i}(y,z,c)$ are of the form%
\[
m_{i}(y,z,c)=m_{i0}(y,c)+m_{i1}(y,c)z,
\]
and $m_{i0}(y,c)$,$~m_{i1}(y,c)$ are polynomials on $\theta_{1}(c)$,
$\theta_{2}(c)$, $\theta_{1}(ac)$, $\theta_{2}(ac)$, $\theta_{1}^{\prime}(c)$,
$\theta_{2}^{\prime}(c)$, $\theta_{1}^{\prime}(ac)$, $\theta_{2}^{\prime}%
(ac)$, $y$, $c$, $a$. The functions $g_{i}(y,z,c)$ in
\eqref{Formula C0sombrero} are positive linear combinations of $(z-y)\theta
_{1}(c)$, $(z-y)\theta_{2}(c)$, $y\theta_{1}(ac)$ and $y\theta_{2}(ac)$, with
the concrete form given in Section \ref{Seccion Formulas}. Define
\[
F_{1}(y,z,c):=\frac{\tilde{C}_{0}(y,z,c)}{e^{g_{12}(y,z,c)}}.
\]
Let us show first that there is not a sequence $(z_{n},c_{n})$ with
$z_{n}>b^{\ast}(c_{n})$ such that $C_{0}(b^{\ast}(c_{n}),z_{n},c_{n})=0$,
$c_{n}\rightarrow\infty$ and $z_{n}\rightarrow\infty$. From the definitions of
the exponents $g_{i}$ given in Section \ref{Seccion Formulas} and the
expressions (\ref{Taylor titas}), we have that%
\[
\lim_{n\rightarrow\infty}F_{1}(b^{\ast}(c_{n}),z_{n},c_{n})=\lim
_{n\rightarrow\infty}%
{\displaystyle\sum\limits_{i=12}^{14}}
m_{i}(b^{\ast}(c_{n}),z_{n},c_{n})e^{g_{i}(b^{\ast}(c_{n}),z_{n},c_{n}%
)-g_{12}(b^{\ast}(c_{n}),z_{n},c_{n})},
\]
because the other terms are negligible. We can write%
\[%
\begin{array}
[c]{lll}%
m_{12,0}(y,c)=\frac{64(a-1)a^{2}}{\sigma^{12}}c^{6}+O(c^{5}), &  &
m_{12,1}(y,c)=O(c^{4}),\\
m_{13,0}(y,c)=\frac{64(a-1)^{2}a^{2}}{\sigma^{12}}c^{6}+O(c^{5}), &  &
m_{13,1}(y,c)=\frac{128(a-1)^{2}a^{2}q}{\sigma^{12}}c^{5}+O(c^{4})\\
m_{14,0}(y,c)=-\frac{64(a-1)a^{3}}{\sigma^{12}}c^{6}+O(c^{5}), &  &
m_{14,1}(y,c)=-\frac{64a^{2}(1-3a+2a^{2})q}{\sigma^{12}}c^{5}+O(c^{4})
\end{array}
\]
and%
\[%
\begin{array}
[c]{lll}%
g_{13}(y,z,c)-g_{12}(y,z,c) & = & -q(z-y)\frac{1}{c}+(1+z)O(\frac{1}{c^{2}%
}),\\
g_{14}(y,z,c)-g_{12}(y,z,c) & = & -\left(  \frac{qy}{a}+q(z-y)\right)
\frac{1}{c}+(1+z)O(\frac{1}{c^{2}}).
\end{array}
\]
If $z_{n}\rightarrow\infty$, $c_{n}\rightarrow\infty,$ with $b^{\ast}%
(c_{n})\rightarrow\frac{\mu}{q}$ we deduce that
\[
0=\lim_{n\rightarrow\infty}F_{1}(b^{\ast}(c_{n}),z_{n},c_{n})=\lim
_{n\rightarrow\infty}\frac{64(a-1)a^{2}c_{n}^{6}}{\sigma^{12}}e^{-q\frac
{z_{n}}{c_{n}}}(e^{q\frac{z_{n}}{c_{n}}}-1-q\frac{z_{n}}{c_{n}}).
\]
Firstly, let us assume that $z_{n}=c_{n}\ \alpha_{n}\ $with $\alpha
_{n}\rightarrow\overline{\alpha}\in(0,\infty)$. Then, since $e^{-q\overline
{\alpha}}(e^{q\overline{\alpha}}-1-q\overline{\alpha})>0$ and $a<1,$
\[
0=\lim_{n\rightarrow\infty}F_{1}(b^{\ast}(c_{n}),z_{n},c_{n})=-\infty
\]
which is a contradiction. \ Secondly, let us assume that $z_{n}=c_{n}%
\ \alpha_{n}\ $with $\alpha_{n}\rightarrow\infty$. Then, since
\[
e^{-q\alpha_{n}}(e^{q\alpha_{n}}-1-q\alpha_{n})=1-(1+q\alpha_{n}%
)e^{-q\alpha_{n}}\rightarrow1,
\]
we have%
\[
0=\lim_{n\rightarrow\infty}F_{1}(b^{\ast}(c_{n}),z_{n},c_{n})=-\infty
\]
which is also a contradiction. Finally, let us assume that $z_{n}=c_{n}%
\alpha_{n}\ \ $with $\alpha_{n}\rightarrow0^{+}$.
\[
0=\lim_{n\rightarrow\infty}F_{1}(b^{\ast}(c_{n}),z_{n},c_{n})=\lim
_{n\rightarrow\infty}\tfrac{64(a-1)a^{2}}{\sigma^{12}}q^{2}e^{-q\alpha_{n}%
}(\tfrac{e^{q\alpha_{n}}-1-q\alpha_{n}}{q^{2}\alpha_{n}^{2}})z_{n}^{2}%
c_{n}^{4}=-\infty
\]
which is also a contradiction. Hence, there is not such a sequence
$(z_{n},c_{n})$.

Using the Taylor expansions of $\theta_{1}(\overline{c})$, $\theta
_{2}(\overline{c})\ $and $b^{\ast}(\overline{c})$ at $\overline{c}=\infty$
given in (\ref{Taylor titas}) and Proposition
\ref{Limite de condiciones en cbarra}, we find that the function
\[
H_{1}(z,u)=\left\{
\begin{array}
[c]{ccc}%
\frac{32(a-1)a\left(  a(\mu-qz)^{2}-\mu^{2}\right)  }{\sigma^{12}} &
\text{for} & u=0\\
u^{4}F_{1}(b^{\ast}(\frac{1}{u}),z,\frac{1}{u}) & \text{for} & u>0
\end{array}
\right.
\]
is infinitely continuously differentiable, because it is infinitely
continuously differentiable for $u>0$ and
\[
\lim_{u\rightarrow0^{+}}u^{4}F_{1}(b^{\ast}(\frac{1}{u}),z,\frac{1}{u}%
)=\frac{32(a-1)a\left(  a(\mu-qz)^{2}-\mu^{2}\right)  }{\sigma^{12}}<\infty.
\]
Moreover, its first-order Taylor expansion is given by%
\[%
\begin{array}
[c]{lll}%
H_{1}(z,u) & = & \frac{32(a-1)a\left(  a(\mu-qz)^{2}-\mu^{2}\right)  }%
{\sigma^{12}}\\
&  & -\frac{32(a-1)\left(  -4\mu^{3}+3a\mu(2q^{2}z^{2}-6qz\mu+\mu^{2}%
)+a^{2}(qz-\mu)\left(  2q^{2}z^{2}-\mu^{2}-q(z\mu+3\sigma^{2}\right)  \right)
}{3\sigma^{12}}u+O(u^{2}).
\end{array}
\]
Since, the only zero of $H_{1}(z,0)$ in $[\frac{\mu}{q},\infty)$ is $\frac
{\mu}{q}(1+\frac{1}{\sqrt{a}})$ and
\[
\partial_{z}H_{1}(z,0)=\partial_{z}(\frac{32(a-1)a\left(  a(\mu-qz)^{2}%
-\mu^{2}\right)  }{\sigma^{12}})=\frac{64(1-a)a^{2}q(\mu-qz)}{\sigma^{12}}%
\neq0
\]
for $z\geq\frac{\mu}{q}$, we conclude by the implicit function theorem, that
there exists $\varepsilon>0$ and a unique infinitely continuously
differentiable function $g(u):[0,\varepsilon)\rightarrow\mathbf{R}$ with
$g(0)=\frac{\mu}{q}(1+\frac{1}{\sqrt{a}})$ and $H_{1}(g(u),u)=0$ for
$u\in\lbrack0,\varepsilon)$; also $g(u)$ is the unique zero of $H_{1}%
(\cdot,u)$ in a neighborhood $U$ of $(\frac{\mu}{q}(1+\frac{1}{\sqrt{a}}),0).$
Moreover, the first-order Taylor expansion of $g$ at $u=0$ is given by%
\begin{equation}
g(u)=\frac{\mu}{q}\Big(1+\frac{1}{\sqrt{a}}\Big)+\frac{\left(  2-4\sqrt{a}%
\mu^{2}-6a\right)  \mu^{2}-3(2+\sqrt{a^{3}})q\sigma^{2}}{6q\sqrt{a^{3}}%
}u+O(u^{2}). \label{expansion de z en u}%
\end{equation}
Let us show now that $g(u)$ is the only zero of $H_{1}(\cdot,u)$ in $(b^{\ast
}(1/u),\infty)$ for $u$ small enough. If this were not the case, there should
be a sequence $(z_{n},u_{n})_{n\geq1}$ with $z_{n}>b^{\ast}(1/u_{n}),$
$z_{n}\neq g(u_{n})$ such that $u_{n}\searrow0$ and $H_{1}(z_{n},u_{n})=0$. If
there exists a convergent subsequence $z_{n_{k}}$ with $z_{n_{k}}\rightarrow
z_{0}\in\lbrack\frac{\mu}{q},\infty),$ by continuity $H_{1}(z_{0},0)=0$ and so
$z_{0}=g(0)=\frac{\mu}{q}(1+\frac{1}{\sqrt{a}})$ and this is a contradiction
because $(z_{n_{k}},u_{n_{k}})\notin U$ for $k$ large enough. So
$z_{n}\rightarrow\infty$ and this is also a contradiction. So, from
(\ref{expansion de z en u}), we get the result. \hfill$\blacksquare$

\begin{proposition}
\label{Derivada de V cbarra} There exists a unique zero $x^{\ast}
(\overline{c})$ of $\partial_{\overline{c}}V^{\overline{c}}(x,\overline{c})$
in $(0,\infty)$ for $\overline{c}$ large enough with $x^{\ast}(\overline
{c})>b^{\ast}(\overline{c})$ and $\lim_{\overline{c}\rightarrow\infty}x^{\ast
}(\overline{c})=\frac{\mu}{q}(1+\frac{1}{\sqrt{a}})$. More precisely,
$x^{\ast}(\overline{c})$ is infinitely continuously differentiable for
$\overline{c}$ large enough and its first-order Taylor expansion at
$\overline{c}=\infty$ is given by
\begin{equation}
x^{\ast}(\overline{c})=\frac{\mu}{q}\left(  1+\frac{1}{\sqrt{a}}\right)
+\frac{(1-2\sqrt{a}-3a)\mu^{2}-3(1+\sqrt{a^{3}})q\sigma^{2}}{3\sqrt{a^{3}}
q}\frac{1}{\overline{c}}+O(\frac{1}{\overline{c}^{2}}). \label{xstarConSigma}%
\end{equation}
Moreover, $\lim_{\overline{c}\rightarrow\infty}V^{\overline{c}} (x,\overline
{c})\searrow x$ for $0<x<\frac{\mu}{q}(1+\frac{1}{\sqrt{a}})$ and
$\lim_{\overline{c}\rightarrow\infty}V^{\overline{c}}(x,\overline{c})\nearrow
x$ for $x>\frac{\mu}{q}(1+\frac{1}{\sqrt{a}}).$
\end{proposition}

\noindent\textit{Proof.} From (\ref{Formula value function refraction en b}),
we have that%
\[
V^{\overline{c}}(x,\overline{c})=v(x,\overline{c},b^{\ast}(\overline
{c}))=\Big(B(\overline{c},b^{\ast}(\overline{c}))W_{0}(x,\overline{c}%
)+\frac{a\overline{c}}{q}(1-e^{\theta_{2}(a\overline{c})x})\Big)I_{x<b^{\ast
}(\overline{c})}+\Big(\frac{\overline{c}}{q}+D(\overline{c},b^{\ast}%
(\overline{c}))e^{\theta_{2}(\overline{c})x}\Big)I_{x\geq b^{\ast}
(\overline{c})}%
\]
and from Proposition \ref{Limite de condiciones en cbarra}, we know that
$\lim_{\overline{c}\rightarrow\infty}b^{\ast}(\overline{c})\nearrow\mu
/q$.\newline

Take $x<\mu/q$, then $x<b^{\ast}(\overline{c})$ for $\overline{c}$ large
enough, so we have that%
\[
V^{\overline{c}}(x,\overline{c})=B(\overline{c},b^{\ast}(\overline{c}%
))W_{0}(x,\overline{c})+\frac{a\overline{c}}{q}(1-e^{\theta_{2}(a\overline
{c})x})
\]
and so
\[
\partial_{\overline{c}}V^{\overline{c}}(x,\overline{c})=\dfrac{F_{2}
(x,\overline{c})}{l_{0}(b^{\ast}(\overline{c}),\overline{c})},
\]
where%

\begin{equation}
l_{0}(b,c)=q\left(  (\mu-ac)^{2}+2q\sigma^{2}\right)  \left(  e^{b\theta
_{1}(ac)}(\theta_{1}(ac)-\theta_{2}(c))+e^{b\theta_{2}(ac)}(\theta
_{2}(c)-\theta_{2}(ac))\right)  ^{2}>0, \label{l0(b,c)}%
\end{equation}
and
\[
F_{2}(x,\overline{c}):=%
{\displaystyle\sum\limits_{i=1}^{11}}
l_{i}(x,b^{\ast}(\overline{c}),\overline{c})e^{h_{i}(x,b^{\ast}(\overline
{c}),\overline{c})}.
\]
Here $l_{i}(x,b,c)$ are polynomials on $\theta_{1}(c)$, $\theta_{2}(c)$,
$\theta_{1}(ac)$, $\theta_{2}(ac)$, $\theta_{1}^{\prime}(c)$, $\theta
_{2}^{\prime}(c)$, $\theta_{1}^{\prime}(ac)$, $\theta_{2}^{\prime}(ac)$, $x$,
$b$, $c$, $a$ and $h_{i}(x,b,c)$, $i=1,...,11$, are positive linear
combinations of $b\theta_{1}(ac)$, $x\theta_{1}(ac)$, $b\theta_{2}(ac)$ and
$x\theta_{2}(ac)$ stated in detail in Section \ref{Seccion Formulas}. ~Since
$\lim_{\overline{c}\rightarrow\infty}b^{\ast}(\overline{c})=\mu/q,$ the Taylor
expansion of $F_{2}(x,\overline{c})/\overline{c}^{2}$ at $\overline{c}=\infty$
is given by%
\[
\frac{F_{2}(x,\overline{c})}{\overline{c}^{2}}=\frac{2a^{3}xq(xq-2\mu)}%
{\sigma^{4}}+O\left(  \frac{1}{\overline{c}}\right)  .
\]
Since $xq-2\mu<$ $0$, we have $\partial_{\overline{c}}V^{\overline{c}%
}(x,\overline{c})<0$ for $\overline{c}$ large enough and so $\lim
_{\overline{c}\rightarrow\infty}V^{\overline{c}}(x,\overline{c})=x^{+}$ for
$x<\frac{\mu}{q}$.

Take now $x\geq\mu/q>b^{\ast}(\overline{c})$, then%
\[
V^{\overline{c}}(x,\overline{c})=\Big(\frac{\overline{c}}{q}+D(\overline
{c},b^{\ast}(\overline{c})\Big)e^{\theta_{2}(\overline{c})x}%
\]
and so%
\[
\partial_{\overline{c}}V^{\overline{c}}(x,\overline{c})=\dfrac{F_{3}
(x,\overline{c})}{l_{0}(b^{\ast}(\overline{c}),\overline{c})},
\]
where%
\[
F_{3}(x,\overline{c})=%
{\displaystyle\sum\limits_{i=1}^{8}}
\bar{l}_{i}(x,b^{\ast}(\overline{c}),\overline{c})e^{k_{i}(x,b^{\ast
}(\overline{c}),\overline{c})},
\]
$l_{0}(b,\overline{c})$ is defined in (\ref{l0(b,c)}), $\bar{l}_{i}(x,b,c)$
are polynomials on $\theta_{1}(c)$, $\theta_{2}(c)$, $\theta_{1}(ac)$,
$\theta_{2}(ac)$, $\theta_{1}^{\prime}(c)$, $\theta_{2}^{\prime}(c)$,
$\theta_{1}^{\prime}(ac)$, $\theta_{2}^{\prime}(ac)$, $x$, $b$, $c$, $a$ and
$k_{i}(x,b,c)$, $i=1,...,8$, are positive linear combinations of $b\theta
_{1}(ac)$, $b\theta_{2}(ac)$ and ($x-b)\theta_{2}(c)$ as detailed in Section
\ref{Seccion Formulas}. ~Since $\lim_{\overline{c}\rightarrow\infty}b^{\ast
}(\overline{c})=\mu/q$, the Taylor expansion of $F_{3}(x,\overline{c})/$
$\overline{c}^{2}$ at $\overline{c}=\infty$ is given by
\[%
\begin{array}
[c]{lll}%
\frac{F_{3}(x,\overline{c})}{\overline{c}^{2}} & = & \dfrac{2a^{3}}{\sigma
^{4}}\left(  aq^{2}x^{2}-2aq\mu x+\mu^{2}(a-1)\right) \\
&  & +\dfrac{4a^{2}}{3\sigma^{4}}(-3a(2qx-3\mu)(qx-\mu)\mu+5\mu^{3}
+3q\mu\sigma^{2}-a^{2}(qx-\mu)(q^{2}x^{2}-5qx\mu+4\mu^{2}-3q\sigma^{2}
))\dfrac{1}{\overline{c}}\\[0.3cm]
&  & +O\left(  \frac{1}{\overline{c}^{2}}\right)  .
\end{array}
\]
So from $F_{3}(x^{\ast}(\overline{c}),\overline{c})=0$ we obtain that the
Taylor expansion of $x^{\ast}(\overline{c})$ is given by (\ref{xstarConSigma}).

Moreover, we obtain that for $\overline{c}$ large enough, $\partial
_{\overline{c}}V^{\overline{c}}(x,\overline{c})<0$ for $x\in\lbrack\frac{\mu
}{q},\frac{\mu}{q}(1+\frac{1}{\sqrt{a}}))$ and $\partial_{\overline{c}%
}V^{\overline{c}}(x,\overline{c})>0$ for $x>\frac{\mu}{q}(1+\frac{1}{\sqrt{a}%
})$. So, we conclude the result.\hfill$\blacksquare$

\bigskip

\begin{remark}
\label{remm71} \normalfont Note that
\begin{equation}
\label{thisoo}z^{\ast}(\overline{c})-x^{\ast}(\overline{c})=\dfrac{\sigma^{2}%
}{2}\dfrac{1}{\overline{c}}+O(\dfrac{1}{\overline{c}^{2}}),
\end{equation}
so $z^{\ast}(\overline{c})>x^{\ast}(\overline{c})$ for $\overline{c}$ large
enough, and asymptotic equivalence for these two quantities when $\overline
{c}\to\infty$. At the same time, the inequality $z^{\ast}(\overline{c})\geq
x^{\ast}(\overline{c})$ can easily be seen to hold for any $\overline{c}$ from
the following argument: We have
\[%
\begin{array}
[c]{lll}%
V^{\overline{c}}(x,\overline{c})-V^{\overline{c}-h}(x,\overline{c}-h) & = &
V^{\overline{c}}(x,\overline{c})-V^{\overline{c}}(x,\overline{c}%
-h)+V^{\overline{c}}(x,\overline{c}-h)-V^{\overline{c}-h}(x,\overline{c}-h)\\
& \geq & V^{\overline{c}}(x,\overline{c})-V^{\overline{c}}(x,\overline{c}-h),
\end{array}
\]
since, by Proposition \ref{Cambios con cbarra}, $V^{\overline{c}}(x,c)$ is
non-decreasing in $\overline{c}$. So, dividing by $h$ and taking the limit as
$h$ goes to zero, we get
\[
\partial_{\overline{c}}V^{\overline{c}}(x,\overline{c})\geq\left.
V_{c}^{\overline{c}}(x,c)\right\vert _{c=\overline{c}}\text{.}%
\]
Hence $\partial_{\overline{c}}V^{\overline{c}}(z^{\ast}(\overline
{c}),\overline{c})\geq\left.  V_{c}^{\overline{c}}(z^{\ast}(\overline
{c}),c)\right\vert _{c=\overline{c}}=0$ and then the value $x^{\ast}%
(\overline{c})$ where $\partial_{\overline{c}}V^{\overline{c}}(\cdot
,\overline{c})$ changes from negative to positive satisfies $x^{\ast
}(\overline{c})\leq z^{\ast}(\overline{c})$.
\end{remark}

\begin{remark}
\label{Rem7.1}\normalfont One observes from \eqref{zstar} that for very small
values of $a$, the coefficient of $1/\overline{c}$ in the asymptotic expansion
is positive, so that the limit $\mu(1+1/\sqrt{a})/q$ is approached from the
right, whereas for larger values of $a$ that coefficient is negative and the
limit is approached from the left as $\overline{c}$ becomes large, see also
the numerical illustrations in Section \ref{Numerical examples}. \newline It
may also be instructive to derive the higher-order limiting behavior of
$x^{\ast}(\overline{c})$ established in the previous proposition in a direct
way for the deterministic case discussed in Section \ref{sec2}. Concretely,
including one more term in the expansion \eqref{lowerorder} gives
\[
x+\frac{2axq\mu-ax^{2}q^{2}+\mu^{2}(1-a)}{2aq\,\overline{c}}+\frac{\mu
^{3}+3a\mu^{2}(\mu-xq)+a^{2}(xq-4\mu)(\mu-xq)^{2}}{6a^{2}q\;\overline{c}^{2}%
}+O\left(  \frac{1}{\overline{c}^{3}}\right)  ,
\]
and substituting $x=\frac{\mu}{q}(1+\frac{1}{\sqrt{a}})+\frac{a_{0}}%
{\overline{c}}$ (for an $a_{0}\in{\mathbb{R}}$ to be identified) into this
expression gives
\[
\frac{3\sqrt{a^{3}}\,q\mu\,a_{0}+\mu^{3}(2\sqrt{a}+3a-1)}{3a^{2}q\overline
{c}^{3}}+O\left(  \frac{1}{\overline{c}^{4}}\right)  .
\]
This fraction equals zero for $a_{0}=\frac{(1-2\sqrt{a}-3a)\mu^{2}}%
{3\sqrt{a^{3}}q}$, so that we obtain
\[
x^{\ast}(\overline{c})=\frac{\mu}{q}\Big(1+\frac{1}{\sqrt{a}}\Big)+\frac
{(1-2\sqrt{a}-3a)\mu^{2}}{3\sqrt{a^{3}}q\;\overline{c}}+O\left(  \frac
{1}{\overline{c}^{2}}\right)  ,
\]
which exactly corresponds to \eqref{xstarConSigma} for $\sigma
=0$.\newline

\noindent The latter formula shows that in the deterministic case indeed the
limit $\mu(1+1/\sqrt{a})/q$ is approached from the right for $a<1/9$ and from
the left for $a>1/9$ as $\overline{c}\rightarrow\infty$.
\end{remark}

\section{Optimal strategies for $\overline{c}$ large
\label{Optimal strategies for c large enough}}

%
%

In the next proposition, we show that for $\overline{c}$ large enough, there
exists a unique solution of (\ref{Ecuacion diferencial de z0}) with boundary
conditions (\ref{Condicion de Optimo en cbarra}) and that $\overline{\zeta
}^{\prime}<0$ and $\overline{\gamma}^{\prime}>0$ in a neighborhood of
$\overline{c}$. {We emphasize again that for all results in this section, $a$ is assumed to be strictly smaller than 1.}

\begin{proposition}
\label{existencia curvas c grande}For $\overline{c}$ large enough, we can find
$\underline{c}\in\lbrack0,\overline{c})$ such that there exists a unique
solution $(\overline{\gamma}(c),\overline{\zeta}(c))$ of
(\ref{Ecuacion diferencial de z0}) with boundary conditions
(\ref{Condicion de Optimo en cbarra}) in $[\underline{c},\overline{c}]$%
,$\ $and $\overline{\gamma}$ is strictly increasing and $\overline{\zeta}%
\ $\ is strictly decreasing in $[\underline{c},\overline{c}],$ respectively.
\end{proposition}

\textit{Proof.} In order to prove that there exists a unique solution $(\overline{\gamma
}(c),\overline{\zeta}(c))$ of (\ref{Ecuacion diferencial de z0}) in
$[\underline{c},\overline{c}]$ for some $\underline{c}<\overline{c}$, it suffices to show that
\[
C_{11}(b^{\ast}(\overline{c}),z^{\ast}(\overline{c}),\overline{c}%
)\neq0~\text{and }C_{22}(b^{\ast}(\overline{c}),z^{\ast}(\overline
{c}),\overline{c})\neq0\text{ }%
\]
for $\overline{c}$ large enough. Combining (\ref{expansion de bestrella en c}%
) and (\ref{zstar}) with the formulas of $C_{11}(y,z,c)$ and $C_{22}(y,z,c)$
given in Section \ref{Seccion Formulas}, we obtain that
\[
C_{11}(b^{\ast}(\overline{c}),z^{\ast}(\overline{c}),\overline{c}%
)=-\frac{32(1-a)^{2}aq}{\sigma^{10}}\,\overline{c}^{5}+O(\overline{c}^{4}),
\]

\[
C_{22}(b^{\ast}(\overline{c}),z^{\ast}(\overline{c}),\overline{c})e^{(z^{\ast
}(\overline{c})-b^{\ast}(\overline{c}))\theta_{1}(\overline{c})}%
=\frac{32(1-a)q\mu}{\sqrt{a}\,\sigma^{10}}\,\overline{c}^{3}+O(\overline{c}^{2}),
\]
and so
\[
C_{11}(b^{\ast}(\overline{c}),z^{\ast}(\overline{c}),\overline{c})<0~\text{and
}C_{22}(b^{\ast}(\overline{c}),z^{\ast}(\overline{c}),\overline{c})>0\text{ }%
\]
for $\overline{c}$ large enough.

In order to prove that $\overline{\gamma}(c)$ is increasing and $\overline
{\zeta}(c)$ is decreasing in $[\underline{c},\overline{c}]$ for $\overline{c}$
large enough and some $\underline{c}<\overline{c}$, we use the differential
equations (\ref{Ecuacion diferencial de z0}) at $c=\overline{c}$ and the
Taylor expansion of $C_{ij}$ to show that%

\[%
\begin{array}
[c]{lll}%
\gamma^{\prime}(\overline{c}) & = & \dfrac{aq\sigma^{2}+\mu^{2}}{2aq}\,\dfrac
{1}{\overline{c}^{2}}+O(\frac{1}{\overline{c}^{3}}),\\[0.3cm]
\zeta^{\prime}(\overline{c}) & = & -\dfrac{3q\sigma^{4}}{4}\,\dfrac{1}%
{\overline{c}^{4}}+O(\frac{1}{\overline{c}^{5}})
\end{array}
\]
for $\overline{c}$ large enough, so we have the result.\hfill  $\blacksquare$\\

In the following theorem, we show that the value function of the two-curve
strategy $W^{\overline{\gamma},\overline{\zeta}}$ given by the solutions of
$(\overline{\gamma},\overline{\zeta})$ of (\ref{Ecuacion diferencial de z0}) with
boundary condition (\ref{Condicion de Optimo en cbarra}) is the optimal value
function in $[0,\infty)\times\lbrack\underline{c},\overline{c}]$ for
$\overline{c}$ large enough and some $\underline{c}<\overline{c}$. So, the
optimal strategy is a two-curve strategy.

\begin{theorem}
\label{Teorema verificacion asintotico}In the case $\overline{c}$
$>q\sigma^{2}/(2\mu)$, there exists a $\overline{c}$ large enough and some
$\underline{c}<\overline{c}$ such that $W^{\overline{\gamma},\overline{\zeta}}=V$ in
$[0,\infty)\times\lbrack\underline{c},\overline{c}]$.
\end{theorem}

{ \textit{Proof.} By Propositions \ref{existencia curvas c grande} there
exists }$\overline{c}$ large enough and some $\underline{c}<\overline{c}$ such
that $\overline{\zeta}^{\prime}(c)\neq0$ and so, by Proposition \ref{z creciente},
$W^{\overline{\gamma},\overline{\zeta}}${ is (2,1)-differentiable in } $[0,\infty
)\times\lbrack\underline{c},\overline{c}]$. Using Proposition \ref{Verificacion Curva}, in order to prove the result, it is sufficient to show
that%
\begin{equation}
\partial_{x}W^{\overline{\gamma},\overline{\zeta}}(x,c)\geq1\text{ for }x\in
\lbrack0,\overline{\gamma}(c))\text{ and }\partial_{x}W^{\overline{\gamma}%
,\overline{\zeta}}(x,c)\leq1\text{ for }x\in\lbrack\overline{\gamma}(c),\overline
{\zeta}(c)] \label{Desigualdad Wx}%
\end{equation}
and%
\[
\partial_{c}W^{\overline{\gamma},\overline{\zeta}}(x,c)\leq0\text{ for }x\in
\lbrack0,\overline{\zeta}(c))\text{ }%
\]
for $c\in\lbrack\underline{c},\overline{c})$.

We have from Proposition \ref{smooth pasting} $\ $that $\partial_{x}%
(W^{\overline{\gamma},\overline{\zeta}})(\overline{\gamma}(c),c)=1$ for $c\in\lbrack
\underline{c},\overline{c}]$, and the Taylor expansion of $\partial
_{xx}v^{\overline{c}}(x)$ at $\overline{c}=\infty$ is given by%
\[
\left\{
\begin{array}
[c]{ll}%
-\dfrac{q}{ac}+\dfrac{q(qx-2\mu)}{a^{2}c^{2}}+O(\frac{1}{c^3}) & \text{if
}x<b^{\ast}(\overline{c}),\\[0.2cm]
-\dfrac{q}{c}+\dfrac{q(qx-2\mu)}{c^{2}}+O(\frac{1}{c^3}) & \text{if }x\geq
b^{\ast}(\overline{c}),
\end{array}
\right.
\]
which is negative. So, $\partial_{xx}v^{\overline{c}}(x)<0$ for $\overline{c}$
large enough. Since $\partial_{xx}(W^{\overline{\gamma},\overline{\zeta}})(x,c)$ is
continuous, there exists a $\underline{c}<\overline{c}$ such that
$\partial_{xx}(W^{\overline{\gamma},\overline{\zeta}})(x,c)<0$ in $(x,c)$ for
$c\in\lbrack\underline{c},\overline{c}]$ and $x\in\lbrack0,\overline{\zeta}(c)]$.
We conclude that (\ref{Desigualdad Wx}) holds for $c\in\lbrack\underline
{c},\overline{c}]$ and $\overline{c}$ large enough.

Let us show that for $\overline{c}$ large enough and some $\underline{c}<$
$\overline{c}$, it holds that $\partial_{c}W^{\overline{\gamma},\overline{\zeta}}(x,c)\leq0$
for $c\in\lbrack\underline{c},\overline{c}]$ and $0\leq x\leq\overline{\zeta
}(c)$. We prove first that $\partial_{c}W^{\overline{\gamma},\overline{\zeta}}%
(x,c)\leq0$ for $x\in\lbrack\overline{\gamma}(c),\overline{\zeta}(c)].$ We
have that%
\[%
\begin{array}
[c]{lll}%
\partial_{c}W^{\overline{\gamma},\overline{\zeta}}(x,c)=\partial_{c}H^{\overline{\gamma}%
,\overline{\zeta}}(x,c) & = & \frac{d}{dc}(f_{20}(\overline{\gamma}(c),x,c))+\frac
{d}{dc}\left(  f_{21}(\overline{\gamma}(c),x,c)\right)  A^{\overline{\gamma}%
,\overline{\zeta}}(c))+f_{21}(\overline{\gamma}(c),x,c)(A^{\overline{\gamma},\overline{\zeta
}})^{\prime}(c)\\
& = & f_{21}(\overline{\gamma}(c),x,c)\left(  -b_{0}(\overline{\gamma
}(c),x,\overline{\gamma}^{\prime}(c),c)-b_{1}(\overline{\gamma}(c),x,\overline
{\gamma}^{\prime}(c),c)A^{\overline{\gamma},\overline{\zeta}}(c))+(A^{\overline{\gamma}%
,\overline{\zeta}})^{\prime}(c)\right)  ,
\end{array}
\]
and by Lemma \ref{Lema Denominador Positivo}, $f_{21}(y,x,c)>0$ for $x>y$, so
we should prove that
\[
G(x,\overline{c}):=-b_{0}(\gamma(\overline{c}),x,\gamma^{\prime}(\overline
{c}),\overline{c})-b_{1}(\gamma(\overline{c}),x,\gamma^{\prime}(\overline
{c}),\overline{c})A^{\overline{\gamma},\overline{\zeta}}(\overline{c}))+(A^{\overline{\gamma
},\overline{\zeta}})^{\prime}(\overline{c})<0
\]
for $x\in\lbrack\overline{\gamma}(c),\overline{\zeta}(c)].$ By Proposition
\ref{smooth pasting}, $0=\partial_{c}H^{\overline{\gamma},\overline{\zeta}}%
(x,\overline{c})=\partial_{cx}H^{\overline{\gamma},\overline{\zeta}}(x,\overline{c})=0,$
so we have that $G(\overline{\zeta}(\overline{c}),\overline{c})=\partial
_{x}G(\overline{\zeta}(\overline{c}),\overline{c})=0$. Then it is sufficient to
prove that $\partial_{xx}G(x,c)<0$ for $x\in\lbrack\overline{\gamma
}(c),\overline{\zeta}(c)]$. We will first show that $\partial_{xx}%
G(x,\overline{c})<0$ for $x\in\lbrack\overline{\gamma}(\overline{c}%
),\overline{\zeta}(\overline{c})]\ $for $\overline{c}$ large enough, and then
the result follows for $c\in\lbrack\underline{c},\overline{c}]$ for some
$\underline{c}<\overline{c}$ by continuity arguments in a compact set. Using
that $\overline{\gamma}(\overline{c})=b^{\ast}(\overline{c}),~\overline{\zeta
}(\overline{c})=z^{\ast}(\overline{c})$, (\ref{expansion de bestrella en c})
and (\ref{zstar}), we obtain that the Taylor expansion at $\overline{c}%
=\infty$ of
\[
h(x,\overline{c}):=\dfrac{\partial_{xx}G(x,\overline{c})}{e^{-(x-\gamma
(\overline{c}))\theta_{1}(\overline{c})-\gamma(\overline{c})\theta
_{1}(a\overline{c})-\gamma(\overline{c})\theta_{2}(\overline{c})}}%
\]
is given by%

\[
h(x,\overline{c})=\dfrac{2(a(qx-\mu)^{2}-\mu^{2})}{a^{2}q\sigma^{4}}%
+\dfrac{12q\mu\sigma^{2}-10\mu^{3}-4a^{2}(qx-\mu)^{2}(qx+2\mu)+6a(q^{2}%
x^{2}\mu+3\mu^{2})}{3a^{3}q\sigma^{4}}\dfrac{1}{\overline{c}}+O(\dfrac
{1}{\overline{c}^{2}}),
\]
and the Taylor expansion of $h(z^{\ast}(\overline{c}),\overline{c})$ at
$\overline{c}=\infty$ is given by
\[
h(z^{\ast}(\overline{c}),\overline{c})=-\dfrac{2\mu}{a^{\frac{3}{2}}\sigma
^{2}}\dfrac{1}{\overline{c}}+O(\dfrac{1}{\overline{c}^{2}})\text{.}%
\]
Since%
\[
\partial_{x}h(x,\overline{c})=\dfrac{4(a(\overline{c}-qx-\mu)(qx-\mu)+qx\mu
)}{a^{2}\sigma^{4}}\dfrac{1}{\overline{c}}+O(\dfrac{1}{\overline{c}^{2}})
\]
is positive and $h(z^{\ast}(\overline{c}),\overline{c})<0$ for $\overline{c}$
large enough, we conclude that $\partial_{xx}G(x,\overline{c})<0\ $for
$x\in\lbrack\overline{\gamma}(\overline{c}),\overline{\zeta}(\overline{c})]$.
Let us show that for $\overline{c}$ large enough and some $\underline{c}<$
$\overline{c}$, it holds that $\partial_{c}W^{\overline{\gamma},\overline{\zeta}%
}(x,c)\leq0$ for $x\in\lbrack0,\overline{\gamma}(c)]$ and $c\in\lbrack
\underline{c},\overline{c}]$. We can write
\[%
\begin{array}
[c]{lll}%
\partial_{c}W^{\overline{\gamma},\overline{\zeta}}(x,c) & = & \partial_{c}%
(f_{10}(x,c)+f_{11}(x,c)A^{\overline{\gamma},\overline{\zeta}}(c))\\
& = & f_{11}(x,c)\left(  \frac{\partial_{c}f_{10}(x,c)}{f_{11}(x,c)}%
+\frac{\partial_{c}f_{11}(x,c)}{f_{11}(x,c)}A^{\overline{\gamma},\overline{\zeta}%
}(c)+(A^{\overline{\gamma},\overline{\zeta}})^{\prime}(c)\right)
\end{array}
\]
where $f_{11}(x,c)>0$, so we should prove that
\[
G_{1}(x,\overline{c}):=\frac{\partial_{c}f_{10}(x,c)}{f_{11}(x,c)}%
+\frac{\partial_{c}f_{11}(x,c)}{f_{11}(x,c)}A^{\overline{\gamma},\overline{\zeta}%
}(c)+(A^{\overline{\gamma},\overline{\zeta}})^{\prime}(c)<0
\]
for $x\in\lbrack0,\overline{\gamma}(c)].$ We have shown that $\partial
_{c}W^{\overline{\gamma},\overline{\zeta}}(\overline{\gamma}(c),c)<0,$ so we have
$G_{1}(\overline{\gamma}(c),\overline{c})<0$; then, it suffices to prove that
$\partial_{x}G_{1}(x,c)>0$ for $x\in\lbrack0,\overline{\gamma}(c)]$. We will
see first that $\partial_{x}G_{1}(x,\overline{c})>0$ for $x\in\lbrack
0,\overline{\gamma}(\overline{c})]\ $for $\overline{c}$ large enough, then the
result follows for $c\in\lbrack\underline{c},\overline{c}]$ with some
$\underline{c}<\overline{c}$ by continuity arguments in a compact set. Using
that $\overline{\gamma}(\overline{c})=b^{\ast}(\overline{c})$,
(\ref{expansion de bestrella en c}) and (\ref{zstar}), we obtain that the
Taylor expansion at $\overline{c}=\infty$ of
\[
h_{1}(x,\overline{c}):=e^{x\theta_{1}(a\overline{c})}\partial_{x}%
G_{1}(x,\overline{c})
\]
is given by%
\[
h_{1}(x,\overline{c})=\dfrac{x(2\mu-qx)}{\sigma^{2}}\dfrac{1}{\overline{c}%
}+O(\dfrac{1}{\overline{c}^{2}}),
\]
which is positive for $\overline{c}$ large enough and $x\leq b^{\ast
}(\overline{c})<\frac{\mu}{q}<\frac{2\mu}{q}.$ \hfill $\blacksquare$

\bigskip

\section{Numerical examples \label{Numerical examples}}

In this section we will consider some numerical illustrations for the case
$q=0.1$, $\mu=4$ and $\sigma=2$.

\subsection{Bounded Case}

Let us first consider the case with an upper bound $\overline{c}=3$ for the
dividend rate. In this case we are able to derive the optimal value function
and the optimal strategies for the problem with drawdown constraints $a=0.2$,
$a=0.5$ and $a=0.8$. Indeed they are of two-curve type as conjectured in
Remark \ref{Conjetura dos curvas}. The obtained value function and optimal
dividend strategies will then also allow us to compare them with the ones for
the (already previously known) extreme cases $a=0$ (the classical dividend
problem without any constraint) and $a=1$ (the dividend problem with
ratcheting constraint).\newline

\noindent In order to obtain the optimal value functions $V_{a}^{\overline{c}%
}$ for each set of parameters, we proceed as follows:

\begin{enumerate}
\item We check that there exists a unique zero $z^{\ast}(\overline{c})$ of
$C_{0}(b^{\ast}(\overline{c}),\cdot,\overline{c})$ in $(b^{\ast}(\overline
{c}),\infty)$.

\item We obtain the curves $\overline{\gamma}$ and $\overline{\zeta}$ solving
numerically, by the Euler method, the system of ordinary differential
equations (\ref{Ecuacion diferencial de z0}) with boundary condition
(\ref{Condicion de Optimo en cbarra}).

\item We check numerically that the pair $(\overline{\gamma},\overline{\zeta
})$ satisfies condition (\ref{Condicion 2 para z derivable}) for $c\in
\lbrack0,\overline{c}]$. So, by Proposition
\ref{Unicidad de ecuacion diferencial}, we are approximating the unique
solution $(\overline{\gamma},\overline{\zeta})$. We also verify that
$\overline{\zeta}$ is non-decreasing.

\item We check that the function $W^{\overline{\gamma},\overline{\zeta}}$
defined in (\ref{Definicion de Wz}) satisfies the conditions of Theorem
\ref{Caracterizacion Continua}. Hence $W^{\overline{\gamma},\overline{\zeta}}$
is the optimal value function $V_{a}^{\overline{c}}$ and the optimal strategy
is indeed a two-curve strategy given by $(\overline{\gamma},\overline{\zeta}
)\in\mathcal{B}$.\ \ \textit{\ \ \ }
\end{enumerate}

\begin{figure}[tbh]
\begin{center}
\includegraphics[width=6cm]{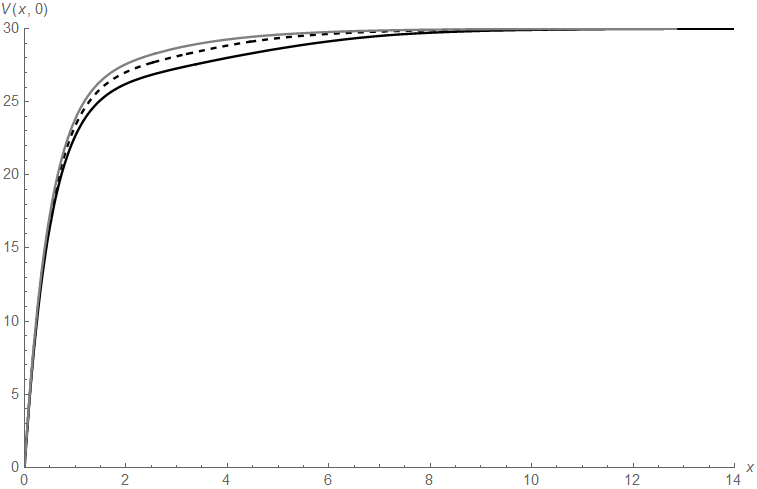}
\end{center}
\caption{{\protect\small $V_{0}^{3}(x,0)$ (gray solid), $V_{0.5}^{3}(x,0)$
(dashed) and $V_{1}^{3}(x,0)$ (black solid) as a function of $x$.}}%
\label{Fig1}%
\end{figure}

Figure \ref{Fig1} depicts the graphs of $V_{a}^{\overline{c}}(x,0)$ with
$\overline{c}=3$ for $a=0$ (no restrictions, gray solid), $a=0.5$ (dashed) and
$a=1$ (ratcheting, black solid) as a function of $x$. One can nicely see how
the drawdown case is - in terms of performance - a compromise between the
unconstrained case and the stronger constraint of ratcheting.

\begin{figure}[tbh]
\begin{center}
\includegraphics[width=7cm]{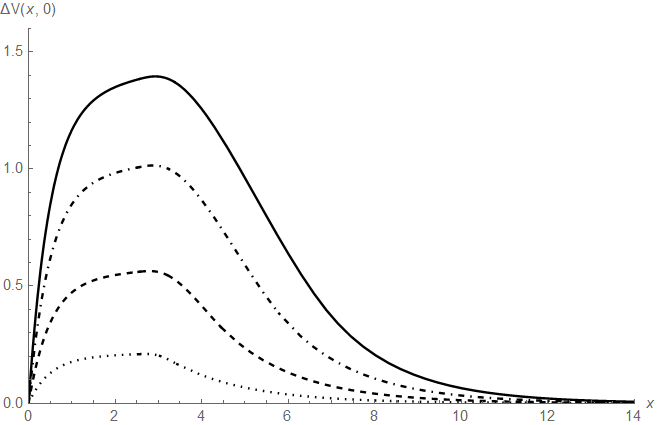}
\end{center}
\caption{{\protect\small $V_{0}^{3} (x)-V_{a}^{3}(x,0)$ for $a=0.2$ (dotted),
$a=0.5$ (dashed), $a=0.8$ (dot-dashed) and $a=1$ (solid).}}%
\label{Fig2}%
\end{figure}\begin{figure}[ptb]
\begin{subfigure}{.32\textwidth}
		\centering
		\includegraphics[width=5.5cm]{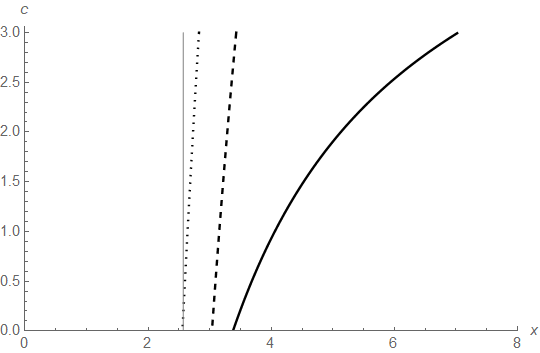}
		\caption{$a=0.2$}
		\label{Fig31}
	\end{subfigure} \begin{subfigure}{.32\textwidth}
		\centering
		\includegraphics[width=5.5cm]{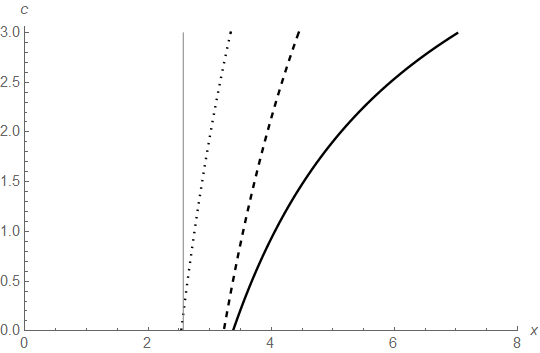}
		\caption{$a=0.5$}
		\label{Fig32}
	\end{subfigure}
\begin{subfigure}{.32\textwidth}
		\centering
		\includegraphics[width=5.5cm]{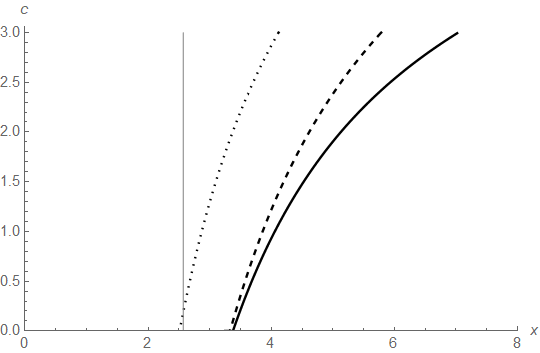}
		\caption{$a=0.8$}
		\label{Fig33}
	\end{subfigure}
\caption{{\protect\small Optimal drawdown curves $(\overline{\gamma}(c),c)$
(dotted) and $(\overline{\zeta}(c),c)\ $(dashed) for $a=0.2,0.5,0.8$, together
with the optimal threshold of the unconstrained problem ($a=0$, solid gray)
and the optimal ratcheting curve $\left(  \overline{\xi}(c),c\right)  $
($a=1$, solid black).}}%
\label{fig:test}%
\end{figure}

In order to see the impact of the drawdown restriction more clearly, in Figure
\ref{Fig2} we plot the difference between $V_{0}^{\overline{c}}(x)$ (the
unconstrained value function) and $V_{a}^{\overline{c}}(x,0)$ as a function of
$x$ for increasingly restrictive drawdown levels $a=0.2$ (dotted), $a=0.5$
(dashed), $a=0.8$ (dot-dashed) and finally $a=1$ (ratcheting, solid). One
observes that in particular for smaller values of $x$, the relaxation of
ratcheting towards the drawdown constraint improves the performance of the
resulting strategy quite a bit, although the relative gap between the
performance of the ratcheting and the unconstrained case is anyway not so big
(cf.\ Figure \ref{Fig1}). The latter speaks in favor of the consideration of
such strategies, as ratcheting and drawdown may be important for shareholders
from a psychological point of view, and the efficiency loss when introducing
these constraints is quite minor. In particular, if for a given initial
surplus level $x$ one has a target efficiency loss one is willing to accept,
results like Figure \ref{Fig2} can help to identify the corresponding drawdown
coefficient $a$ that can still guarantee such a performance. \newline

\begin{figure}[ptb]
\hspace*{3cm}\includegraphics[width=6cm]{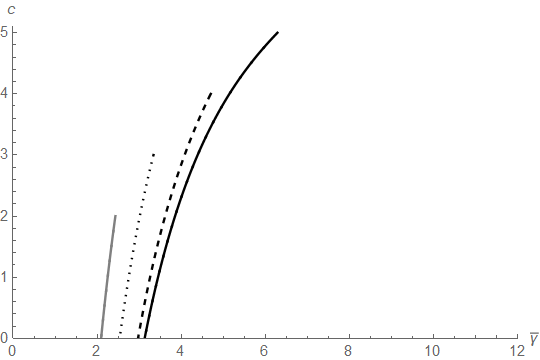} \hspace*{0.5cm}
\includegraphics[width=6cm]{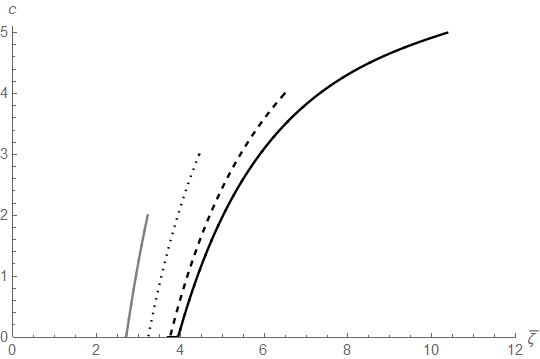} \caption{{\protect\small The
curves $\overline{\gamma}^{\overline{c}}(c)$ (left) and $\overline{\zeta
}^{\overline{c}}(c)$ (right) for $a=0.5$: $\overline{c}=2$ (solid gray),
$\overline{c}=3$ (dotted), $\overline{c}=4=\mu$ (dashed) and $\overline{c}=5$
(solid black).}}%
\label{Figlast}%
\end{figure}

\begin{figure}[ptb]
\begin{center}
\includegraphics[width=5.5cm]{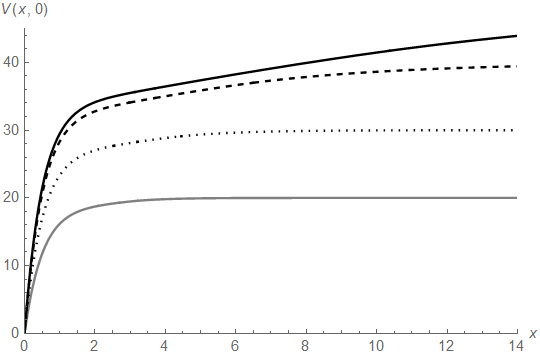}
\end{center}
\caption{{\protect\small $V_{0.5}^{2}(x,0)$ (solid gray), $V_{0.5}^{3}(x,0)$
(dotted), $V_{0.5}^{4}(x,0)$ (dashed) and $V_{0.5}^{5}(x,0)$ (black solid) as
a function of $x$.}}%
\label{Figlast2}%
\end{figure}

\noindent In terms of the nature of the optimal strategy (which indeed turns
out to be of two-curve type), Figures \ref{Fig31}, \ref{Fig32} and \ref{Fig33}
show the optimal drawdown curves $(\overline{\gamma}(c),c)$ (dotted) and
$(\overline{\zeta}(c),c)\ $(dashed) for $a=0.2$, $a=0.5$ and $a=0.8$,
respectively. In all the plots we also depict the optimal threshold of the
unconstrained dividend problem $a=0$ (solid gray) and the optimal ratcheting
curve $\left(  \overline{\xi}(c),c\right)  $ for $a=1$ (solid black). To that
end, recall from Asmussen and Taksar \cite{asmtak} that the optimal threshold
for $a=0$ is given by
\[
\frac{1}{\theta_{1}(0)-\theta_{2}(0)}\log\left(  \frac{\theta_{2}%
(0)\;(\theta_{2}(0)-\theta_{2}(\overline{c}))}{\theta_{1}(0)\;(\theta
_{1}(0)-\theta_{2}(\overline{c}))}\right)  \text{,}%
\]
whereas the optimal strategy in the ratcheting case is given by a one-curve
strategy which is obtained numerically according to the results in
\cite{AAM21}. One can nicely see how the two curves $(\overline{\gamma}(c),c)$
and $(\overline{\zeta}(c),c)$ move towards the right as $a$ increases,
interpolating between the unconstrained and the ratcheting case. Note that the resulting two-curve shapes are somewhat reminiscent of some
figures obtained in Guo and Tomecek \cite{Guo} for other types of singular
control problems, where also a smooth-fit principle was established.\\

\noindent Also notice that the location of these curves can vary considerably as the
maximally allowed dividend rate $\overline{c}$ changes. Figure \ref{Figlast}
depicts $\overline{\gamma}^{\overline{c}}(c)$ and $\overline{\zeta}%
^{\overline{c}}(c)$ for $a=0.5$ for $\overline{c}$ growing from 2 to 5. In
particular, when $\overline{c}$ is larger, the necessary surplus level $x$ to
switch to higher dividend rates is larger as well. Figure \ref{Figlast2} shows
the corresponding value functions for these increasing values of $\overline
{c}$ ($a=0.5$). Recall that while the drawdown constraint is not a major
efficiency loss when compared to the unconstrained case for the same
$\overline{c}$ (cf.\ Figure \ref{Fig1} for the case $\overline{c}=3$), the
size of $\overline{c}$ itself naturally has a considerable impact on the size
of the value function. 

\subsection{Boundary Conditions}

\begin{figure}[ptb]
\begin{subfigure}{.32\textwidth}
		\centering
		\includegraphics[width=5.5cm]{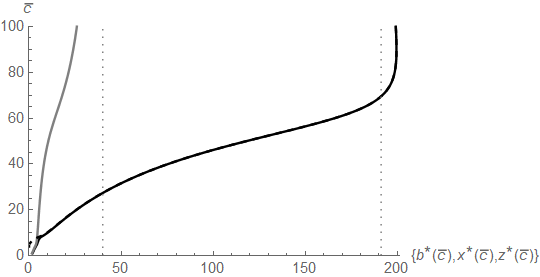}
		\caption{$a=0.07$}
		\label{Fig41}
	\end{subfigure} \begin{subfigure}{.32\textwidth}
		\centering
		\includegraphics[width=5.5cm]{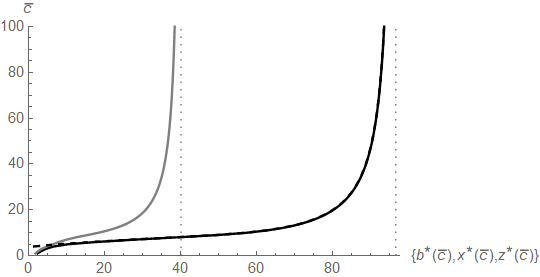}
		\caption{$a=0.5$}
		\label{Fig42}
	\end{subfigure}
\begin{subfigure}{.32\textwidth}
		\centering
		\includegraphics[width=5.5cm]{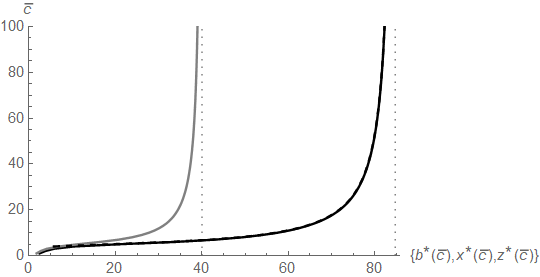}
		\caption{$a=0.8$}
		\label{Fig43}
	\end{subfigure}
\caption{{\protect\small The boundary condition values $b^{\ast}(\overline
{c})$ (grey), $z^{\ast}(\overline{c})$ (solid) and $x^{\ast}(\overline{c})$
(dashed) as a function of $\overline{c}$ for different values of $a$.}}%
\label{fig:test2}%
\end{figure}

\begin{figure}[ptb]
\begin{subfigure}{.32\textwidth}
		\centering
		\includegraphics[width=5.5cm]{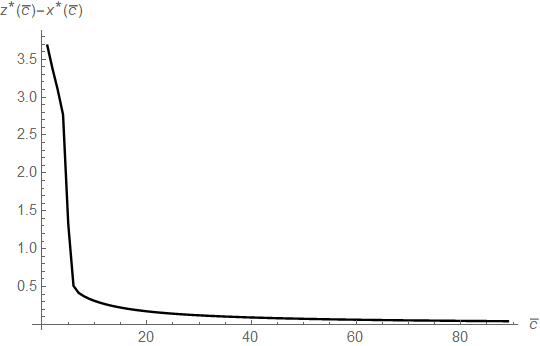}
		\caption{$a=0.07$}
		\label{Fig6a1}
	\end{subfigure} \begin{subfigure}{.32\textwidth}
		\centering
		\includegraphics[width=5.5cm]{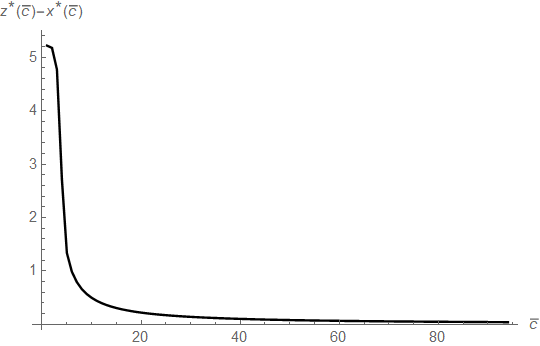}
		\caption{$a=0.5$}
		\label{Fig6a2}
	\end{subfigure}
\begin{subfigure}{.32\textwidth}
		\centering
		\includegraphics[width=5.5cm]{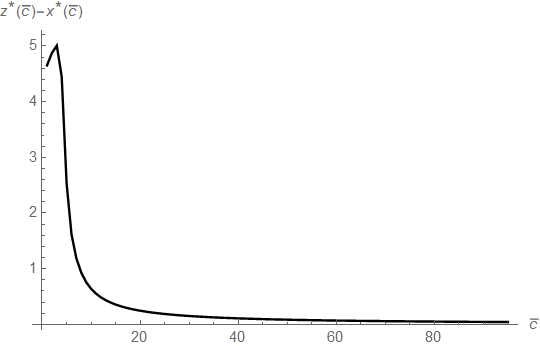}
		\caption{$a=0.8$}
		\label{Fig6a3}
	\end{subfigure}
\caption{{\protect\small The difference $z^{\ast}(\overline{c})-$ $x^{\ast
}(\overline{c})$ as a function of $\overline{c}$ for different values of $a$%
.}}%
\label{fig:test3}%
\end{figure}

Let us now investigate the situation when the maximally allowed dividend rate
$\overline{c}$ becomes large. In addition to $a=0.5$ and $a=0.8$, we now also
consider a smaller drawdown level $a=0.07$ (in order to illustrate the
different monotonicity for small values of $a$, cf.\ Remark \ref{Rem7.1}). One
finds numerically that there exists a unique zero $z^{\ast}(\overline{c})$ of
$C_{0}(b^{\ast}(\overline{c}),\cdot,\overline{c})$ in $(b^{\ast}(\overline
{c}),\infty)$ for any $\overline{c}\geq0$. We also have found that there
exists a unique zero $x^{\ast}(\overline{c})$ in $(0,\infty)$ of
$\partial_{\overline{c}}V^{\overline{c}}(\cdot,\overline{c})$ for
$\overline{c}\geq5.17$ for $a=0.07$; for $\overline{c}\geq3.45$ for $a=0.5$
and for $\overline{c}\geq2.52$ for $a=0.8$. Recall that we have proved in
Propositions \ref{Limite de condiciones en cbarra},
\ref{Limites de condiciones en z barra} and \ref{Derivada de V cbarra} that
$\lim_{\overline{c}\rightarrow\infty}b^{\ast}(\overline{c})=\mu/q$ and
$\lim_{\overline{c}\rightarrow\infty}z^{\ast}(\overline{c})=\lim_{\overline
{c}\rightarrow\infty}x^{\ast}(\overline{c})=\mu(1+1/\sqrt{a})/q$.\newline
Figure \ref{fig:test2} shows the curves of the boundary conditions $(b^{\ast
}(\overline{c}),\overline{c})$, $(z^{\ast}(\overline{c}),\overline{c})$ and
$(x^{\ast}(\overline{c}),\overline{c})$ for $a=0.07$, $a=0.5$ and $a=0.8$
respectively. In the case $a=0.07$ one sees how the limit $\mu(1+1/\sqrt
{a})/q=191.2$ (vertical dotted line) is indeed approached from the right as
$\overline{c}\rightarrow\infty$, whereas for $a=0.5$ and $a=0.8$ the
respective limits $96.57$ and $84.72$ (vertical dotted line) are approached
from the left, cf.\ Remark \ref{Rem7.1}. It is important to keep in mind that
these plots only depict the boundary value for each choice of $\overline{c}$,
and are not to be confused with the optimal drawdown curves in Figure
\ref{fig:test}. Note that $x^{\ast}(\overline{c})$ and $z^{\ast}(\overline
{c})$ are -- already for moderate values of $\overline{c}$ -- almost
identical, with $z^{\ast}(\overline{c})>x^{\ast}(\overline{c})$, see Figure
\ref{fig:test3} for a graph of the difference $z^{\ast}(\overline{c})-$
$x^{\ast}(\overline{c})$ for $a=0.07$, $a=0.5$ and $a=0.8$ respectively. From
the latter, one nicely sees $z^{\ast}(\overline{c})>x^{\ast}(\overline{c})$
(cf.\ Remark \ref{remm71}) as well as the asymptotic equivalence
\eqref{thisoo} of the two quantities. \newline

\begin{figure}[ptb]
\begin{subfigure}{.32\textwidth}
		\centering
		\includegraphics[width=5.5cm]{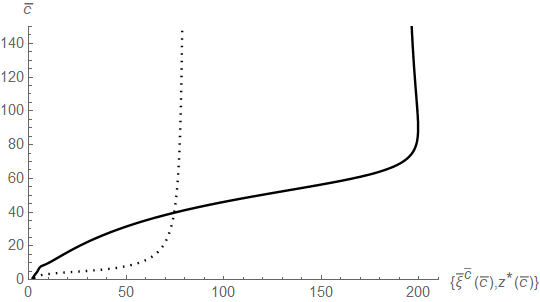}
		\caption{$a=0.07$}
		\label{Fig51}
	\end{subfigure} \begin{subfigure}{.32\textwidth}
		\centering
		\includegraphics[width=5.5cm]{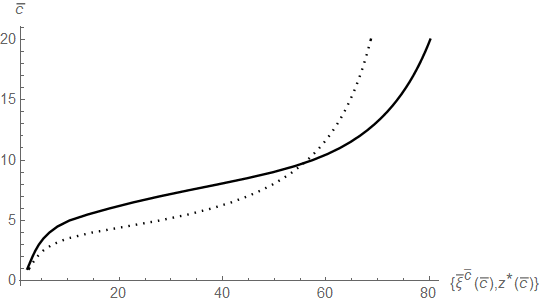}
		\caption{$a=0.5$}
		\label{Fig52}
	\end{subfigure}
\begin{subfigure}{.32\textwidth}
		\centering
		\includegraphics[width=5.5cm]{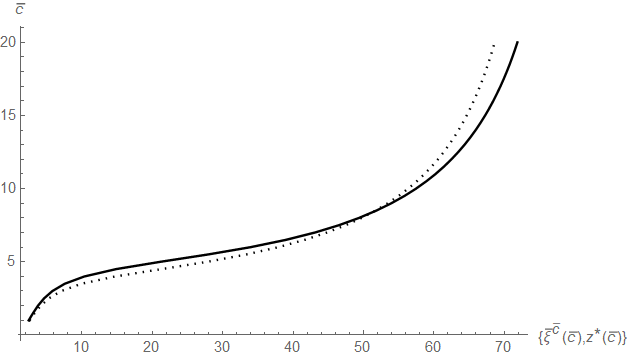}
		\caption{$a=0.8$}
		\label{Fig53}
	\end{subfigure}
\caption{{\protect\small The boundary values $z^{\ast}(\overline{c})\ $
(solid) and the optimal ratcheting boundary value $\overline{\xi}%
^{\overline{c}}(\overline{c})$ (dotted) as $\overline{c}$ grows, for different
values of $a$.}}%
\label{fig:test4}%
\end{figure}

In Figure \ref{fig:test} we saw that the curve $\overline{\zeta}^{\overline
{c}}(c)$ is to the left of the ratcheting curve $\overline{\xi}^{\overline{c}%
}(c)$. At the same time for large values of $\overline{c}$, we know that
$\overline{\zeta}^{\overline{c}}(c)$ must be to the right of $\overline{\xi
}^{\overline{c}}(c)$, as
\[
\lim_{\overline{c}\longrightarrow\infty}\overline{\xi}^{\overline{c}%
}(\overline{c})=\frac{2\mu}{q}<\frac{\mu}{q}\left(  1+\frac{1}{\sqrt{a}%
}\right)  =\lim_{\overline{c}\longrightarrow\infty}z^{\ast}(\overline
{c})\text{.}%
\]
It is therefore of interest to see when this crossing for the limiting value
takes place. Figure \ref{fig:test4} depicts $z^{\ast}(\overline{c})$ (solid)
and $\overline{\xi}^{\overline{c}}(\overline{c})$ (dotted) for $a=0.07$,
$a=0.5$ and $a=0.8$ respectively. We see that indeed for $\overline{c}$ small,
$z^{\ast}(\overline{c})<\overline{\xi}^{\overline{c}}(\overline{c})$ and for
$\overline{c}$ large, $z^{\ast}(\overline{c})>\overline{\xi}^{\overline{c}%
}(\overline{c})$. Moreover, we obtain numerically that the intersection point
of the curves of $z^{\ast}(\overline{c})\ $and $\overline{\xi}^{\overline{c}%
}(\overline{c})$ occurs at $\overline{c}=39.70$ for $a=0.07$, at $\overline
{c}=9.74$ for $a=0.5$, and at $\overline{c}=8.37$ for $a=0.8$ for the given
set of parameters. That is, on from these values of $\overline{c}$, the
possibility of the drawdown increases the value of surplus on from which one
starts to pay the maximal dividend rate, when compared to pure ratcheting, and
it is intuitive that the difference is less pronounced as $a$ increases.

\section{Conclusions}
\label{secconcl} In this paper we addressed the problem of optimal dividends
under a drawdown constraint. We showed that the value function can be
expressed as the unique viscosity solution of a respective two-dimensional
Hamilton-Jacobi-Bellman equation and derived conditions under which the
optimal strategy is of a two-curve form. We conjecture that these conditions
are in fact always fulfilled {and -- using a smooth-fit principle -- could prove it for large values of current and maximal dividend rate $c$ and $\overline{c}$, respectively. For concrete numerical examples, we also proved the optimality of two-curve strategies numerically for small values of $c$ and $\overline{c}$, and showed how to identify
the resulting optimal curves,} which turns out to be a very challenging and
technical task, involving the numerical solution of a highly involved system
of ordinary differential equations and its boundary conditions. We illustrate
how this can be concretely implemented for a moderate size of $\overline{c}$ ;
for high values of $\overline{c}$ this is difficult numerically because the
formulas involve algebraic sums with terms with exponentials with very large
exponents and the computations require very high numerical precision. We
furthermore showed that, when $\overline{c}$ tends to infinity, the curves
converge to a finite limit, the size of which follows a surprisingly simple
and intriguing formula in terms of the square-root of the drawdown percentage
$a$, and irrespective of the size of the volatility parameter $\sigma$. The
latter fact also allowed to get some intuition on the nature of this limit
from the deterministic limit case $\sigma=0$.\newline Altogether, this paper
for the first time explicitly addressed a drawdown constraint for a control
problem in this context, and it turned out that the resulting strategies
smoothly interpolate between the unconstrained problem and the situation with
ratcheting constraints, allowing to get some quantitative insight in the
efficiency gain when relaxing the ratcheting. It will be interesting to see
whether other dividend -- and more generally control -- problems can be
extended in a similar way. In particular, extending the results of the paper
from the Brownian risk model to a compound Poisson surplus process may be an
interesting endeavour, which would lead to a relaxation of the ratcheting
problem studied in \cite{AAM}. Another future direction of research may be to
extend the approach of this paper to incorporate constraints on the dividend
rate in terms of an average of its previous values, for instance along the
lines of Angoshtari et al.\ \cite{bayr2}.

\section{Appendix: Some Formulas \label{Seccion Formulas}}

In the following, we state some definitions and formulas referred to earlier
in the paper in a compact way.
\begingroup
\allowdisplaybreaks
\begin{align*}
d(y,z,c)  &  =e^{y\theta_{2}(ac)+(z-y)\theta_{1}(c)}\theta_{2}%
(c)-e^{(z-y)\theta_{1}(c)+y\theta_{1}(ac)}\theta_{2}(c)+e^{(z-y)\theta
_{2}(c)+y\theta_{2}(ac)}\theta_{2}(ac)\\
&  -e^{y\theta_{2}(ac)+(z-y)\theta_{1}(c)}\theta_{2}(ac)+e^{(z-y)\theta
_{2}(c)}\left(  -e^{y\theta_{2}(ac)}+e^{y\theta_{1}(ac)}\right)  \theta
_{1}(c)\\
&  +e^{y\theta_{1}(ac)}\left(  -e^{(z-y)\theta_{2}(c)}+e^{(z-y)\theta_{1}%
(c)}\right)  \theta_{1}(ac),
\end{align*}

\begin{align*}
b_{00}(y,z,c)  &  =ae^{(z-y)\theta_{2}(c)+y\theta_{2}(ac)}\theta
_{2}(ac)\left(  \theta_{2}(c)-\theta_{1}(c)\right)  -e^{(z-y)\theta_{1}%
(c)}\theta_{2}(c)\left(  \theta_{1}(c)-\theta_{2}(c)\right) \\
&  -ae^{(z-y)\theta_{1}(c)}\left(  -1+e^{y\theta_{2}(ac)}\right)  \theta
_{2}(c)\left(  \theta_{1}(c)-\theta_{2}(c)\right) \\
&  +ae^{y\theta_{2}(ac)+(z-y)\theta_{1}(c)}\theta_{2}(ac)\left(  \theta
_{1}(c)-\theta_{2}(c)\right)  +e^{(z-y)\theta_{2}(c)}\theta_{1}(c)\left(
\theta_{1}(c)-\theta_{2}(c)\right) \\
&  +ae^{(z-y)\theta_{2}(c)}\left(  -1+e^{y\theta_{2}(ac)}\right)  \theta
_{1}(c)\left(  \theta_{1}(c)-\theta_{2}(c)\right) \\
&  -\left(  \theta_{1}(c)-\theta_{2}(c)\right)  {}^{2}-e^{(z-y)\theta_{1}%
(c)}\left(  c+ac\left(  -1+e^{y\theta_{2}(ac)}\right)  \right)  \left(
\theta_{1}(c)-\theta_{2}(c)\right)  \theta_{2}^{\prime}(c)\\
&  +e^{(z-y)\theta_{2}(c)}(z-y)\left(  \theta_{1}(c)-\theta_{2}(c)\right)
\left(  -ace^{y\theta_{2}(ac)}\theta_{2}(ac)+c\left(  1+a\left(
-1+e^{y\theta_{2}(ac)}\right)  \right)  \theta_{1}(c)\right)  \theta
_{2}^{\prime}(c)\\
&  +a^{2}ce^{(z-y)\theta_{2}(c)+y\theta_{2}(ac)}\left(  \theta_{2}%
(c)-\theta_{1}(c)\right)  \theta_{2}^{\prime}(ac)+a^{2}ce^{(z-y)\theta
_{2}(c)+y\theta_{2}(ac)}y\theta_{2}(ac)\left(  \theta_{2}(c)-\theta
_{1}(c)\right)  \theta_{2}^{\prime}(ac)\\
&  +a^{2}ce^{y\theta_{2}(ac)+(z-y)\theta_{1}(c)}\left(  \theta_{1}%
(c)-\theta_{2}(c)\right)  \theta_{2}^{\prime}(ac)-a^{2}ce^{y\theta
_{2}(ac)+(z-y)\theta_{1}(c)}y\theta_{2}(c)\left(  \theta_{1}(c)-\theta
_{2}(c)\right)  \theta_{2}^{\prime}(ac)\\
&  +a^{2}ce^{y\theta_{2}(ac)+(z-y)\theta_{1}(c)}y\theta_{2}(ac)\left(
\theta_{1}(c)-\theta_{2}(c)\right)  \theta_{2}^{\prime}(ac)+a^{2}%
ce^{(z-y)\theta_{2}(c)+y\theta_{2}(ac)}y\theta_{1}(c)\left(  \theta
_{1}(c)-\theta_{2}(c)\right)  \theta_{2}^{\prime}(ac)\\
&  +e^{(z-y)\theta_{2}(c)}\left(  c+ac\left(  -1+e^{y\theta_{2}(ac)}\right)
\right)  \left(  \theta_{1}(c)-\theta_{2}(c)\right)  \theta_{1}^{\prime}(c)\\
&  +ce^{(z-y)\theta_{1}(c)}(y-z)\left(  \theta_{2}(c)+a\left(  -1+e^{y\theta
_{2}(ac)}\right)  \theta_{2}(c)-ae^{y\theta_{2}(ac)}\theta_{2}(ac)\right)
\left(  \theta_{1}(c)-\theta_{2}(c)\right)  \theta_{1}^{\prime}(c)\\
&  +e^{(z-y)\theta_{1}(c)}\left(  c\left(  1+a\left(  -1+e^{y\theta_{2}%
(ac)}\right)  \right)  \theta_{2}(c)-ace^{y\theta_{2}(ac)}\theta
_{2}(ac)\right)  \left(  -\theta_{2}^{\prime}(c)+\theta_{1}^{\prime}(c)\right)
\\
&  +ce^{(z-y)\theta_{2}(c)}\left(  ae^{y\theta_{2}(ac)}\theta_{2}(ac)+\left(
-1+a-ae^{y\theta_{2}(ac)}\right)  \theta_{1}(c)\right)  \left(  -\theta
_{2}^{\prime}(c)+\theta_{1}^{\prime}(c)\right)  ,
\end{align*}
\endgroup
\begin{align*}
b_{01}(y,z,c)  &  =c\left(  \theta_{1}(c)-\theta_{2}(c)\right) \\
&  \left(  ae^{y\theta_{2}(ac)}\theta_{2}(ac)\left(  -\theta_{2}(c)+\theta
_{2}(ac)\right)  +\left(  \theta_{2}(c)+a\left(  -1+e^{y\theta_{2}%
(ac)}\right)  \theta_{2}(c)-ae^{y\theta_{2}(ac)}\theta_{2}(ac)\right)
\theta_{1}(c)\right)  ,
\end{align*}

\bigskip%

\begin{align*}
b_{10}(y,z,c)  &  =e^{(-y+z)\theta_{2}(c]}\left(  e^{y\theta_{1}%
(ac)}-e^{y\theta_{2}(ac)}\right)  \left(  -\theta_{1}(c)+\theta_{2}(c)\right)
\theta_{1}^{\prime}(c)\\
&  -e^{(-y+z)\theta_{1}(c)}(y-z)\left(  \theta_{1}(c)-\theta_{2}(c)\right)
\left(  -e^{y\theta_{1}(ac)}\theta_{1}(ac)+\left(  e^{y\theta_{1}%
(ac)}-e^{y\theta_{2}(ac)}\right)  \theta_{2}(c)+e^{y\theta_{2}(ac)}\theta
_{2}(ac)\right)  \theta_{1}^{\prime}(c)\\
&  -ae^{(-y+z)\theta_{1}(c)+y\theta_{1}(ac)}\left(  \theta_{1}(c)-\theta
_{2}(c)\right)  \theta_{1}^{\prime}(ac)+ae^{y\theta_{1}(ac)+(-y+z)\theta
_{2}(c)}\left(  \theta_{1}(c)-\theta_{2}(c)\right)  \theta_{1}^{\prime}(ac)\\
&  -ae^{(-y+z)\theta_{1}(c)+y\theta_{1}(ac)}y\theta_{1}(ac)\left(  \theta
_{1}(c)-\theta_{2}(c)\right)  \theta_{1}^{\prime}(ac)\\
&  +ae^{y\theta_{1}(ac)+(-y+z)\theta_{2}(c)}y\theta_{1}(ac)\left(  \theta
_{1}(c)-\theta_{2}(c)\right)  \theta_{1}^{\prime}(ac)\\
&  -e^{(-y+z)\theta_{2}(c)}\left(  \left(  -e^{y\theta_{1}(ac)}+e^{y\theta
_{2}(ac)}\right)  \theta_{1}(c)+e^{y\theta_{1}(ac)}\theta_{1}(ac)-e^{y\theta
_{2}(ac)}\theta_{2}(ac)\right)  \left(  \theta_{1}^{\prime}(c)-\theta
_{2}^{\prime}(c)\right) \\
&  +e^{(-y+z)\theta_{1}(c)}\left(  e^{y\theta_{1}(ac)}\theta_{1}(ac)+\left(
-e^{y\theta_{1}(ac)}+e^{y\theta_{2}(ac)}\right)  \theta_{2}(c)-e^{y\theta
_{2}(ac)}\theta_{2}(ac)\right)  \left(  \theta_{1}^{\prime}(c)-\theta
_{2}^{\prime}(c)\right) \\
&  +e^{(-y+z)\theta_{1}(c)}\left(  e^{y\theta_{1}(ac)}-e^{y\theta_{2}%
(ac)}\right)  \left(  \theta_{1}(c)-\theta_{2}(c)\right)  \theta_{2}^{\prime
}(c)\\
&  +e^{(-y+z)\theta_{2}(c)}(-y+z)\left(  \theta_{1}(c)-\theta_{2}(c)\right)
\left(  \left(  -e^{y\theta_{1}(ac)}+e^{y\theta_{2}(ac)}\right)  \theta
_{1}(c)+e^{y\theta_{1}(ac)}\theta_{1}(ac)-e^{y\theta_{2}(ac)}\theta
_{2}(ac)\right)  \theta_{2}^{\prime}(c)\\
&  +ae^{(-y+z)\theta_{1}(c)+y\theta_{2}(ac)}\left(  \theta_{1}(c)-\theta
_{2}(c)\right)  \theta_{2}^{\prime}(ac)\\
&  +ae^{(-y+z)\theta_{2}(c)+y\theta_{2}(ac)}\left(  -\theta_{1}(c)+\theta
_{2}(c)\right)  \theta_{2}^{\prime}(ac)+ae^{(-y+z)\theta_{1}(c)+y\theta
_{2}(ac)}y\left(  \theta_{1}(c)-\theta_{2}(c)\right)  \theta_{2}(ac)\theta
_{2}^{\prime}(ac)\\
&  +ae^{(-y+z)\theta_{2}(c)+y\theta_{2}(ac)}y\left(  -\theta_{1}(c)+\theta
_{2}(c)\right)  \theta_{2}(ac)\theta_{2}^{\prime}(ac)\\
&  -ae^{(-y+z)\theta_{2}(c)}y\theta_{1}(c)\left(  \theta_{1}(c)-\theta
_{2}(c)\right)  \left(  e^{y\theta_{1}(ac)}\theta_{1}^{\prime}(ac)-e^{y\theta
_{2}(ac)}\theta_{2}^{\prime}(ac)\right) \\
&  +ae^{(-y+z)\theta_{1}(c)}y\left(  \theta_{1}(c)-\theta_{2}(c)\right)
\theta_{2}(c)\left(  e^{y\theta_{1}(ac)}\theta_{1}^{\prime}(ac)-e^{y\theta
_{2}(ac)}\theta_{2}^{\prime}(ac)\right)
\end{align*}

and%

\begin{align*}
b_{11}(y,z,c)  &  =-(\text{$\theta_{1}(c)$}-\text{$\theta_{2}(c)$})\\
&  \left(  e^{y\text{$\theta_{1}(ac)$}}(\theta_{1}(ac)-\theta_{1}%
(c))(\theta_{1}(ac)-\theta_{2}(c))+e^{y\text{$\theta_{2}(ac)$}}(\theta
_{2}(c)-\theta_{2}(ac))(\theta_{2}(ac)-\theta_{1}(c))\right)
\end{align*}

\[
f_{10}(x,c)=\frac{ca}{q}(1-e^{\theta_{2}(ac)x}),
\]%
\[
f_{11}(x,c)=e^{\theta_{1}(ac)x}-e^{\theta_{2}(ac)x},
\]

\begin{align*}
f_{20}(y,x,c)  &  =\frac{c}{q(\text{$\theta_{2}(c)$}-\text{$\theta_{1}(c)$}%
)}\\
&  (\text{$\theta_{2}(c)+(a-1)$}e^{\theta_{1}(c)(x-y)}\text{$\theta_{2}(c)+a$%
}e^{y\theta_{2}(ac)}(-e^{\theta_{2}(c)(x-y)}\theta_{2}(ac)+e^{\theta
_{1}(c)(x-y)}(\theta_{2}(ac)-\theta_{2}(c)))\\
&  +\theta_{1}(c)(-1+e^{\theta_{2}(c)(x-y)}(1+a(e^{y\theta_{2}(ac)}-1)))),
\end{align*}

\begin{align*}
f_{21}(y,x,c)  &  =\frac{1}{\text{$\theta_{1}(c)$}-\text{$\theta_{2}(c)$}}\\
&  (e^{y\text{$\theta_{2}(ac)$}}(e^{\theta_{1}(c)(x-y)}(\text{$\theta
_{2}(c)-\theta_{2}(ac))+$}e^{\theta_{2}(c)(x-y)}(\theta_{2}(ac)-\text{$\theta
_{1}(c)$}))\\
&  +e^{y\text{$\theta_{1}(ac)$}}(e^{\theta_{2}(c)(x-y)}(\text{$\theta
_{1}(c)-\theta_{1}(ac))+$}e^{\theta_{1}(c)(x-y)}(\theta_{1}(ac)-\text{$\theta
_{2}(c)$}))\\
&  =\frac{d(y,x,c)}{\text{$\theta_{1}(c)$}-\text{$\theta_{2}(c)$}}%
\end{align*}

\[
C_{0}(y,z,c)=b_{11}(y,c)\partial_{z}\left(  \frac{b_{00}(y,z,c)}%
{d(y,z,c)}\right)  -b_{01}(y,c)\partial_{z}\left(  \frac{b_{10}(y,z,c)}%
{d(y,z,c)}\right)
\]

\[
C_{11}(y,z,c)=b_{11}(y,c)\partial_{y}\left(  \frac{(e^{(z-y)\theta_{1}%
(c)}-e^{(z-y)\theta_{2}(c)})b_{01}(y,c)}{d(y,z,c)}\right)  -~b_{01}%
(y,c)\partial_{y}\left(  \frac{(e^{(z-y)\theta_{1}(c)}-e^{(z-y)\theta_{2}%
(c)})b_{11}(y,c)}{d(y,z,c)}\right)
\]

\[
C_{10}(y,z,c)=b_{01}(y,c)\partial_{y}\left(  \frac{b_{10}(y,z,c)}%
{d(y,z,c)}\right)  -b_{11}(y,c)\partial_{y}\left(  \frac{b_{00}(y,z,c)}%
{d(y,z,c)}\right)
\]%
\[
C_{21}(y,z,c)=\partial_{y}C_{0}(y,z,c)=\partial_{y}\left(  b_{11}%
(y,c)\partial_{z}\left(  \frac{b_{00}(y,z,c)}{d(y,z,c)}\right)  -b_{01}%
(y,c)\partial_{z}\left(  \frac{b_{10}(y,z,c)}{d(y,z,c)}\right)  \right)
\]%
\[
C_{22}(y,z,c)=\partial_{z}C_{0}(y,z,c)=\partial_{z}\left(  b_{11}%
(y,c)\partial_{z}\left(  \frac{b_{00}(y,z,c)}{d(y,z,c)}\right)  -b_{01}%
(y,c)\partial_{z}\left(  \frac{b_{10}(y,z,c)}{d(y,z,c)}\right)  \right)
\]

\[
C_{20}(y,z,c)=-\partial_{c}C_{0}(y,z,c)=-\partial_{c}\left(  b_{11}%
(y,c)\partial_{z}\left(  \frac{b_{00}(y,z,c)}{d(y,z,c)}\right)  -b_{01}%
(y,c)\partial_{z}\left(  \frac{b_{10}(y,z,c)}{d(y,z,c)}\right)  \right)
\]

\[
\left(
\begin{array}
[c]{c}%
g_{_{1}}(y,z,c)\\
g_{2}(y,z,c)\\
g_{3}(y,z,c)\\
g_{4}(y,z,c)\\
g_{_{5}}(y,z,c)\\
g_{_{6}}(y,z,c)\\
g_{_{7}}(y,z,c)\\
g_{_{8}}(y,z,c)\\
g_{_{9}}(y,z,c)\\
g_{_{10}}(y,z,c)\\
g_{_{11}}(y,z,c)\\
g_{_{12}}(y,z,c)\\
g_{_{13}}(y,z,c)\\
g_{_{14}}(y,z,c)\\
g_{_{15}}(y,z,c)\\
g_{_{16}}(y,z,c)
\end{array}
\right)  =\left(
\begin{array}
[c]{cccc}%
0 & 1 & 0 & 2\\
0 & 2 & 0 & 2\\
1 & 0 & 0 & 2\\
1 & 1 & 0 & 2\\
0 & 2 & 1 & 2\\
1 & 0 & 1 & 1\\
1 & 1 & 1 & 1\\
1 & 1 & 1 & 2\\
0 & 1 & 2 & 0\\
0 & 2 & 2 & 0\\
0 & 2 & 2 & 1\\
1 & 0 & 2 & 0\\
1 & 1 & 2 & 0\\
1 & 1 & 2 & 1\\
0 & 2 & 1 & 1\\
0 & 1 & 1 & 1
\end{array}
\right)  \cdot\left(
\begin{array}
[c]{c}%
(z-y)\theta_{1}(c)\\
(z-y)\theta_{2}(c)\\
y\theta_{1}(ac)\\
y\theta_{2}(ac)
\end{array}
\right)
\]

\[
\left(
\begin{array}
[c]{c}%
h_{_{1}}(x,b,c)\\
h_{2}(x,b,c)\\
h_{3}(x,b,c)\\
h_{4}(x,b,c)\\
h_{_{5}}(x,b,c)\\
h_{_{6}}(x,b,c)\\
h_{_{7}}(x,b,c)\\
h_{_{8}}(x,b,c)\\
h_{_{9}}(x,b,c)\\
h_{_{10}}(x,b,c)\\
h_{_{11}}(x,b,c)
\end{array}
\right)  =\left(
\begin{array}
[c]{cccc}%
2 & 0 & 0 & 0\\
1 & 1 & 0 & 0\\
0 & 0 & 2 & 0\\
0 & 0 & 1 & 1\\
1 & 0 & 1 & 0\\
0 & 1 & 1 & 0\\
1 & 1 & 1 & 0\\
0 & 1 & 2 & 0\\
1 & 0 & 0 & 1\\
2 & 0 & 0 & 1\\
1 & 0 & 1 & 1
\end{array}
\right)  \cdot\left(
\begin{array}
[c]{c}%
b\theta_{1}(ac)\\
x\theta_{1}(ac)\\
b\theta_{2}(ac))\\
x\theta_{2}(ac)
\end{array}
\right)
\]

\[
\left(
\begin{array}
[c]{c}%
k_{_{1}}(x,b,c)\\
k_{2}(x,b,c)\\
k_{3}(x,b,c)\\
k_{4}(x,b,c)\\
k_{_{5}}(x,b,c)\\
k_{_{6}}(x,b,c)\\
k_{_{7}}(x,b,c)\\
k_{_{8}}(x,b,c)
\end{array}
\right)  =\left(
\begin{array}
[c]{ccc}%
2 & 0 & 0\\
2 & 0 & 1\\
0 & 2 & 0\\
1 & 1 & 0\\
1 & 1 & 1\\
2 & 1 & 1\\
0 & 2 & 1\\
1 & 2 & 1
\end{array}
\right)  \cdot\left(
\begin{array}
[c]{c}%
b\theta_{1}(ac)\\
b\theta_{2}(ac)\\
(x-b)\theta_{2}(c))
\end{array}
\right)  .
\]


\begin{thebibliography}{99}                                                                                               %




\bibitem {AAM21}Albrecher, H., Azcue, P. and Muler N. (2022). Optimal
ratcheting of dividends in a Brownian model. \textit{SIAM J. Financial Math.},
to appear.

\bibitem {AAM}Albrecher, H., Azcue, P. and Muler N. (2020), Optimal ratcheting
of dividends in insurance.\textit{\ SIAM Journal on Control and Optimization},
58(4), 1822--1845.

\bibitem {ABB}Albrecher H., B\"{a}uerle N. and Bladt M. (2018). Dividends:
From refracting to ratcheting.\textit{ Insurance Math. Econom.} \textbf{83,} 47--58.

\bibitem {AT}Albrecher, H. and Thonhauser, S. (2009). Optimality results for
dividend problems in insurance. \textit{RACSAM-Revista de la Real Academia de
Ciencias Exactas, Fisicas y Naturales. Serie A. Matematicas} \textbf{103},
No.2, 295--320.

\bibitem {bayr}Angoshtari, B., Bayraktar, E. and Young, V.R. (2019) Optimal
dividend distribution under drawdown and ratcheting constraints on dividend
rates. \textit{SIAM Journal on Financial Mathematics} \textbf{10}, 2, 547--577.

\bibitem {bayr2}Angoshtari, B., Bayraktar, E. and Young, V.R. (2020) Optimal
Consumption under a Habit-Formation Constraint. \textit{SIAM Journal on
Financial Mathematics} \textbf{13}, 1, 321--352.

\bibitem {asmtak}Asmussen, S.and Taksar, M. (1997). Controlled diffusion
models for optimal dividend pay-out. \textit{Insurance: Mathematics and
Economics} \textbf{20}, 1, 1-15.

\bibitem {avanzi}Avanzi, B. (2009). Strategies for dividend distribution: A
review. \textit{North American Actuarial Journal} \textbf{13}, 2, 217-251.



\bibitem {AM Libro}Azcue P. and Muler N. (2014). \textit{Stochastic
Optimization in Insurance: a Dynamic Programming Approach. }Springer Briefs in
Quantitative Finance. Springer.

\bibitem {BorodinSalminem}Borodin, A. N. and Salminen, P. (2002).
\textit{Handbook of Brownian Motion---Facts and Formulae}. 2nd ed.
Birkh\"{a}user, Basel.

\bibitem {Brinker}Brinker, L.V. (2021). \textit{Stochastic Optimisation of
Drawdowns via Dynamic Reinsurance Controls.} Doctoral Dissertation,
Universit\"at zu K\"oln.

\bibitem {BrinkerS}Brinker, L.V. and Schmidli, H. (2022) Optimal discounted
drawdowns in a diffusion approximation under proportional reinsurance.
\textit{Journal of Applied Probability}, to appear.

\bibitem {Claisse}Claisse, J., Talay, D. and Tan, X. (2016). A pseudo-Markov
property for controlled diffusion processes. \textit{SIAM Journal on Control
and Optimization}, \textbf{54}, 2, 1017-1029.




\bibitem {chen}Chen, X., Landriault, D., Li, B. and Li, D. (2015). On
minimizing drawdown risks of lifetime investments. \textit{Insurance:
Mathematics and Economics} \textbf{65}, 46--54.





\bibitem {defin}De Finetti, B. (1957). Su un'Impostazione Alternativa della
Teoria Collettiva del Rischio. \textit{Transactions of the 15th Int. Congress
of Actuaries} \textbf{2}, 433--443.

\bibitem {dybvig}Dybvig, P.H. (1995). Dusenberry's ratcheting of consumption:
optimal dynamic consumption and investment given intolerance for any decline
in standard of living.\textit{The Review of Economic Studies} \textbf{62}, 2, 287--313.

\bibitem {thon}Eisenberg, J., Grandits, P. and Thonhauser, S. (2014). Optimal
consumption under deterministic income. \textit{Journal of Optimization Theory
and Applications} \textbf{160},1, 255-279.

\bibitem {elie}Elie, R. and Touzi, N. (2008). Optimal lifetime consumption and
investment under a drawdown constraint. \textit{Finance and Stochastics}
\textbf{12}, 3, 299--330.

\bibitem {Ger69}Gerber, H.U. (1969). {Entscheidungskriterien fuer den
zusammengesetzten {P}oisson-Prozess}. \textit{Schweiz. Aktuarver. Mitt.}
(1969), No.1, 185--227.

\bibitem {Ger72}Gerber, H. U. (1972). Games of economic survival with
discrete-and continuous-income processes. Operations Research \textbf{20}, 1, 37--45.

\bibitem {Gerber2004}Gerber, H. U. and Shiu, E.S.W. (2004). Optimal dividends:
analysis with Brownian motion. \textit{North American Actuarial Journal},
\textbf{8}, 1, 1--20.







\bibitem {Guo}Guo, X. and Tomecek, P. (2009). A class of singular control
problems and the smooth fit principle. \textit{SIAM Journal on Control and
Optimization} \textbf{47}, 6, 3076-3099.

\bibitem {Jeanblanc}Jeanblanc-Picqu\'e, M. and Shiryaev, A. (1995)
Optimization of the flow of dividends. \textit{Uspekhi Mat. Nauk} \textbf{50},
2(302), 25--46.



\bibitem {Kar}Kardaras, C., Obloj, J. and Platen, E. (2017). The num\'eraire
property and long-term growth optimality for drawdown-constrained investments.
\textit{Mathematical Finance} \textbf{27}, 1, 68--95.

\bibitem {Land}Landriault, D., Li, B. and Zhang, H. (2017). On magnitude,
asymptotics and duration of drawdowns for L\'evy models. \textit{Bernoulli}
\textbf{23}, 1, 432--458.

\bibitem {LoeffenRenaud}Loeffen, R.L. and Renaud, J. F. (2010). De Finetti's
optimal dividends problem with an affine penalty function at ruin.
\textit{Insurance: Mathematics and Economics} \textbf{46},1, 98-108.



















\bibitem {Radner}Radner, R. Shepp, L. (1996) Risk vs.profit potential: a model
for corporate strategy. \textit{J. Econom. Dynamics Control} \textbf{20}, 1373--1393.



\bibitem {Schmidli book 2008}Schmidli, H. (2008). \textit{Stochastic Control
in Insurance}. Springer, New York.

\bibitem {shreve}Shreve, S.E., Lehoczky J.P. and Gaver, D.P. (1984) Optimal
consumption for general diffusions with absorbing and reflecting
barriers.\textit{ SIAM J. Control Optim.} \textbf{22}, 1, 55--75.










\end{thebibliography}
\end{document}